\documentclass[a4j,11pt,dvipdfmx]{article}

\usepackage{graphicx}
\usepackage{amsmath} 
\usepackage{amssymb}
\usepackage{amsthm}
\usepackage{mathrsfs}

\usepackage{color}
\usepackage{amsthm} 
\usepackage{subfigure} 
\usepackage{caption} 
\usepackage{bm} 
\usepackage{framed}
\usepackage{here}

\newtheoremstyle{mystyle}
{}
{}
{\normalfont}
{}
{\bfseries}
{}
{\newline}
{}
\theoremstyle{mystyle}

\newtheorem{Thm}{Theorem}
\newtheorem*{Thm*}{Theorem}
\newtheorem{Prop}{Proposition}
\newtheorem{lemm}{Lemma}
\newtheorem{Cor}{Corollary}

\newtheorem{Rem}{Remark}
\newtheorem*{Rem*}{Remark}
\newcommand{\1}{\mbox{1}\hspace{-0.25em}\mbox{l}}

\newcommand{\norm}[1]{\left\|#1\right\|}

\allowdisplaybreaks[4] 

\title{Asymptotics and practical aspects of testing normality with kernel methods}
\author{Natsumi Makigusa$^{1}$ and Kanta Naito$^{2}$\\
\\
Graduate School of Science and Engineering, Chiba University$^{1}$ \\ 
and \\
Graduate School of Science, Chiba University$^{2}$} 
\date{\today}
\begin{document}
\maketitle
\begin{abstract}
This paper is concerned with testing normality in a Hilbert space based on the maximum mean discrepancy.
Specifically, we discuss the behavior of the test from two standpoints: asymptotics and practical aspects.
Asymptotic normality of the test under a fixed alternative hypothesis is developed, which implies that the test has consistency.  
Asymptotic distribution of the test under a sequence of local alternatives is also derived, from which asymptotic null distribution of the test is obtained.
A concrete expression for the integral kernel associated with the null distribution is derived under the use of the Gaussian kernel, allowing the implementation of a reliable approximation of the null distribution.    
Simulations and applications to real data sets are reported with emphasis on high-dimension low-sample size cases. 
\end{abstract}
\section{Introduction}
For a probability distribution $P$, the test of the null hypothesis $H_0:P=N(m_0,\Sigma_0)$ against an alternative hypothesis $H_1: P \neq N(m_0,\Sigma_0)$ based on data $Y_1,\dots, Y_n \overset{i.i.d.}{\sim} P$ is known as testing normality, where $N(m_0,\Sigma_0)$ is the normal distribution with $m_0=\mathbb{E}[Y_1]$ and $\Sigma_0=V[Y_1]$.
Hereinafter we interchangeably use the terms normal distribution and Gaussian distribution.
Testing normality has  traditionally been an important problem in statistical science.
The literature on testing normality in a Euclidean space is huge, so we will not attempt a complete bibliography.
Notable papers include the Shapiro--Wilk test \cite{Shapiro and Wilk}, the Kolmogorov--Smirnov test \cite{Durbin}, the test based on skewness and kurtosis proposed by Mardia \cite{Mardia}, the omnibus test that uses a combination of skewness and kurtosis \cite{Mardia and Foster}, the test based on the empirical characteristic function \cite{Epps and Pulley}, a large comparative study by Romeu and Ozturk \cite{Romeu and Ozturk}.
Many others are given in the reference lists of these papers.

Testing normality has been generalized to a Hilbert space, and the approach with a reproducing kernel Hilbert space is known to be efficient for this problem.  
In paticular, an important application of such an approach is to provide methods for testing normality in a high-dimensional space.
For a given kernel function $k$, the Maximum Mean Discrepancy (MMD) between a distribution $P$ and the normal distribution $N(m_0,\Sigma_0)$ is defined as 
\[
\Delta=\left\|
\mathbb{E}_{Y \sim P}[k(\cdot,Y)] -\mathbb{E}_{Y \sim N(m_0,\Sigma_0)}[k(\cdot,Y)] 
\right\|_{H(k)}.
\]
The estimator of MMD has been proposed as a test statistic for $H_0$ in \cite{Kellner and Celisse}, where a bootstrap method was utilized to obtain the null distribution of the test.

However, \cite{Kellner and Celisse} did not comprehensively derive the asymptotic behavior of the test statistic.
In particular, asymptotic non-null behaviors of the MMD test were not investigated.
Furthermore, some practically important issues, such as a concrete expression of the test statistic, the integral equation associated with the asymptotic null distribution, the moments, and the approximation of the null distribution, were not sufficiently addressed, although a fast bootstrap for the approximation of the null distribution was proposed.  
Therefore, the aim of this research is to clarify the properties of the MMD test from two perspectives:
an asymptotic investigation under the situation $n \to \infty$, and a practical implementation the MMD test.
We first obtain the asymptotic distribution of the test by $\widehat{\Delta}^2$, a consistent estimator of $\Delta^2$, under a fixed alternative distribution $Q$. 
We also consider a sequence of local alternative distributions $P=P_n=(1-1/\sqrt{n})N(m_0,\Sigma_0)+(1/\sqrt{n})Q$, and address the asymptotic distribution of $\widehat{\Delta}^2$ under this sequence, by which the asymptotic null distribution of MMD test is derived.
For practical purposes, we consider the MMD test under the case where the utilized kernel $k$ is a Gaussian kernel and the Hilbert space $\mathcal{H}$ is $\mathbb{R}^d$.
Under this setting, we give the integral equation associated with the asymptotic distributions of the MMD test.
The first and second moments of the asymptotic null distribution can be obtained using the above expression of the integral equation, which yields an efficient and reliable approximation of the null distribution by exploiting a classical approximation method using a single weighted chi-squared distribution.

The rest of this paper is structured as follows.
Section \ref{s2} describes our framework and the testing normality. 
Section \ref{s3} develops asymptotics for the test by $\widehat{\Delta}^2$.
In particular,  the test by $\widehat{\Delta}^2$ under $H_1$ is addressed in Section \ref{s3.3}.
Furthermore, the behavior of $\widehat{\Delta}^2$ under a sequence of local alternative hypotheses is clarified in Section \ref{s3.4}, by which the asymptotic null distribution is developed.  
Section \ref{s4} examines the test statistic when kernel $k$ is a Gaussian kernel and the Hilbert space $\mathcal{H}$ is $\mathbb{R}^d$.
In particular, the form of the test statistic is shown in Section \ref{s4.1}, and the integral equation for obtaining the necessary eigenvalue for the asymptotic distribution is investigated in Section \ref{s4.2}.
Several approximations of the $\alpha$-quantile of the asymptotic null distribution of the test statistic are discussed in Section \ref{s4.3}, where we show that a single weighted chi-squared approximation works efficiently. 
Simulation results on the power of this test are reported in Section \ref{s4.4}.
Section \ref{s5} contains the results of applications to real data sets, including high-dimension low-sample size data.
Conclusions are given in Section \ref{s6}, and all proofs and calculations for theoretical results are provided in Section \ref{s7}.

\section{Setting and the MMD test}\label{s2}
Let $\mathcal{H}$ be a separable Hilbert space, and $(\mathcal{H},\mathcal{A})$ be a measurable space．
Let $Y_1,\dots,Y_n \in \mathcal{H}$ denote a sample of independent and identically distributed (i.i.d.) random variables drawn from an unknown distribution $P$.
Let the inner product of $\mathcal{H}$ be $\left<\cdot,\cdot\right>_{\mathcal{H}}$, and the associated norm $\|\cdot\|_{\mathcal{H}}=\sqrt{\left<\cdot,\cdot\right>_{\mathcal{H}}}$．
Our goal is to test whether $Y_i$ is a Gaussian random variable of $\mathcal{H}$ (see \cite{Mikusinski} for details).

Let us define the null hypothesis $H_0: P =$ Gaussian and the alternative hypothesis $H_1: P \neq$ Gaussian.
Following \cite{Gretton2007}, the gap between two distributions $P$ and $Q$ on $\mathcal{H}$ is measured by
\begin{equation}\label{E1}
\Delta(P,Q)=\sup_{f \in \mathcal{F}} \left|\mathbb{E}_{Y \sim P} [f(Y)] -\mathbb{E}_{Z \sim Q} [f(Z)]\right|,
\end{equation}
where $\mathcal{F}$ is a class of real-valued functions on $\mathcal{H}$.
Regardless of $\mathcal{F}$, $(\ref{E1})$ always defines a pseudo-metric on probability distributions.
In this paper, let $\mathcal{F}$ be the unit ball of a reproducing kernel Hilbert space $H(k)$ associated with a characteristic kernel $k:\mathcal{H} \times \mathcal{H} \to \mathbb{R}$ (see \cite{Aronszajn} and  \cite{Fukumizu} for details).
In addition, assume that $\mathbb{E}_{Y \sim P}[\sqrt{k(Y,Y)}] < \infty$ and $\mathbb{E}_{Y \sim Q}[\sqrt{k(Y,Y)}] < \infty$．
Then, $\Delta(P,Q)$ can be expressed as the gap between the Hilbert space embeddings of $P$ and $Q$ (see \cite{Gretton2012} for details):
\begin{equation*}
\Delta(P,Q)=\|\mu(P)-\mu(Q)\|_{H(k)}.
\end{equation*}
This distance between distributions $P$ and $Q$ is called the MMD．

In this paper, we aim to measure the gap between $P$ and $N(m_0,\Sigma_0)$ by
\[
\Delta^2=\left\|\mu(P)-\mu(N(m_0,\Sigma_0))\right\|_{H(k)}^2,
\]
where $m_0$ and $\Sigma_0$ are the mean and covariance operator of $P$, respectively (see \cite{Aronszajn} for details).
Here, $\Delta^2$ can be estimated by 
\[
\widehat{\Delta}^2=\left\|\frac{1}{n} \sum_{i=1}^{n}k(\cdot,Y_i) - \mu(N(\widehat{m},\widehat{\Sigma}))\right\|^2_{H(k)},
\]
where 
\begin{equation}\label{Em_hat}
\widehat{m}=\frac{1}{n}\sum_{i=1}^{n}Y_i
\end{equation}
and
\begin{equation}\label{ESigma_hat}
\widehat{\Sigma}=\frac{1}{n}\sum_{i=1}^{n}(Y_i-\widehat{m})^{\otimes 2}=\frac{1}{n}\sum_{i=1}^{n}\left<Y_i-\widehat{m},\cdot\right>_{\mathcal{H}}(Y_i-\widehat{m}).
\end{equation}

In this paper, we derive the asymptotic distribution of $\widehat{\Delta}^2$ under $H_1$ and a sequence of local alternative hypotheses, both of which are not derived in the previous study \cite{Kellner and Celisse}.

\section{Asymptotics}\label{s3}
	In this section, we develop the asymptotic distributions of the test by $\widehat{\Delta}^2$ under an alternative hypothesis $H_1$ and a sequence of local alternatives.
	Further, the asymptotic null distribution is obtained as by-product. 
\subsection{Preliminaries}
	Assume that $\mu(N(\cdot,\cdot))$ is twice continuously Fr\'echet differentiable ($C^2$-class, see Section 1 in \cite{Abraham and Robbin} for details) in 
	\[
	B((m_0,\Sigma_0),\varepsilon)=\left\{x \in \mathcal{H} \oplus HS(\mathcal{H}) ~|~\left\|(m_0,\Sigma_0) -x\right\|_{\mathcal{H} \oplus HS(\mathcal{H})}< \varepsilon \right\},
	\]
	where $\mathcal{H} \oplus HS(\mathcal{H})$ is a direct sum space of $\mathcal{H}$ and $HS(\mathcal{H})$．
	Let $D_{(m_0,\Sigma_0)} \mu(N(\cdot,\cdot))$ be the Fr\'echet derivative of $\mu(N(\cdot,\cdot))$ at a point $(m_0,\Sigma_0)$ (see Section 1 in \cite{Abraham and Robbin} for details), and let
	\begin{equation}\label{E3}
	f(x)=k(\cdot,x) -D_{(m_0,\Sigma_0)}\mu(N(x-m_0,(x-m_0)^{\otimes 2} -\Sigma_0))-\mu(N(m_0,\Sigma_0)).
	\end{equation}
	First, we note that $\widehat{\Sigma}$ can be alternatively expressed as
	\begin{align*}
	\frac{1}{n}\sum_{i=1}^{n}(Y_i-\widehat{m})^{\otimes 2} 
	=\frac{1}{n}\sum_{i=1}^{n}(Y_i-m_0)^{\otimes 2} -(\widehat{m}-m_0)^{\otimes 2}.
	\end{align*}
	From Taylor's theorem addressed in Section 1.5 in \cite{Abraham and Robbin}, 
	\begin{align}
	&\sqrt{n}\mu(N(\widehat{m},\widehat{\Sigma}))\nonumber\\
	&=\sqrt{n}\left(\mu(N(m_0,\Sigma_0))+D_{(m_0,\Sigma_0)}\mu(N(\widehat{m}-m_0,\widehat{\Sigma}-\Sigma_0)\right.\nonumber\\
	&\left.~~~~~~~~~~~~~~~~~~~~~~~~~~~~~~~~~~~~~~~~~~~~~~+R_2((m_0,\Sigma_0),(\widehat{m}-m_0,\widehat{\Sigma}-\Sigma_0))\right)\nonumber\\
	&=\sqrt{n}\mu(N(m_0,\Sigma_0))
	+\sqrt{n}D_{(m_0,\Sigma_0)} \mu\left(N\left(\frac{1}{n}\sum_{i=1}^{n}(Y_i-m_0),\frac{1}{n}\sum_{i=1}^{n}((Y_i-m_0)^{\otimes 2} -\Sigma_0)\right)\right)\nonumber\\
	&~~~~~+D_{(m_0,\Sigma_0)}\mu((0,\sqrt{n}(\widehat{m}-m_0)^{\otimes 2})) +\sqrt{n}R_2((m_0,\Sigma_0),(\widehat{m}-m_0,\widehat{\Sigma}-\Sigma_0)), \label{E5}
	\end{align}
	where $\displaystyle R_2((m_0,\Sigma_0),(\widehat{m}-m_0,\widehat{\Sigma}-\Sigma_0))=\int^{1}_{0} (1-s) D^2_{(m_0,\Sigma_0)} \mu(N(\widehat{m}-m_0,\widehat{\Sigma}-\Sigma_0)^2) ds$.
	The self-adjoint Hilbert--Schmidt operator $S_k$ is defined as 
	\begin{equation}\label{E6}
	S_k : L^2(\mathcal{H},N(m_0,\Sigma_0)) \to L^2(\mathcal{H},N(m_0,\Sigma_0)),~~g \mapsto \int_{\mathcal{H}} \left<f(\cdot),f(y)\right>_{H(\overline{k})}g(y)dN(m_0,\Sigma_0)(y)
	\end{equation}
	(see Sections VI.1, VI.3 and VI.6 in \cite{Reed and Simon} for details). 
	The eigenvalue $\lambda_{\ell}$ of $S_k$ satisfies 
	\begin{equation}\label{E30}
	\int_{\mathcal{H}} \left<f(x),f(y)\right>_{H(k)} \Psi_\ell(y) dN(m_0,\Sigma_0) (y) =\lambda_\ell\Psi_\ell(x), 
	\end{equation}
	where $\Psi_{\ell}$ is the eigenfunction corresponding to $\lambda_{\ell}$.
	These $\Psi_{\ell}$ form an orthonormal system in $L^2(\mathcal{H},N(m_0,\Sigma_0))$ as follows:
	\begin{equation}
	\int_{\mathcal{H}} \Psi_i(y) \Psi_j(y) dN(m_0,\Sigma_0) (y) =\delta_{ij}.\label{E31}
	\end{equation}
	For a distribution $Q$ on $\mathcal{H}$, let
	\[
	\eta(Q) =\mu(Q) -\mu(N(m_0,\Sigma_0)),
	\]
	and let
	\[
	\eta_{\ell} (Q) =\int_{\mathcal{H}} \left<\eta(Q),f(y)\right>_{H(k)} \Psi_{\ell} (y) dN(m_0,\Sigma_0) (y),\ \ell=1,2,\ldots.
	\]
	Consider a class of distributions of the alternative hypotheses for which element $Q$ satisfies 
	\begin{equation}\label{E15}
	\mathbb{E}_{X \sim Q}\left[k(X,X)\right]< \infty,
	\end{equation}
	\begin{equation}\label{E16}
	\mathbb{E}_{X \sim Q}\left[\left\|(X-m_0,(X-m_0)^{\otimes 2}-\Sigma_0)\right\|^2_{\mathcal{H} \oplus HS(\mathcal{H})}\right]< \infty
	\end{equation}
	and
	\begin{equation}\label{E17}
	\sum_{\ell=1}^{\infty}
	\frac{\eta^2_{\ell}(Q)}{\lambda_{\ell}}<\infty,
	\end{equation}
	where $\lambda_{\ell}$ and $\Psi_{\ell}(y)$ are  the eigenvalue and eigenfunction of $S_k$ in (\ref{E6}).
	In other words, let 
	\[
	\mathcal{A}_k=\{Q~|~
	Q \neq \text{Gaussian},~~\mathbb{E}_{X \sim Q}[X]=m_0,~~V_{X \sim Q}[X]=\Sigma_0,~~(\ref{E15}),~(\ref{E16}) \text{ and } (\ref{E17})\text{ are held}
	\}
	\]
	be the set of distributions of alternative hypotheses.

\subsection{Asymptotic nonnull distribution}\label{s3.3}
In this section, we investigate the asymptotic distribution of $\widehat{\Delta}^2$ under $H_1$ and prove the consistency of this test.
In what follows, the symbol $``\overset{\mathcal{D}}{\to}"$ designates convergence in distribution.
\begin{Thm}\label{T2}
	Let $Q \in \mathcal{A}_k$.
	Then under $H_1:P = Q$, as $n \to \infty$,
	\[
	\sqrt{n}(\widehat{\Delta}^2-\Delta^2)  \xrightarrow{\mathcal{D}} N\left(0,v^2\right),
	\]
	where $v^2=4\left<V[f(Y_1)](\mu(P)-\mu(N(m_0,\Sigma_0))),\mu(P)-\mu(N(m_0,\Sigma_0))\right>_{H(k)}$.
\end{Thm}
\noindent The asymptotic variance $v^2$ in Theorem \ref{T2} is guaranteed to exist by the definition of $\mathcal{A}_k$. 
See Remark \ref{Rem3} in Section \ref{s7} for details.
\begin{Rem}\label{Rem1}
	We see by Theorem \ref{T2} that
	\[
	\frac{\sqrt{n}(\widehat{\Delta}^2-\Delta^2)}{v}  \xrightarrow{\mathcal{D}} N(0,1).
	\]
	Thus, we can evaluate the power of the test by $n\widehat{\Delta}^2$ as
	\begin{align*}
	\mathbb{P}(n\widehat{\Delta}^2 \geq t_{\alpha}~|~H_1) 
	&=\mathbb{P}(n(\widehat{\Delta}^2-\Delta^2) \geq  t_{\alpha} -n\Delta^2~|~H_1)\\
	&=\mathbb{P}\left(\left.
	\frac{\sqrt{n}(\widehat{\Delta}^2-\Delta^2)}{v} \geq \frac{t_{\alpha}}{\sqrt{n}v}-\frac{\sqrt{n}\Delta^2}{v} ~\right|~H_1
	\right)\\
	& \approx 1-\Phi\left(\frac{t_{\alpha}}{\sqrt{n}v}-\frac{\sqrt{n}\Delta^2}{v}\right)\\
	& \to 1
	\end{align*}
	as $n \to \infty$, where $t_{\alpha}$ is ($1-\alpha$)-quantile of the $n\widehat{\Delta}^2$ distribution under $H_0$, and $\Phi$ is the distribution function of the standard Gaussian.
	Therefore, this test is consistent.
\end{Rem}

\subsection{Asymptotic distribution under contiguous alternatives}\label{s3.4}
In this section, we develop the asymptotic distribution of $\widehat{\Delta}^2$ under a sequence of local alternative distributions $P=P_n(Q)=(1-1/\sqrt{n}) N(m_0,\Sigma_0) +(1/\sqrt{n})Q$, for $Q \in \mathcal{A}_k$.
The proof is based on the asymptotic result of degenerate V-statistics (see Section 5.5 in \cite{Serfling}). 
Further, we derive the asymptotic distribution of $\widehat{\Delta}^2$ under $H_0$ using Theorem \ref{T3}.
\begin{Thm}\label{T3}
	Assume $\left<f(x),f(y)\right>_{H(k)} \in L^2(\mathcal{H} \times \mathcal{H},N(m_0,\Sigma_0) \times N(m_0,\Sigma_0))$.
	Then, under a local alternative hypothesis $P=P_n(Q)$, as $n \to \infty$,
	\[
	n\widehat{\Delta}^2 \xrightarrow{\mathcal{D}} W\equiv \sum_{\ell=1}^{\infty} \lambda_{\ell} W^2_{\ell},
	\]
	where $W_{\ell}~(\ell=1,2,\dots)$ are independent and $\displaystyle W_\ell \sim N\left(\frac{\eta_\ell(Q)}{\lambda_{\ell}},1\right)$ and $\lambda_{\ell}$ is the eigenvalue of $S_k$ in (\ref{E6}). 
\end{Thm}
Using Theorem \ref{T3}, the asymptotic null distribution can be obtained as follows: 
\begin{Cor}[Asymptotic null distribution]\label{C1}
	Assume $\left<f(x),f(y)\right>_{H(k)} \in L^2(\mathcal{H} \times \mathcal{H},N(m_0,\Sigma_0) \times N(m_0,\Sigma_0))$.
	Then under $H_0:P=$ Gaussian, as $n \to \infty$,
	\[
	n\widehat{\Delta}^2 \xrightarrow{\mathcal{D}}Z\equiv\sum_{\ell=1}^{\infty} \lambda_{\ell} Z^2_{\ell},
	\]	
	where $Z_{\ell} \overset{i.i.d.}{\sim} N(0,1)$.
\end{Cor}
Theorem 2 and Corollary 1 reveal that the local power of the test by $n\widehat{\Delta}^2$ is essentially dominated by the noncentrality parameters.
Now, it follows that 
\[
\eta_{\ell}(Q)=\int_{\mathcal{H}}\left<\mathbb{E}_{X \sim Q}[f(X)],f(y)\right>_{H(k)} \Psi_{\ell}(y)dN(m_0,\Sigma_{0})(y)
=\lambda_{\ell} \mathbb{E}_{X \sim Q}[\Psi_{\ell}(X)],
\]
by which we obtain 
\[
\mathbb{E}[W]=\sum_{\ell=1}^{\infty}\lambda_{\ell}\left(
1+\frac{\eta_{\ell}(Q)^2}{\lambda_{\ell}^2}
\right)=\sum_{\ell=1}^{\infty}\lambda_{\ell}(1+\mathbb{E}_{X \sim Q}[\Psi_{\ell}]^2).
\]
Hence, the local power results in the  magnitude of $\mathbb{E}_{X \sim Q}[\Psi_{\ell}(X)]$.

The index of performance of the test by $n\widehat{\Delta}^2$ discussed in \cite{Gregory} becomes 
\[
e(n\widehat{\Delta}^2,Q)=\frac{\mathbb{E}[W]}{\sqrt{V[Z]}}=\frac{\sum_{\ell=1}^{\infty}\lambda_{\ell}(1+\mathbb{E}_{X \sim Q}[\Psi_{\ell}(X)]^2)}{\sqrt{2\sum_{\ell=1}^{\infty} \lambda_{\ell}^2}}.
\]
Theoretical comparison with another test with the same asymptotic null distribution can be demonstrated using the relative efficiency calculated by the ratio of the above indices.
The relative efficiency is the limiting ratio of sample sizes needed to give the same asymptotic local power (see \cite{Gregory} and \cite{Serfling}).

\section{Practical aspects}\label{s4}
In this section, we discuss $\widehat{\Delta}^2$ and $f(x)$ of $(\ref{E3})$ when $\mathcal{H}=\mathbb{R}^d$ and $k(\cdot,\cdot)$ is the Gaussian kernel: 
\begin{equation}\label{Gaussian kernel}
k(\underline{t},\underline{s})=\exp\left(-\sigma\|\underline{t}-\underline{s}\|^2_{\mathbb{R}^d}\right),~\sigma>0.
\end{equation}
Let $\underline{Y}_1,\dots,\underline{Y}_n \in \mathbb{R}^d$ denote a sample of i.i.d. random variables drawn from an unknown distribution $P$.
The estimator $\widehat{\underline{m}}$ is (\ref{Em_hat}) and $\widehat{\Sigma}$ is (\ref{ESigma_hat}) with $\widehat{\Sigma}=(1/n)\sum_{i=1}^{n}(\underline{Y}_i-\widehat{\underline{m}})(\underline{Y}_i-\widehat{\underline{m}})^T$, where $T$ denotes transpose.

Further, we derive the integral equation associated with (\ref{E30}).
It follows from  Proposition 4.2 of \cite{Kellner and Celisse} that
\begin{align*}
\mu(N(\underline{m},\Sigma))(\cdot)
=|I_d+2\sigma\Sigma|^{-1/2} \exp\left(-\sigma(\cdot-\underline{m})^T(I_d+2\sigma\Sigma)^{-1} (\cdot-\underline{m})\right),
\end{align*}
where $\underline{m}=
\begin{bmatrix}
m_1&  \cdots& m_d
\end{bmatrix}^T  \in \mathbb{R}^d
$
,
\[
\Sigma=
\begin{bmatrix}
\sigma_{11}&\cdots  &\sigma_{1d}  \\ 
\vdots& \ddots  & \vdots \\ 
\sigma_{1d}& \cdots & \sigma_{dd}
\end{bmatrix}
\]
and $I_d$ is the $d \times d$ identity matrix.

\subsection{Practical form of the test statistic}\label{s4.1}
We see from Proposition 4.2 of \cite{Kellner and Celisse} that
\[
\|\mu(N(\underline{m},\Sigma))\|^2_{H(k)}=|I_d+4\sigma\Sigma|^{-1/2},
\]
from which the test statistic $\widehat{\Delta}^2$ is
\begin{align*}
\widehat{\Delta}^2
&=\left\|\frac{1}{n}\sum_{i=1}^{n}k(\cdot,\underline{Y}_i) -\mu(N(\widehat{\underline{m}},\widehat{\Sigma}))(\cdot)\right\|^2_{H(k)}\\
&=\frac{1}{n^2} \sum_{i,j=1}^{n} \exp\left(-\sigma\|\underline{Y}_i-\underline{Y}_j\|^2_{\mathbb{R}^d}\right)\\
&~~~~~-|I_d+2\sigma\widehat{\Sigma}|^{-1/2} \frac{2}{n}\sum_{i=1}^{n}\exp\left(-\sigma(\underline{Y}_i-\widehat{\underline{m}})^T(I_d+2\sigma\widehat{\Sigma})^{-1} (\underline{Y}_i-\widehat{\underline{m}})\right)+|I_d+4\sigma\widehat{\Sigma}|^{-1/2}.
\end{align*}

\subsection{Integral equation}\label{s4.2}
We need to obtain $\left<f(\underline{x})(\cdot),f(\underline{y})(\cdot)\right>_{H(k)}$ in the integral equation (\ref{E30}).
To do this, we aim to calculate (\ref{E3}) for the Gaussian kernel $k(\cdot,\cdot)$ to obtain the derivative of $\mu(N(\underline{m},\Sigma))(\cdot)$ with $
\begin{bmatrix}
\underline{m}& \text{vech}(\Sigma) 
\end{bmatrix}^T $.
Here the symbol vech is defined as
\[
\text{vech}(A)=
\begin{bmatrix}
a_{11}& \cdots & a_{dd} & a_{12} & \cdots & a_{1d} &a_{23}&\cdots &a_{d-1,d}
\end{bmatrix}^T 
\]
for a $d \times d$ symmetric matrix $A=(a_{ij})_{1 \leq i,j \leq d}.$
The derivative of $\mu(N(\underline{m},\Sigma))$ with $\underline{m}$ is easily obtained as
\begin{equation}\label{E32}
\frac{\partial}{\partial \underline{m}} \mu(N(\underline{m},\Sigma))(\cdot)=2\sigma\mu(N(\underline{m},\Sigma))(\cdot)(I_d+2\sigma\Sigma)^{-1}(\cdot-\underline{m}).
\end{equation}
Straightforward calculations yield the derivative of $\mu(N(\underline{m},\Sigma))(\cdot)$ with $\sigma_{ij}$ as
\begin{align}\label{E33}
&\frac{\partial}{\partial \sigma_{ij}}\mu(N(\underline{m},\Sigma))(\cdot)\nonumber\\
&=\sigma \mu(N(\underline{m},\Sigma))(\cdot) \text{tr}\left(
\frac{\partial \Sigma}{\partial \sigma_{ij}}  (I_d+2\sigma\Sigma)^{-1} \left(
2\sigma (\cdot-\underline{m})(\cdot-\underline{m})^T  (I_d+2\sigma\Sigma)^{-1}-I_d
\right)
\right).
\end{align}
Therefore, we have
\begin{align*}
&\frac{\partial}{\partial \text{vech}(\Sigma)} \mu(N(\underline{m},\Sigma))(\cdot)\\
&=2\sigma \mu(N(\underline{m},\Sigma))(\cdot)\text{vech}\left\{
 (I_d+2\sigma\Sigma)^{-1} (2\sigma(\cdot-\underline{m}) (\cdot-\underline{m})^T(I_d+2\sigma\Sigma)^{-1} -I_d)
\right\}\\
&~~~~~-\sigma \mu(N(\underline{m},\Sigma))(\cdot) \text{vech}\left\{
\text{diag} \left\{
(I_d+2\sigma\Sigma)^{-1} (2\sigma(\cdot-\underline{m}) (\cdot-\underline{m})^T(I_d+2\sigma\Sigma)^{-1} -I_d)
\right\}
\right\}.
\end{align*}
Let
\begin{align*}
A(\underline{m},\Sigma)(\cdot)
&=\text{vech}\left\{2
(I_d+2\sigma\Sigma)^{-1} (2\sigma(\cdot-\underline{m}) (\cdot-\underline{m})^T(I_d+2\sigma\Sigma)^{-1} -I_d)
\right.\\
&~~~~~\left.-
\text{diag} \left\{
(I_d+2\sigma\Sigma)^{-1} (2\sigma(\cdot-\underline{m}) (\cdot-\underline{m})^T(I_d+2\sigma\Sigma)^{-1} -I_d)
\right\}
\right\}
\end{align*}
and let
\begin{eqnarray}\label{DEF_MT_V}
V=I_d+2\sigma\Sigma_0.
\end{eqnarray}
Note that $V$ is positive definite.
Then, we get an expansion of $\mu(N(\widehat{\underline{m}},\widehat{\Sigma}))(\cdot)$
\begin{align*}
&\mu(N(\widehat{\underline{m}},\widehat{\Sigma}))(\cdot)\\
&=\mu(N(\underline{m}_0,\Sigma_0))(\cdot)+D_{(\underline{m}_0,\Sigma_0)}\mu(N(\widehat{\underline{m}}-\underline{m}_0,\widehat{\Sigma}-\Sigma_0))+O_p\left(\frac{1}{\sqrt{n}}\right)\\
&=\mu(N(\underline{m}_0,\Sigma_0))(\cdot)+2\sigma\mu(N(\underline{m}_0,\Sigma_0))(\cdot)(\cdot-\underline{m}_0)^TV^{-1}(\widehat{\underline{m}}-\underline{m}_0)\\
&~~~~~+\sigma \mu(N(\underline{m}_0,\Sigma_0))(\cdot
)A(\underline{m}_0,\Sigma_0)(\cdot)^T\text{vech}(\widehat{\Sigma}-\Sigma_0)+O_p\left(\frac{1}{\sqrt{n}}\right)\\
&=\frac{1}{n}\sum_{i=1}^{n}\mu(N(\underline{m}_0,\Sigma_0))(\cdot) \left\{1+2\sigma(\cdot-\underline{m}_0)^TV^{-1}(\underline{Y}_i-\underline{m}_0)\right.\\
&~~~~~+\left.\sigma A(\underline{m}_0,\Sigma_0)(\cdot)^T\text{vech}\left((\underline{Y}_i-\underline{m}_0)(\underline{Y}_i-\underline{m}_0)^T-\Sigma_0\right)\right\}+O_p\left(\frac{1}{\sqrt{n}}\right).
\end{align*}
Therefore, we have a practical form of (\ref{E3}) under the case $\mathcal{H}=\mathbb{R}^d$ and $k(\underline{t},\underline{s})=\exp(-\sigma\left\|\underline{t}-\underline{s}\right\|^2_{\mathbb{R}^d})$ as
\begin{align}
f(\underline{x})(\cdot)
&=k(\cdot,\underline{x}) -\mu(N(\underline{m}_0,\Sigma_0))(\cdot)\left\{
1+2\sigma(\cdot-\underline{m}_0)^TV^{-1} (\underline{x}-\underline{m}_0)\right.\nonumber\\
&~~~~~+\left.\sigma A(\underline{m}_0,\Sigma_0)(\cdot)^T \text{vech}((\underline{x}-\underline{m}_0)(\underline{x}-\underline{m}_0)^T-\Sigma_0)
\right\}. \label{E40}
\end{align}
We aim to calculate $\left<f(\underline{x})(\cdot), f(\underline{y})(\cdot)\right>_{H(k)}$ in the integral equation (\ref{E30}) using the fact that the inner product of the reproducing kernel Hilbert space corresponding to the Gaussian kernel is given as
\[
\left<g,h\right>_{H(k)} =\sqrt{(4\pi\sigma)^d} \int_{\mathbb{R}^d} \widehat{g}(\underline{t}) \overline{\widehat{h}(\underline{t})} \exp\left(\frac{1}{4\sigma} \underline{t}^T\underline{t}\right) d\underline{t},~~\text{for}~h,g \in H(k),
\]
where $\widehat{g}(\underline{t})$ is the Fourier transform of $g$ (see Theorem 10.12 of \cite{Wendland}).
First, $\widehat{f(\underline{x})}(\underline{t})$ is calculated  using the characteristic function of the normal distribution to obtain 
\begin{align}\label{E34}
\widehat{f(\underline{x})}(\underline{t}) &=\frac{1}{\sqrt{(4\pi\sigma)^d}} \left\{\exp\left(-i\underline{x}^T\underline{t}-\frac{1}{4\sigma}\underline{t}^T\underline{t}\right)
-\exp\left(-i\underline{m}_0^T\underline{t}-\frac{1}{4\sigma}\underline{t}^TV\underline{t}\right)\right.\nonumber\\
&~~~~~\left.
\times \left(
1-i\underline{t}^T(\underline{x}-\underline{m}_0) -\frac{1}{2} \underline{t}^TB(\underline{x})\underline{t}
\right)
\right\},
\end{align}
where
\begin{eqnarray}
B(\underline{x})=(\underline{x}-\underline{m}_0)(\underline{x}-\underline{m}_0)^T-\Sigma_0. \label{DEF_Bx}
\end{eqnarray}
Using (\ref{E34}), we have 
\begin{align}
&\left<f(\underline{x})(\cdot),f(\underline{y}) (\cdot)\right>_{H(k)}\nonumber\\
&=\sqrt{(4\pi\sigma)^d}\int_{\mathbb{R}^d}\widehat{f}(\underline{x})(\underline{t})\overline{\widehat{f}(\underline{y})(\underline{t})} \exp\left(\frac{1}{4\sigma}\underline{t}^T\underline{t}\right)d\underline{t}\nonumber\\
&=\exp\left(-\sigma \|\underline{x}-\underline{y}\|^2_{\mathbb{R}^d}\right)\nonumber\\
&~~~~~-|V|^{-1/2} \exp\left(-\sigma(\underline{x}-\underline{m}_0)^TV^{-1}(\underline{x}-\underline{m}_0) \right)
\Bigl\{
1+2\sigma(\underline{x}-\underline{m}_0)^TV^{-1} (\underline{y}-\underline{m}_0)\nonumber \\
&~~~~~+\sigma \text{tr}\left[
V^{-1}\{2\sigma(B(\underline{x}) +\Sigma_0) V^{-1} -I_d \}B(\underline{y})\right] \Bigr\}\nonumber\\
&~~~~~-|V|^{-1/2} \exp\left(-\sigma(\underline{y}-\underline{m}_0)^TV^{-1}(\underline{y}-\underline{m}_0) \right)
\Bigl\{
1+2\sigma(\underline{y}-\underline{m}_0)^TV^{-1} (\underline{x}-\underline{m}_0)\nonumber \\
&~~~~~+\sigma \text{tr}\left[
V^{-1}\{2\sigma(B(\underline{y}) +\Sigma_0) V^{-1} -I_d \}B(\underline{x})\right] \Bigr\}\nonumber\\
&~~~~~+|2V- I_d|^{-1/2} \Bigl\{
1+\sigma \text{tr}\left[
(2(\underline{x}-\underline{m}_0)(\underline{y}-\underline{m}_0)^T-B(\underline{x})-B(\underline{y})) 
(2V- I_d)^{-1}
\right]
\nonumber\\
&~~~~~
+\sigma^2\Bigl(\text{tr}\left[B(\underline{x}) (2V- I_d)^{-1}\right] \text{tr}\left[B(\underline{y}) (2V - I_d)^{-1}\right]\nonumber\\
&~~~~~+2\text{tr}\left[B(\underline{x}) (2V- I_d)^{-1} B(\underline{y}) (2V- I_d)^{-1}\right]
\Bigr)
\Bigr\}.\label{E47}
\end{align}
Note that $2V-I_d=I_d+4\sigma\Sigma_0$ is invertible by its positive definiteness.
We would solve the integral equation (\ref{E30}) with (\ref{E47}) to clarify all eigenvalues involved in the asymptotic null distribution developed in  Corollary \ref{C1}, but it is not easy to solve the equation at this stage.
Hence, we defer the problem as a future project.

\subsection{Approximation of the distribution}\label{s4.3}

In this section, we discuss methods to approximate the null distribution of the MMD test.
The asymptotic null distribution of the MMD test was obtained in Corollary \ref{C1} as an infinite sum of weighted chi-squared random variables with one degree of freedom. 
For the case using a Gaussian kernel, we derived the integral kernel with eigenvalues that are the weights appearing in the asymptotic distribution.
Each eigenvalue of this integral kernel is hard to obtain at this stage, but the sum of the eigenvalues and the sum of the squared eigenvalues can be derived easily.
In fact, we see from the general theory of Hilbert spaces that the asymptotic mean and variance are, respectively, obtained as
\begin{eqnarray}
E[Z]&=&\sum_{i=1}^{\infty}\lambda_{i}=\int_{{\mathbb R}^d}\left<f(\underline{x}),f(\underline{x})\right>_{H(k)} dN(m_{0},\Sigma_{0})(\underline{x}), \label{E_Z}\\
V[Z]&=&2\sum_{i=1}^{\infty}\lambda_{i}^{2}=2\int_{{\mathbb R}^{d}}\int_{{\mathbb R}^{d}}\left<f(\underline{x}),f(\underline{y})\right>^{2}_{H(k)}dN(m_{0},\Sigma_{0})(\underline{x})dN(m_{0},\Sigma_{0})(\underline{y}). \label{V_Z}
\end{eqnarray}
Our approximation of the null distribution is based on a classical method using a single weighted chi-squared distribution, as discussed in \cite{Hirotsu1}, \cite{Hirotsu2}.
The method aims to approximate the distribution of the sum of weighted chi-square random variables by using a single weighted chi-squared random variable of the form $c\chi_{r}^{2}$.
Suppose we have appropriate estimates $\widehat{E[Z]}$ and $\widehat{V[Z]}$ of (\ref{E_Z}) and (\ref{V_Z}), respectively.
The method is to fit the first two cumulants of the $c\chi_{r}^{2}$ distribution to those of $Z$, which implies that 
\begin{eqnarray}
c=\frac{\widehat{V[Z]}}{2\widehat{E[Z]}}\ \ \ \mbox{and}\ \ \ r=\frac{2\widehat{E[Z]}^{2}}{\widehat{V[Z]}}. \label{cr_formula}
\end{eqnarray}
It is known that there is no mathematical validity to this approximation, since there is no convergence result: it is heuristic in this sense.
However, this approximation sometimes works very well as reported in \cite{Hirotsu1}, \cite{Hirotsu2}.
The same basic approximation by a two-parameter gamma distribution was also discussed in \cite{Gretton2009}, where computational efficiency of the method was emphasized.
Therefore, it is worth checking if this approximation works for our MMD test, especially in high-dimension low-sample size cases.   
In the following sections, we propose two methods to obtain $\widehat{E[Z]}$ and $\widehat{V[Z]}$.

\subsubsection{Single Weighted Chi-Squared Approximation I}\label{s431}
The expression derived in (\ref{E47}) motivates us to obtain concrete formulas for $E[Z]$ and $V[Z]$ in (\ref{E_Z}) and (\ref{V_Z}), respectively.
First, we consider (\ref{E_Z}).
Straightforward calculations using properties of Gaussian density as well as formulas of expectation of quadratic forms yield the following result.

\begin{Prop}\label{P1}
\begin{eqnarray}
E[Z]&=&1-|2V-I_d|^{-1/2}\Big\{1+2\sigma \text{tr}[(V+2\sigma \Sigma_0)^{-1}\Sigma_0]  \nonumber \\
\ &\ &\hspace{1cm}+2\sigma^2
\{\text{tr}[(V+2\sigma \Sigma_0)^{-1}\Sigma_0]\}^2+4\sigma^2\text{tr}[\{(V+2\sigma \Sigma_0)^{-1} \Sigma_0\}^2]
\Big\}, \label{E_Z_formula}
\end{eqnarray}
where $V$ is the matrix given in (\ref{DEF_MT_V}).
\end{Prop}

The proof of Proposition \ref{P1} is in Section 7.

Next, we consider $V[Z]$.
We see from direct computations that (\ref{V_Z}) essentially consists of the expectations of the products of quadratic forms with Gaussian random vectors:

\begin{eqnarray}
Q_{1}(A)&=&E[\underline{X}^{T}A\underline{X}], \label{Q_1} \\
Q_{2}(A,B)&=&E[\underline{X}^{T}A\underline{X}\cdot\underline{X}^{T}B\underline{X} ], \label{Q_2} \\
Q_{3}(A,B,C)&=&E[\underline{X}^{T}A\underline{X}\cdot\underline{X}^{T}B\underline{X}\cdot\underline{X}^{T}C\underline{X}], \label{Q_3} \\
Q_{4}(A,B,C,D)&=&E[\underline{X}^{T}A\underline{X}\cdot\underline{X}^{T}B\underline{X}\cdot\underline{X}^{T}C\underline{X}\cdot\underline{X}^{T}D\underline{X}], \label{Q_4} 
\end{eqnarray} 
where $A$, $B$, $C$, and $D$ are $2d \times 2d$ matrices, and $\underline{X}$ is distributed as $N_{2d}(\underline{0},I_{2d})$.
Note that the expressions of (\ref{Q_1}), (\ref{Q_2}), (\ref{Q_3}) and (\ref{Q_4}) are all given in Section 9.6 in \cite{Schott} as functions of $A$, $B$, $C$ and $D$, so we omit those expressions here.

We see by tedious but straightforward calculations with (\ref{Q_1}), (\ref{Q_2}), (\ref{Q_3}) and (\ref{Q_4}) arranged in Section 7 that $V[Z]$ is finally expressed as follows.

 \begin{Prop}\label{P2}
 	\begin{align}
 &V[Z] \nonumber \\
 &=2|I_d+8\sigma \Sigma_0|^{-1/2} \nonumber \\
 &~~~~-4|V|^{-1/2} |V+4\sigma\Sigma_0|^{-1/2}\Big\{
  1
  +\frac{1}{2}\sigma^2\{\text{tr}[V^{-1} \Sigma_0]\}^2
  +\sigma^2 \text{tr}[\{V^{-1} \Sigma_0\}^2]  \nonumber \\
 &~~~~+\frac{1}{2}\sigma^2\{\text{tr}[(V+4\sigma \Sigma_0)^{-1} \Sigma_0]\}^2+\sigma^2 \text{tr}[\{(V+4\sigma \Sigma_0)^{-1} \Sigma_0\}^2]
  +\sigma\text{tr}[V^{-1} \Sigma_0]  \nonumber \\
 &~~~~-\sigma \text{tr}[(V+4\sigma \Sigma_0)^{-1} \Sigma_0]
  -\sigma^2 \text{tr}[V^{-1} \Sigma_0] \text{tr}[(V+4\sigma \Sigma_0)^{-1} \Sigma_0]
  \Big\} \nonumber \\
 &~~~~+2|V+2\sigma \Sigma_0|^{-1} \Big\{
  1+8\sigma^2 \text{tr}[\{(V+2\sigma \Sigma_0)^{-1} \Sigma_0\}^2]
  +12\sigma^4 \{\text{tr}[\{(V+2\sigma \Sigma_0)^{-1} \Sigma_0\}^2]\}^2 \nonumber \\
  &~~~~~
  +24\sigma^4 \text{tr}[\{(V+2\sigma \Sigma_0)^{-1} \Sigma_0\}^4]
  \Big\}		, \label{V_Z_formula}
 	\end{align}
 	where $V$ is the matrix given in (\ref{DEF_MT_V}).
 \end{Prop}


It is worth noting that $E[Z]$ in Proposition 1 and $V[Z]$ in Proposition 2 both depend only on $\Sigma_{0}$ and the scale parameter $\sigma$ of the kernel.
Therefore, once we have an estimate $\widehat{\Sigma}$ of $\Sigma_{0}$, substituting $\widehat{\Sigma}$ into the expressions in Propositions 1 and 2 provides $\widehat{E[Z]}$ and $\widehat{V[Z]}$ for a given $\sigma$.

The null distribution of $n\widehat{\Delta}^{2}$ is approximated by $\widehat{c}\chi_{\widehat{r}}^{2}$, where $\widehat{c}$ and $\widehat{r}$ can be obtained via (\ref{cr_formula}) with $\widehat{E[Z]}$ and $\widehat{V[Z]}$ calculated by substituting $\widehat{\Sigma}$ into $\Sigma_{0}$ in (\ref{E_Z_formula}) and (\ref{V_Z_formula}).
Although the formulas (\ref{E_Z_formula}) and $(\ref{V_Z_formula})$ look very long, the approximation to the percentile based on this method can be calculated quickly and performs fairly well as seen in Section \ref{s433}.

\subsubsection{Single Weighted Chi-Squared Approximation II}\label{s432}
Since the eigenvalue $\lambda_{\ell}$ satisfies (\ref{E30}), we can consider an alternative approach to estimating $E[Z]$ and $V[Z]$ based on equation (\ref{E30}).

Once we obtain $\widehat{m}$ and $\widehat{\Sigma}$ from the given dataset, we generate $\underline{X}_1,\dots,\underline{X}_L \overset{i.i.d.}{\sim} N(\widehat{m},\widehat{\Sigma})$ to compose the Gram matrix
\[
\mathcal{G}_{L}=\frac{1}{L}
\begin{bmatrix}
\left<f(\underline{X}_1),f(\underline{X}_1)\right>_{H(k)}& \cdots & \left<f(\underline{X}_1),f(\underline{X}_L)\right>_{H(k)} \\ 
\vdots& \ddots &\vdots  \\ 
\left<f(\underline{X}_L),f(\underline{X}_1)\right>_{H(k)}&\cdots  & \left<f(\underline{X}_L),f(\underline{X}_L)\right>_{H(k)}
\end{bmatrix}, 
\]
where each component can be calculated using (\ref{E47}).
Then, $\lambda_{\ell}$ for $\ell=1,...,L$ in (\ref{E30}) can be estimated by the eigenvalues $\widehat{\lambda}_{\ell}$ of $\mathcal{G}_{L}$ with descending order.
This gives another approximation test by $\overline{c}\chi^{2}_{\overline{r}}$, where $\overline{c}$ and $\overline{r}$ can be obtained via (\ref{cr_formula}) with $\widehat{E[Z]}=\text{tr} \mathcal{G}_{L}$ and $\widehat{V[Z]}=2\text{tr}  \mathcal{G}_{L}^2$.
The approximation of the null distribution by this method also performs well, but we need to adopt a relatively large value of $L$.

\begin{Rem}\label{Rem2}
Gretton et al. \cite{Gretton2009} discussed an approximation of the infinite sum of weighted chi-squared random variables, which also appeared as the asymptotic null distribution of the two-sample kernel goodness of fit test.
Their method, called {\it Spec}, uses $\sum_{\ell=1}^{L}\widehat{\lambda}_{\ell}Z_{\ell}^{2}$ as the approximation, where $Z_{\ell}$ is an independent zero-mean Gaussian random variable, and $\widehat{\lambda}_{\ell}$ eigenvalue of $\mathcal{G}_{L}$.
This method has a validity as shown in Gretton et al. (\cite{Gretton2009}, Theorem 1); however, it involves not only generating $\underline{X}_1,\dots,\underline{X}_L \overset{i.i.d.}{\sim} N( \widehat{m},\widehat{\Sigma})$ but also an additional step for generating $Z_{\ell}$s.
Further, this approach requires to calculate each eigenvalue $\hat{\lambda}_{\ell}$, which is computationally heavier than calculating the trace of the Gram matrix.
According to our experiments reported in the sequent sections, this method did not show fine accuracy in one sample testing normality.
In our comparative studies of critical points, we denote this method  $\sum \widehat{\lambda}_{\ell} \chi^2_1$.    
\end{Rem}

\subsubsection{Accuracy of the approximation}\label{s433}
We compared the critical points obtained by $n\widehat{\Delta}^2_{\text{fast}}$, $\widehat{c}\chi^2_{\widehat{r}}$, $\overline{c}\chi^2_{\overline{r}}$, and  $\sum \widehat{\lambda}_{\ell} \chi^2_{1}$ with that of $n\widehat{\Delta}^2$, where $n\widehat{\Delta}^2_{\text{fast}}$ designates the fast parametric bootstrap discussed in \cite{Kellner and Celisse} with 10,000 bootstrap iterations, and $\widehat{c}\chi^2_{\widehat{r}},~\overline{c}\chi^2_{\overline{r}}$, and $\sum \widehat{\lambda}_{\ell} \chi^2_{1}$ are methods described in Sections \ref{s431} and \ref{s432} with $L=1000$ for $\overline{c} \chi^2_{\overline{r}}$ and $L=500$ for $\sum \widehat{\lambda}_{\ell} \chi^2_{1}$.
The critical point of $n\widehat{\Delta}^2$ was determined by calculating $n\widehat{\Delta}^2$ 10,000 times under $H_0$.
Let us denote $t_{\alpha}(T_n)$ as the upper 100$\alpha$-percentile of the approximation method $T_n:T_n \in \{ n\widehat{\Delta}^2_{\text{fast}} ,~\widehat{c}\chi^2_{\widehat{r}},~\overline{c}\chi
^2_{\overline{r}},~\sum \widehat{\lambda}_{\ell} \chi^2_{1}\}$.
Let $t_{\alpha}$ be the upper 100$\alpha$-percentile of $n\widehat{\Delta}^2$ under $H_0$, determined by 10,000 simulations.
We introduce a measure for accuracy of $T_n$ defined as 
\[
\mathcal{D}(T_n)=|t_{0.1}(T_n)-t_{0.1}|+|t_{0.05}(T_n) -t_{0.05}|+|t_{0.01}(T_n)-t_{0.01}|,
\]
by which we can confirm how well $T_n$ approximates the critical points of $n\widehat{\Delta}^2$ for practically important significance levels $\alpha=0.1,~0.05,~0.01$.
Results of the comparative studies are listed in Tables {\ref{T4}} and {\ref{Table5}}, in which the obtained critical points $t_{\alpha}(T_n)$s with the smallest $\mathcal{D}(T_n)$ are italicized for each combination $(d,~\sigma,~n)$.

We observe from the case $d=10$ and $n=50$ in Table \ref{T4} that there remain differences between $t_{\alpha}$ and $t_{\alpha}(\widehat{c}\chi^2_{\widehat{r}})$, and $t_{\alpha}(\overline{c}\chi^2_{\overline{r}})$ and $t_{\alpha}(\sum \widehat{\lambda}_{\ell} \chi^2_{1})$.
However, the differences between $t_{\alpha}(n\widehat{\Delta}^2_{\text{fast}})$ and $t_{\alpha}$ are marginally large.
When the sample size $n$ grows as $n=200,~500$, it is observed that $t_{\alpha}(\widehat{c}\chi^2_{\widehat{r}})$ and $t_{\alpha}(\overline{c}\chi^2_{\overline{r}})$ get closer to $t_{\alpha}$.
In the case $d=10$, $t_{\alpha}(\overline{c}\chi^2_{\overline{r}})$ performs best, with italics for 5 out of 9 cases.
$\sum \widehat{\lambda}_{\ell} \chi^2_{1}$ performs well for 3 cases for $n=50$.

For $d=300$ in Table \ref{Table5}, the results of $t_{\alpha}(\widehat{c}\chi^2_{\widehat{r}})$ are the best in all 9 cases.
The accuracy of $t_{\alpha}(\overline{c}\chi^2_{\overline{r}})$ is marginally inferior to $t_{\alpha}(\widehat{c}\chi^2_{\widehat{r}})$.
The values of $t_{\alpha}(n\widehat{\Delta}^2_{\text{fast}})$ and $t_{\alpha}(\sum \widehat{\lambda}_{\ell} \chi^2_{1})$ tend to be larger than $t_{\alpha}(\overline{c}\chi^2_{\overline{r}})$.
Furthermore, from Table \ref{Table4}, computation of  $t_{\alpha}(\widehat{c}\chi^2_{\widehat{r}}) $ is very fast even in the case of $d=300$.
From the perspectives of accuracy and computation, we strongly recommend $t_{\alpha} (\widehat{c}\chi^2_{\widehat{r}})$ to approximate critical points to testing normality by $n\widehat{\Delta}^2$.
\begin{table}[!h]
	\caption{Approximation of critical points in the case of $d=10$.
	For each $(\sigma,~n)$, the values of the method with the smallest $\mathcal{D}$ are italicized.
}
	\centering
	\label{T4}
	\begin{tabular}{cccc|c|cccc}
		\hline 
		$d$& $\sigma$ & $n$ &$\alpha$ & $n \widehat{\Delta}^2$ & $n\widehat{\Delta}^2_{\text{fast}}$& $\widehat{c} \chi^2_{\widehat{r}}$ &$\overline{c}\chi^2_{\overline{r}}$& $\sum \widehat{\lambda}_{\ell} \chi^2_{1}$\\ 
		\hline 
		&  &  & 10\% & 0.65912& 0.70042&0.61394&  0.62546&\textit{0.62840}\\ 
		&  &50  & 5\% &  0.67780& 0.76462 &0.62934  &  0.64506 &\textit{0.65503}\\ 
		&  &  & 1\%  &  0.71336&0.88697 &  0.65893 & 0.68292&\textit{0.70689}\\ 
		\cline{3-9}
		&  &  & 10\% &0.65930 &0.72457 &0.65051  & \textit{0.66434}& 0.66774\\ 
		& $d^{-3/4}$  &200  & 5\% &  0.67473 &0.75533 & 0.66423 & \textit{0.68228}&  0.69170\\ 
		&  &  & 1\%  &   0.70517&0.82534 & 0.69047 &\textit{0.71678}&0.73832\\ 
		\cline{3-9}
		& &  & 10\% &  0.65866 &0.68174 &   0.66273 &  \textit{0.65964}&0.70107\\ 
		&  &500  & 5\% &  0.67432 &0.70393 &  0.67607& \textit{0.67718}&0.72434\\ 
		&  &  & 1\%  & 0.70556&0.74606 & 0.70157& \textit{0.71091}&0.76823\\
		\cline{2-9} 
		&  &  & 10\% &0.50757&0.53110 &0.46204&0.47147&\textit{0.47384}\\ 
		&  &50  & 5\% &0.52639&0.58462 &   0.47582&0.48856&\textit{0.49616}\\ 
		&  &  & 1\% &0.56313&0.69568 &  0.50241&0.52172&\textit{0.54386}\\ 
		\cline{3-9}
		&  &  & 10\% &0.50482&0.56084 &0.49569&\textit{0.50664}&0.50692\\ 
		10&$d^{-7/8}$  &200  & 5\% &0.51882&0.58878 &0.50815&\textit{0.52217}&0.52740\\ 
		&  &  & 1\%  &0.54976&0.65123 & 0.53206&\textit{0.55213}& 0.56914\\ 
		\cline{3-9}
		&  && 10\% &0.50323&0.52011& 0.50843 &\textit{0.50132}&0.53878\\ 
		&  &500  & 5\% &0.51874&0.53945&  0.52062& \textit{0.51646}&0.55981\\ 
		&  &  & 1\%  &0.54953& 0.57648 &  0.54401&\textit{0.54568}& 0.59976\\
		\cline{2-9} 
		&  &  & 10\% &0.35914& 0.37574 &0.32004&0.32841&\textit{0.33003}\\ 
		&  &50  & 5\% &0.37596&0.41776 &0.33107&0.34199&\textit{0.34750}\\ 
		&  &  & 1\%  & 0.40931&0.51342 &  0.35242&0.36846&\textit{0.38651}\\ 
		\cline{3-9}
		&  &  & 10\% &0.35518&0.40127 &   0.34713&\textit{0.35483}&0.35383\\ 
		&$d^{-1}$   &200  & 5\% &0.36719&0.42454&0.35718& \textit{0.36702}&0.36980\\ 
		&  &  & 1\%  &0.39313&0.47690& 0.37655&\textit{0.39064}&0.40436\\ 
		\cline{3-9}
		& &  & 10\% & 0.35361& 0.36538 & \textit{0.35854}& 0.35008&0.38070\\ 
		&  &500  & 5\% &0.36639&0.38082 &\textit{0.36844}&0.36194&0.39686\\ 
		&  &  & 1\%  &0.39102&0.41033 &\textit{0.38748}&0.38492&0.43111\\
		\hline
	\end{tabular} 
\end{table}
\begin{table}[!h]
	\caption{Approximation of critical points in the case of $d=300$.
	For each $(\sigma,~n)$, the values of the method with the smallest $\mathcal{D}$ are italicized.	
	}
	\centering
	\label{Table5}
	\begin{tabular}{cccc|c|cccc}
		\hline 
		$d$& $\sigma$ & $n$ &$\alpha$ & $n \widehat{\Delta}^2$ &$n\widehat{\Delta}^2_{\text{fast}}$ & $\widehat{c} \chi^2_{\widehat{r}}$ &$\overline{c}\chi^2_{\overline{r}}$& $\sum \widehat{\lambda}_{\ell} \chi^2_{1}$\\ 
		\hline 
		&  &  & 10\% & 0.98561& 1.19361 &\textit{0.97748}&1.03114&1.05425\\ 
		&  &50  & 5\% &0.98597&1.27905 &\textit{0.97789}&1.04768&1.07876\\ 
		&  &  & 1\%  &0.98656& 1.44868 &\textit{0.97865}&1.07917&1.12624\\ 
		\cline{3-9}
		&  &  & 10\% & 0.98762&1.10135 &\textit{0.98497}&1.04055&1.06310\\ 
		& $d^{-3/4}$ &200  & 5\% & 0.98777& 1.14276& \textit{0.98502}&1.05723&1.08773\\ 
		&  &  & 1\%  & 0.98807&1.22131 &\textit{0.98512}&1.08900&1.13519\\ 
		\cline{3-9}
		&  &  & 10\% &0.98787&1.06465 &\textit{0.98664}&1.04276&1.06491\\ 
		&  &500  & 5\% &0.98797& 1.08913 &\textit{0.98666}&1.05948&1.08951\\ 
		&  &  & 1\%  &0.98816&1.13737 &\textit{0.98671}& 1.09131& 1.13753\\
		\cline{2-9} 
		&  &  & 10\% &0.76002&0.86448 &\textit{0.74657}&0.77779&0.79896\\ 
		&  &50  & 5\% &0.76236&0.92857 &\textit{0.74802}&0.79052&0.81836\\ 
		&  &  & 1\%  &0.76659& 1.04667 &\textit{0.75076}&0.81477&0.85305\\ 
		\cline{3-9}
		&  &  & 10\% &0.76932&0.83429 &\textit{0.76045}&0.79675&0.81600\\ 
		300& $d^{-7/8}$  &200& 5\% &0.77056&0.86446 &\textit{0.76074}&0.80959&0.83529\\ 
		&  &  & 1\%  &0.77262&0.92360 &\textit{0.76129}&0.83403&0.87112\\ 
		\cline{3-9}
		& &  & 10\% & 0.77039& 0.81813 &\textit{0.76531}&0.80477&0.82143\\ 
		&  &500  & 5\% &0.77117&0.83744 &\textit{0.76545}&0.81771&0.84125\\ 
		&  &  & 1\%  &0.77259&0.87495 &\textit{0.76573}&0.84235&0.87719\\
		\cline{2-9} 
		&  &  & 10\% & 0.31768&0.35479 &\textit{0.32202}&0.33189&0.34164\\ 
		&  &50  & 5\% &0.31984&0.38115 &\textit{0.32324}&0.33769&0.35025\\ 
		&  &  & 1\%  &0.32402&  0.42868&\textit{ 0.32555}&0.34875&0.36623\\ 
		\cline{3-9}
		&  &  & 10\% &0.32391&0.34482 &\textit{0.32019}&0.33250&0.34227\\ 
		& $d^{-1}$ &200  & 5\% &0.32510&0.35807 &\textit{0.32045}&0.33796&0.35075\\ 
		&  &  & 1\%  &0.32701&0.38199 &\textit{0.32095}&0.34836&0.36615\\ 
		\cline{3-9}
		&  &  & 10\% &0.32447&0.34101 &\textit{0.32155}&0.33623&0.34390\\ 
		&  &500  & 5\% &0.32521&0.34939 &\textit{0.32168}&0.34170&0.35238\\ 
		&  &  & 1\%  &0.32658&0.36466 &\textit{0.32194}&0.35213&0.36736\\
		\hline		
	\end{tabular} 
\end{table}
\begin{table}[!h]
	\centering
	\caption{Comparison of execution times (sec) of $t_{\alpha}(n\widehat{\Delta}^2_{\text{fast}})$, $t_{\alpha}(\widehat{c}\chi^2_{\widehat{r}})$, $t_{\alpha}(\overline{c}\chi^2_{\overline{r}})$, and $t_{\alpha}(\sum \widehat{\lambda}_{\ell} \chi^2_{1})$, where $\alpha=0.05$.
	Computations were coded in {\tt R} and implemented under a Windows machine with Intel(R) Core(TM) i7-8700 CPU @ 3.20GHz and 16.0GB memory.
	}
	\label{Table4}
	\begin{tabular}{cc|cccc}
		\hline 
		$d$&$n$ &$n\widehat{\Delta}^2_{\text{fast}}$  & $\widehat{c}\chi^2_{\widehat{r}}$ & $\overline{c}\chi^2_{\overline{r}} $& $\sum\widehat{\lambda}_{\ell} \chi^2_1$\\ 
		\hline 
		& 50 & 0.87 &0.04& 3.02 & 1.22\\ 
		10& 200 &  3.01 & 0.02 &  3.05 & 1.18\\ 
		& 500 &  15.20& 0.02  &  2.99 & 1.22\\ 
		\hline
		& 50 &  447.97&  0.13& 6.38 & 2.22\\ 
		300& 200 & 572.39 & 0.12 & 6.38 & 2.16\\ 
		& 500 &  907.39& 0.12 &  6.46& 2.26\\ 	
		\hline 
	\end{tabular} 
\end{table}

\subsection{Simulation}\label{s4.4}

In this section,  we investigate the performance of $n\widehat{\Delta}^2$ under a specific alternative hypothesis.
In particular, a Monte Carlo simulation is carried out to see the power of the test against a uniform distribution and an exponential distribution, both of which are standardized.
Two cases are implemented: independent components and correlated components with the correlation matrix $R=\left((1/2)^{|i-j|} \1_{|i-j|\leq 5}\right)_{1 \leq i,j \leq d}$.
The rejection point is determined on the basis of 10,000 simulations of $n\widehat{\Delta}^2$ under the standard normal distribution. 
Then, the estimated power of $n\widehat{\Delta}^2$ can be obtained by counting how many times $n\widehat{\Delta}^2$ exceeds the rejection point in 1000 iterations under each alternative distribution.
We execute the above for $n=200, ~300,~ 400,$ and 500 and $d=10$ and 300. 
The case of $d=10$ addresses the usual testing of multi-normality for large samples.
In contrast, the performance of $n\widehat{\Delta}^2$ for high-dimension data is investigated for the case of $d=300$.
We focus on $d=300$.

The Gaussian kernel (\ref{Gaussian kernel}) was used throughout.
It is known that the selection of the value of $\sigma$ involved in the Gaussian kernel affects the performance.
We utilize $\sigma$ depending on dimension $d$.
The results of simulations are presented in Tables \ref{T1} and \ref{Table2}.
\begin{table}[!h]
	\centering
	\caption{Test power for each sample size and each parameter $\sigma$ (independent case).}
	\label{T1}
	\begin{tabular}{cccccccccc}
		\hline
		\multicolumn{10}{c}{Uniform}\\
		\hline
		$d$&\multicolumn{4}{c}{10}&&\multicolumn{4}{c}{300}\\
		\cline{2-5}
		\cline{7-10}
		$n$&200&300&400&500&&200&300&400&500\\
		\cline{2-5}
		\cline{7-10}
		$\sigma=d^{-3/4}$&1&1&1&1&&0.206&0.420&0.596&0.794\\
		$\sigma=d^{-7/8}$&0.983&1&1&1&&0.004&0.009&0.014&0.015\\
		$\sigma=d^{-1}$&0.863&0.999&1&1&&0&0&0.001&0.001\\
		\hline
		\multicolumn{10}{c}{Exponential}\\
		\hline
		$d$&\multicolumn{4}{c}{10}&&\multicolumn{4}{c}{300}\\
		\cline{2-5}
		\cline{7-10}
		$n$&200&300&400&500&&200&300&400&500\\
		\cline{2-5}
		\cline{7-10}
		$\sigma=d^{-3/4}$&1&1&1&1&&0.144&0.581&0.956&0.998\\
		$\sigma=d^{-7/8}$&1&1&1&1&&0.675&0.931&0.995&1\\
		$\sigma=d^{-1}$&1&1&1&1&&0.901&0.988&1&1\\
		\hline
	\end{tabular}
	\vspace{1cm}
	\centering
	\caption{Test power for each sample size and each parameter $\sigma$ (correlate case).}
	\label{Table2}
	\begin{tabular}{cccccccccc}
		\hline
		\multicolumn{10}{c}{Uniform}\\
		\hline
		$d$&\multicolumn{4}{c}{10}&&\multicolumn{4}{c}{300}\\
		\cline{2-5}
		\cline{7-10}
		$n$&200&300&400&500&&200&300&400&500\\
		\cline{2-5}
		\cline{7-10}
		$\sigma=d^{-3/4}$&0.747&0.958&0.996&1&&0.156&0.273&0.367&0.493\\
		$\sigma=d^{-7/8}$&0.537&0.848&0.980&0.997&&0.023&0.030&0.031&0.027\\
		$\sigma=d^{-1}$&0.290&0.618&0.885&0.973&&0.005&0.004&0.005&0.003\\
		\hline
		\multicolumn{10}{c}{Exponential}\\
		\hline
		$d$&\multicolumn{4}{c}{10}&&\multicolumn{4}{c}{300}\\
		\cline{2-5}
		\cline{7-10}
		$n$&200&300&400&500&&200&300&400&500\\
		\cline{2-5}
		\cline{7-10}
		$\sigma=d^{-3/4}$&1&1&1&1&&0.312&0.853&0.991&1\\
		$\sigma=d^{-7/8}$&1&1&1&1&&0.788&0.978&1&1\\
		$\sigma=d^{-1}$&1&1&1&1&&0.925&0.995&1&1\\
		\hline
	\end{tabular}
\end{table}
\begin{figure}[!h]
	\centering
	\begin{minipage}{0.49\hsize}
		\includegraphics[width=1\linewidth]{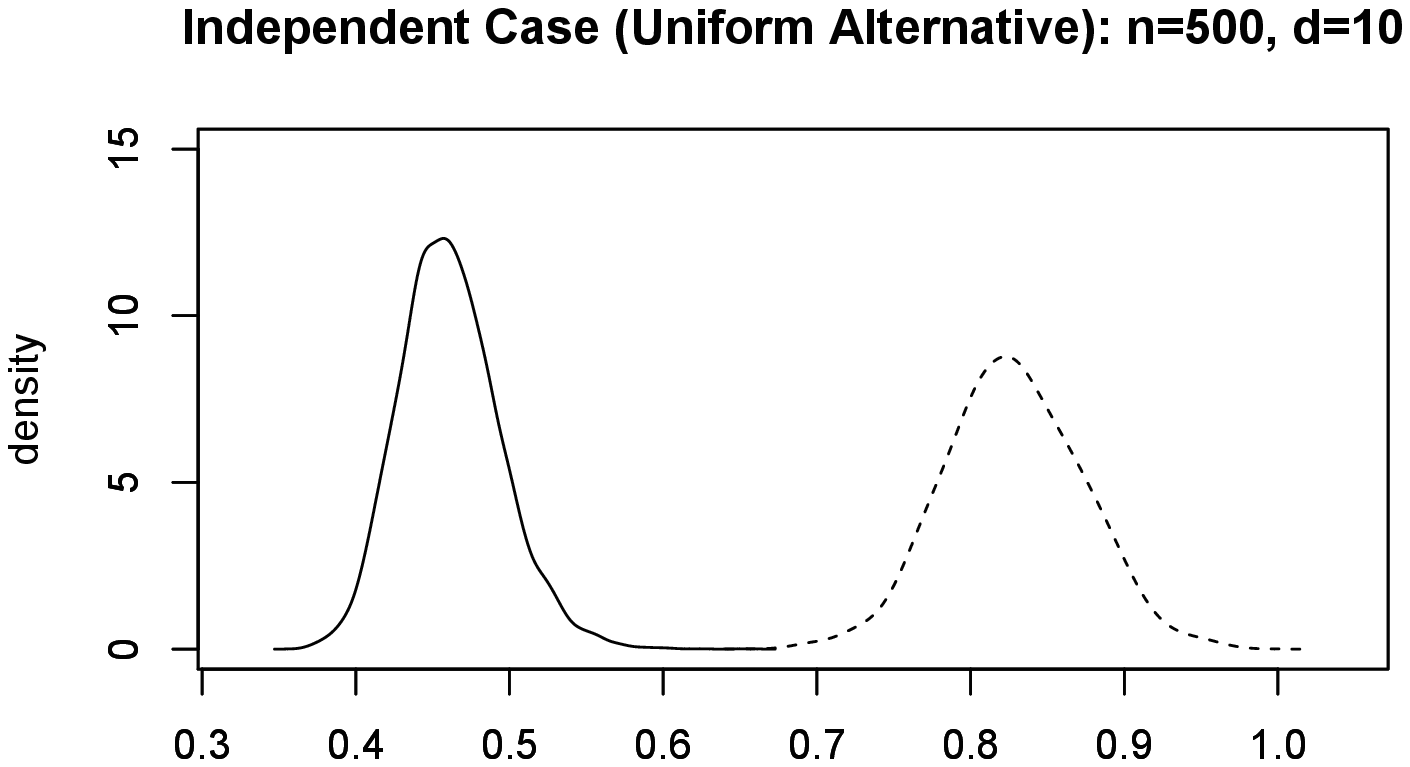}
		\vspace{-2cm}
		\caption*{(a)}
	\end{minipage}
	\begin{minipage}{0.49\hsize}
		\includegraphics[width=1\linewidth]{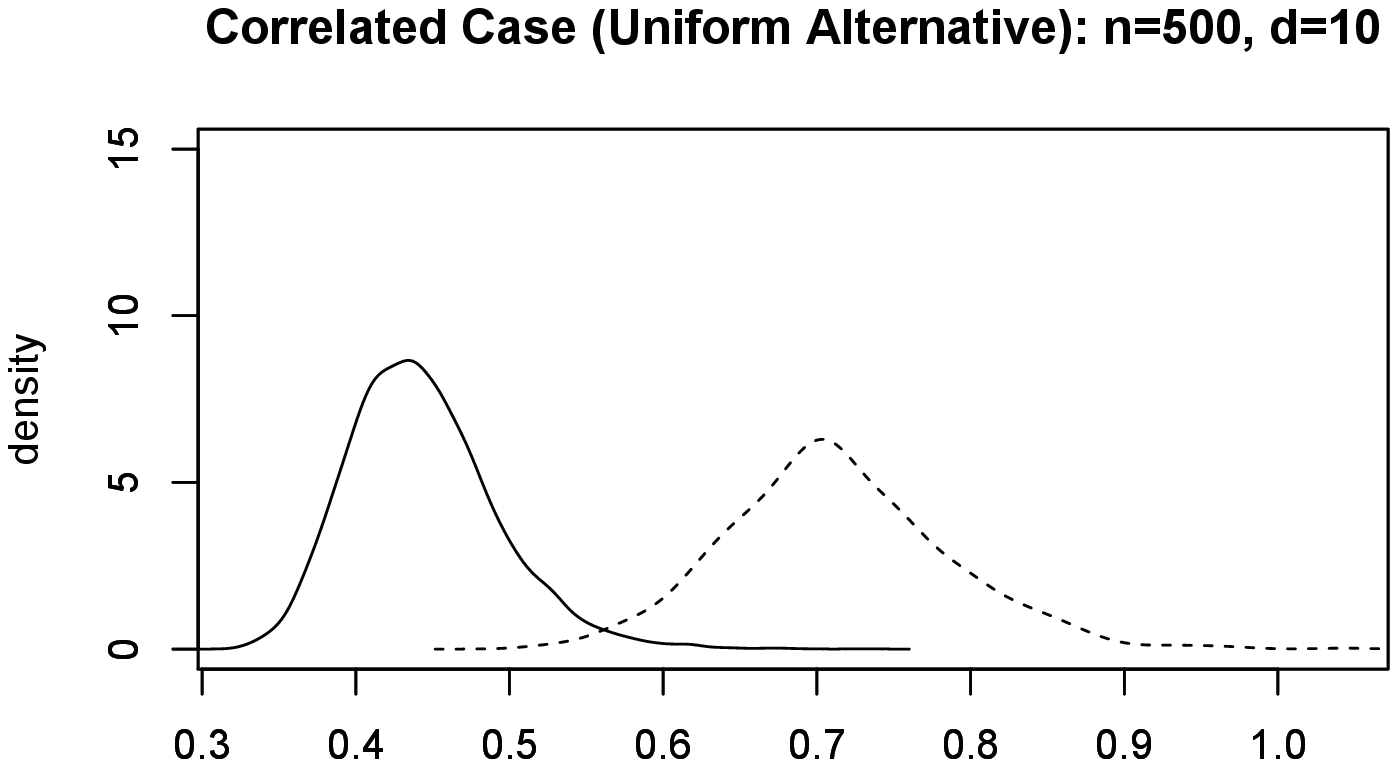}
		\vspace{-2cm}
		\caption*{(b)}
	\end{minipage}
	\begin{minipage}{0.49\hsize}
		\includegraphics[width=1\linewidth]{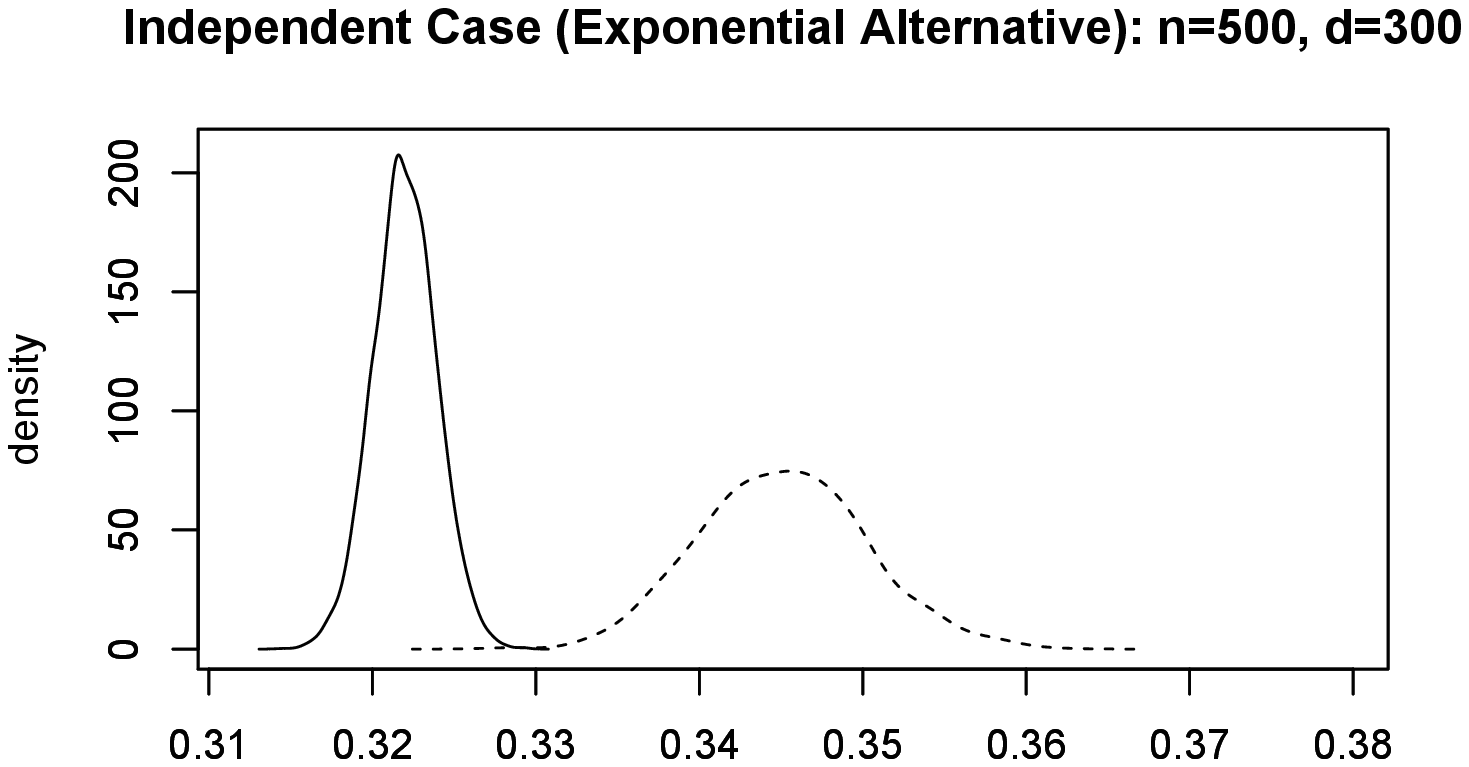}
		\vspace{-2cm}
		\caption*{(c)}
	\end{minipage}
	\begin{minipage}{0.49\hsize}
		\includegraphics[width=1\linewidth]{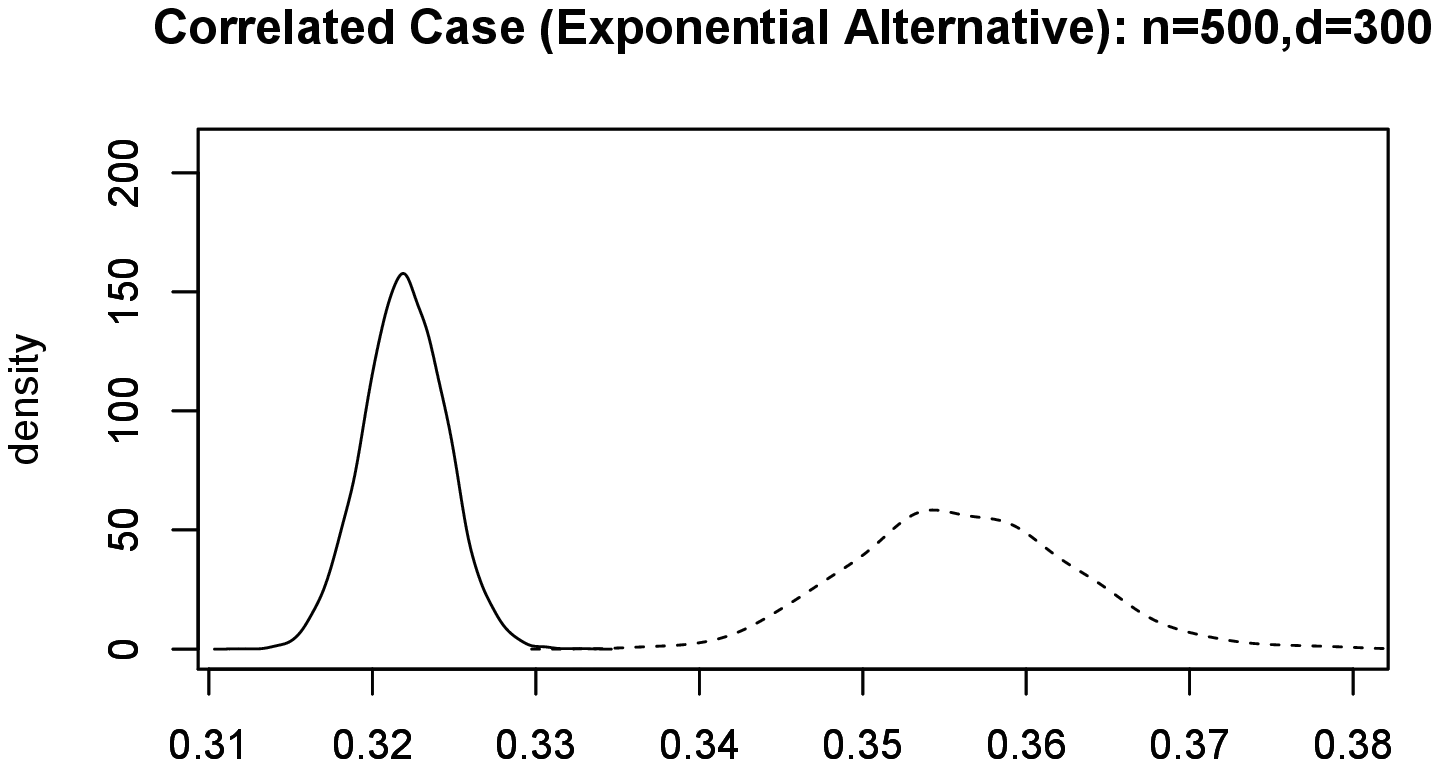}
		\vspace{-2cm}
		\caption*{(d)}
	\end{minipage}
	\caption{Density estimates for various distributions: The panel (a) is the density estimate of the distribution of $n\widehat{\Delta}^2$ for a standardized uniform distribution with $n = 500$, $d = 10$.
	The panel (b) represents the density estimate for a uniform distribution with covariance matrix $\Sigma = R$.
	The panel (c) is the density estimate of the distribution of $n\widehat{\Delta}^2$ for a standardized exponential distribution with $n = 500$, $d = 300$.
	The panel (d) shows the density estimate for an exponential distribution with covariance matrix $\Sigma = R$.
	The uniform distribution uses a Gaussian kernel with $\sigma=d^{-7/8}$ and the exponential distribution uses a Gaussian kernel with $\sigma=d^{-1}$.}
	\label{F1}
\end{figure}
\begin{itemize}
	\item Table \ref{T1} shows that $n\widehat{\Delta}^2$ has high power against both uniform and exponential distributions in the case of $d=10$.
	Such results hold for all simulated values of $\sigma$.
	\item  From Table \ref{T1}, in the case of $d=300$, the power of the test against uniform distribution did not reach the significance level for $\sigma = d^{-7/8}$ and $d^{-1}$. 
	Even for $\sigma=d^{-3/4}$, the high-dimensional test performed poorly with $n=200$  and 300. 
	Thus, that selecting the parameter $\sigma$ affects the behavior of the test.
	
	\item In all cases with $d=300$ for the exponential distribution exhibited in Table \ref{T1}, when $n$ gets bigger, the power of the test approaches 1,
	which reveals that the consistency of the test as mentioned in Remark \ref{Rem1} certainly holds.
	
	\item From Table \ref{T1}, it might be difficult to find the difference between the normal distribution and a compactly supported distribution such as the uniform distribution, because the power against the exponential distribution is higher than that against the uniform distribution for almost all simulated cases.
	
	\item Table \ref{Table2} includes the results for the correlated case; however, the tendency of the results in Table \ref{Table2} is similar to that of the results in Table \ref{T1} for the independent case.      
	
	\item Comparing Table \ref{T1} with Table \ref{Table2}, we learn that the power of the test against the uniform distribution decreases, while the power against the exponential distribution increases for the correlated case.
	To get a deeper understanding of this phenomenon, we focus on the simulated distribution of $n\widehat{\Delta}^2$ as seen in Figure \ref{F1}.
	In Figure \ref{F1}, the kernel density estimates of $n\widehat{\Delta}^2$ are drawn, where the solid line is $n\widehat{\Delta}^2$ under normality and the dashed line is $n\widehat{\Delta}^2$ under an alternative distribution. 
	
	\item  First, we consider the panels (a) and (b) in Figure \ref{F1}.
	For the uniform distribution, the distribution of   $n\widehat{\Delta}^2$ in the correlated case (dashed line) shows bigger variance than the independent case, and the location of the distribution is smaller as compared with the case without correlation.	
	Even under normality (solid line), the variance in the correlated case is slightly bigger than in the independent case, and the location in the correlated case is smaller than the independent case.
	However, the magnitude of the changes in variance and location of $n\widehat{\Delta}^2$ under normality (solid line) is much smaller than those under uniform distribution (dashed line).
	Therefore, the decrease in power against the uniform distribution can be explained by the movement of the distribution of $n\widehat{\Delta}^2$.
	
	\item Next, we focus on the panels (c) and (d) of  Figure \ref{F1}.	
	For the exponential distribution (dashed line), the existence of correlation among components increases the variance and  the location for the distribution of  $n\widehat{\Delta}^2$.
	On the other hand, the existence of correlation among components yields only a slight increase in the variance under the normal distribution (solid line).	
	Hence, the increase in power against the exponential distribution can be explained by the movement of the distribution of $n\widehat{\Delta}^2$.
\end{itemize}

\section{Applications to real data sets}\label{s5}
The kernel normality test was applied to some real data sets, where the scale parameter $\sigma^{-1}$ was determined as the median of $\left\|\underline{X}_i -\underline{X}_j\right\|^2_{\mathbb{R}^d} ,~~i < j$, where $\underline{X}_i$ is the $i$-th data.
The significance level was 0.05, and the critical point $t_{0.05}$ was obtained by 10,000 or 1000 iterations of calculating $n\widehat{\Delta}^2$  based on the sample of size $n$ drawn from $N(\widehat{\underline{m}},\widehat{\Sigma})$, where $\widehat{\underline{m}}$ and $\widehat{\Sigma}$ are the sample mean and the sample variance calculated from the real data set.
We also calculated the critical points $t_{0.05}(\widehat{c}\chi^2_{\widehat{r}})$, $t_{0.05}(\overline{c}\chi^2_{\overline{r}})$ and $t_{0.05}(\sum \widehat{\lambda}_{\ell} \chi^2_1)$ of the approximation tests by 
$\widehat{c}\chi^2_{\widehat{r}}$, $\overline{c}\chi^2_{\overline{r}}$ and $\sum \widehat{\lambda}_{\ell} \chi^2_1$, respectively (see Section \ref{s4.3}).
\subsection{USPS data}
The USPS dataset consists of handwritten digits represented by a $16\times16$ greyscale matrix (\textit{https://www.csie.ntu.edu.tw/\textasciitilde cjlin/libsvmtools/datasets/multiclass.html\#usps}).
Our question is which group is most deviated from normality among a set of groups numbered 0 to 9.
The $t_{0.05}$ value was calculated by a simulation of 10,000 iterations.
It can be seen from Table \ref{Table3} that normality was strongly rejected for all groups.

\subsection{MNIST data}
The MNIST dataset consists of $28 \times 28$ pixels (\textit{http://yann.lecun.com/exdb/mnist}).
The $t_{0.05}$ value was obtained by a simulation with 1000 iterations.
Similar to the USPS data, the normality hypothesis was rejected for all groups.

\subsection{Leukemia data}
The Leukemia dataset contains the gene expression level of leukemia patients with $d=7128$ and $n=72$, among which 25 have Acute Myelogenous Leukemia (AML) and the remaining 47 have Acute Lymphocytic Leukemia (ALL) (\textit{https://web.stanford.edu/~hastie/CASI\_files/DATA/leukemia.html}).
As in the MNIST data, $t_{0.05}$ was obtained with 1000 iterations.
Table 3 shows that the normality hypothesis was rejected.

\subsection{Colon data}
The Colon dataset contains gene expression data from a DNA microarray experiments of colon tissue samples with $d=2000$ and $n=62$ (see \cite{Alon} for details). 
Among 62 samples, 40 are tumor tissues and 22 are normal tissues.
The $t_{0.05}$ value was calculated by simulation with 1000 iterations.
It is seen from Table 3 that the normality hypothesis for tumor tissue data was rejected with $p$-value 0.033, while that on normal tissue data was not rejected with $p$-value 0.235.
This result is different to that reported in \cite{Himeno_Yamada}, where the normality hypothesis on the tumor tissue data was not rejected.

\subsection{Lung Cancer data}

The Lung Cancer dataset contains gene expression level data of lung cancer patients with $d=12533$ and $n=32$ (see \cite{Gordon} for details).
In this data, 16 out of 32 are on adenocarcinoma (ADCA)  and the remaining 16 are on mesothelioma (MPM).
We applied our normality test to two sets of Lung Cancer data: a subset of Lung Cancer data made by adopting a certain screening of variables, and the full data.

\subsubsection{Screened Data}
Screening was applied to the Lung Cancer dataset regarding whether the sample  variance of each variable is bigger than 7000. This left 1911 variable in both ADCA and MPM.
The results of the kernel test on the ADCA and MPM screened data sets with 1911 dimension are presented in Table \ref{Table3}, where $t_{0.05}$ was calculated by a simulation iterated 1000 times.
Table \ref{Table3} shows that the normality hypothesis was strongly rejected for both ADCA and MPM.

\subsubsection{Full Data}
We did not implement the simulation to obtain $t_{0.05}$ of the original full data, since the calculation cost is too large due to the dimension $d=12533$.
Hence, approximation tests by $t_{0.05}(\widehat{c} \chi^2_{\widehat{r}})$, $t_{0.05}(\overline{c}\chi^2_{\overline{r}})$, and $t_{0.05}(\sum \widehat{\lambda}_{\ell} \chi^2_1)$ were exploited for the full data.
Table \ref{Table3} shows that the normality hypothesis on ADCA was rejected by all approximation tests.
For MPM, the value of $t_{0.05}(\widehat{c} \chi^2_{\widehat{r}})$ was small compared with the values of $t_{0.05}(\overline{c} \chi^2_{\overline{r}})$ and $t_{0.05}(\sum \widehat{\lambda}^*_{\ell} \chi^2_{1})$, however the normality hypothesis on MPM was rejected regardless of the examined approximation tests.
The method by $\widehat{c}\chi^2_{\widehat{r}}$ was calculated using only 16 samples, whereas the methods of $\overline{c} \chi^2_{\overline{r}}$ and $\sum\widehat{\lambda}_{\ell} \chi^2_{1}$ were calculated by generating random numbers of size 1000 or 500.
Due to this difference in sample size for calculating  percentiles, it seems that $t_{0.05}(\widehat{c}\chi^2_{\widehat{r}})$ becomes a small value compared with $t_{0.05}(\overline{c}\chi^2_{\overline{r}})$ and $t_{0.05}(\sum \widehat{\lambda}_{\ell} \chi^2_{1})$.

\subsection{Approximation tests}
In the examined applications to real data sets, $t_{0.05}(\widehat{c} \chi^2_{\widehat{r}})$, $t_{0.05}(\overline{c} \chi^2_{\overline{r}})$, and $t_{0.05}(\sum \widehat{\lambda}_{\ell} \chi^2_{1})$ were all similar to the value of $t_{0.05}$. 
Because the calculation cost of $t_{0.05}$ is in fact very large, it is practically sufficient to use a rejection point by the approximation method.
Based on the results reported in Sections \ref{s433} and \ref{s5}, we specifically suggest using the approximation test by $\widehat{c} \chi^2_{\widehat{r}}$.

\begin{table}[!h]
	\centering
	\caption{Summary of the results for the kernel test.}
	\label{Table3}
	\vspace{-0.5cm}
	{\scriptsize
	\begin{tabular}{cccccccccc}
		\multicolumn{10}{c}{ (\scriptsize Symbol meaning, $\text{**}:\leq 10^{-4},~\text{*}:\leq 10^{-3}$, $\dagger:\times10^{-7},~\ddagger:\times 10^{-4},~\bullet:\times 10^{-3},~\star:\times 10^{-10}$)} \\
		\hline
		Data set&  & $n$ & $\sigma$ & $n\widehat{\Delta}^2_P$ & $t_{0.05}$ & $t_{0.05}(\widehat{c}\chi^2_{\widehat{r}})$& $t_{0.05}(\overline{c}\chi^2_{\overline{r}})$& $t_{0.05}(\sum \widehat{\lambda}_{\ell} \chi^2_1)$&$p$-value \\ 
		\hline 
		&  0&  359&0.055& 4.04114 &0.88192&0.87601 &0.90331&0.92557 &**\\ 
		& 1 &  264&1.178& 5.36517 &1.00094&1.00087 &1.06902 &1.10103 &** \\  
		& 2 & 198 &0.025&1.33186& 0.58824 &0.58137&0.60150&0.60195&** \\ 
		&  3& 166 &0.035&1.45951& 0.62455&0.61966 & 0.62526&0.63326 & **\\ 
		&  4& 200 &0.048&2.63348&0.79405&0.79013& 0.83553&0.84957 &** \\ 
		USPS&  5& 160 &0.025&1.24266&0.57249 &0.56551& 0.57048& 0.57925&** \\ 
		($d=256$)&  6& 170 & 0.065& 2.86309 &0.83082&0.82816& 0.85809&0.87814 & **\\ 
		&  7& 147 & 0.080& 2.46502 &0.84663&0.84326& 0.86484&0.88210 &** \\ 
		&  8& 166 &0.036&2.11776&0.64761&0.64207&0.66407 &0.67116 &** \\ 
		&  9& 177 &0.087&2.88034&0.83163&0.82797 &0.84746 &0.86446 &** \\ 		
		\hline 
		&  0&  980& $5.877 \dagger$ & 4.64786 & 0.69158& 0.68850 &0.71383& 0.73699& *\\ 
		& 1 & 1135 & $33.15 \dagger$ & 10.6682& 0.96064&0.96079& 1.02128&1.05030 &* \\  
		& 2 & 1032 & $4.411 \dagger$ & 2.44670 & 0.53862 &0.53728&0.56601 & 0.58320& *\\ 
		&  3& 1010 & $5.323 \dagger$ & 2.89616 & 0.56813&0.56723& 0.59033&0.60320 & *\\ 
		&  4& 982 & $6.352 \dagger$ & 3.65969 & 0.63293& 0.63020& 0.66788&0.68121 & *\\ 
		MNIST&  5& 892 & $5.144 \dagger$ &3.01429& 0.59798 & 0.59439& 0.62334&0.63663 & *\\ 
		($d=784$)&  6& 958 & $6.525 \dagger$ & 4.74724 & 0.68923& 0.68755& 0.71695& 0.74202& *\\ 
		&  7& 1028 & $8.946 \dagger$ & 5.67498 & 0.76170& 0.75874& 0.80097& 0.82295& *\\ 
		&  8& 974 & $4.726 \dagger$ & 2.63448 & 0.51834 &0.51646 &0.54934 & 0.55618& *\\ 
		&  9& 1009 & $7.714 \dagger$ &5.20226 & 0.69505 &0.69277 & 0.73342& 0.76438& *\\ 		
		\hline 
		Leukemia&AML&25&$5.466 \ddagger$&0.15362&0.17208&0.16294&0.16921 &0.17859 &0.206\\
		($d=7128$)&ALL&47&$5.441 \ddagger$& 0.16651 &0.16525&0.158322& 0.16312& 0.17122&0.039\\
		\hline
		Colon&tumor&40&$1.659 \bullet$& 0.26163&0.25473&0.23996 &0.24394 &0.24523 &0.033 \\
		($d=2000$)&normal&22&$1.604 \bullet$&0.20025&0.24011&0.22328& 0.21788&0.21616 &0.235 \\
		\hline
		Lung Cancer&ADCA&16&$21.24 \star$&1.07492&0.54126&0.52619 & 0.54891& 0.58937&** \\
		(Sub; $d=1911$)&MPM&16&$6.306 \star$&0.51774&0.40507&0.39538&0.37131 &0.39814 &0.009 \\
		\hline
		Lung Cancer&ADCA&16&$17.56 \star$&0.92992&-& 0.47393&0.49085&0.49614& -\\
		(Full; $d=12533$)&MPM&16&$5.921\star$&0.50037&-&0.10981&0.39044&0.39612&-\\
		\hline 
	\end{tabular} }
\end{table}

\section{Conclusion}\label{s6}
We derived an asymptotic non-null distribution of the MMD test in Section \ref{s3.3}, which was a normal distribution.
From this asymptotic normality of the test under alternative hypotheses, we found that the MMD test for normality has consistency.  
We developed an asymptotic distribution of the test under a sequence of local alternatives in Section \ref{s3.4}.
This was in the form of an infinite sum of weighted noncentral chi-squared distribution. 
Further, we derived the asymptotic null distribution using the results of the asymptotic distribution under local alternative hypotheses.
We found that the asymptotic null distribution had the form of an infinite sum of weighted central chi-squared distribution, where the weights are the same as those in the asymptotic distribution of the test under a sequence of local alternatives.
Section \ref{s4} examined the test statistic when kernel $k$ is a Gaussian kernel and the Hilbert space $\mathcal{H}$ is $\mathbb{R}^d$.
The $\alpha$-quantiles in the asymptotic null distribution have been well-approximated by a single weighted chi-squared distribution.
In the simulation of the power reported in Section \ref{s4.4}, we found that the power of the test against the exponential distribution approached 1 as $n \to \infty$.
We saw in Section \ref{s5} that the MMD test works for high-dimension low-sample size real data sets, and we recommended using the approximation test by $\widehat{c} \chi^2_{\widehat{r}}$.
\section{Proof}\label{s7}
\subsection{Lemmas for Theorem \ref{T2}}
This section presents a series of lemmas all of which are necessary to obtain theorems as well as theoretical formulas.
\begin{lemm}\label{L5}
	Assume that $A:\mathcal{H} \oplus HS(\mathcal{H}) \to H(k)$ is bounded linear operator (see Section I.1 of \cite{Reed and Simon} for details).
	Then
	\[
	\left\|A(0,(\widehat{m}-m_0)^{\otimes 2})\right\|_{H(k)}=O_p\left(\frac{1}{n}\right).
	\]
\end{lemm}
\begin{lemm}\label{L6}
	Assume that $\mu(N(\cdot,\cdot)) \in C^2(B((m_0,\Sigma_0),\varepsilon),H(k))$.
	Then
	\[
	\left\|
	\int^{1}_{0} (1-s) D^2_{(m_0,\Sigma_0) +s(\widehat{m}-m_0,\widehat{\Sigma}-\Sigma_0)}\mu(N(\widehat{m}-m_0,\widehat{\Sigma}-\Sigma_0)^2) ds
	\right\|_{H(k)}
	=O_p\left(\frac{1}{n}\right),
	\]
	where
	$D^2_{(m_0,\Sigma_0) +s(\widehat{m}-m_0,\widehat{\Sigma}-\Sigma_0)}\mu(N(\cdot,\cdot))$ is twice Fr\'echet derivative of $\mu(N(\cdot,\cdot))$
	and
	\[
	(\widehat{m}-m_0,\widehat{\Sigma}-\Sigma_0)^2=((\widehat{m}-m_0,\widehat{\Sigma}-\Sigma_0),(\widehat{m}-m_0,\widehat{\Sigma}-\Sigma_0)) \in (\mathcal{H} \oplus HS(\mathcal{H})) \oplus (\mathcal{H}\oplus HS(\mathcal{H}))
	\]
	and the integral is the Bochner integral in Chapter III of \cite{Mikusinski}.
\end{lemm}
\subsection{Proof of Theorem \ref{T2}}
Let us first expand the following quantity
\begin{align*}
&\sqrt{n}(\widehat{\Delta}^2-\Delta^2)\\
&=\sqrt{n}\left\{
\left\|
\frac{1}{n}\sum_{i=1}^{n}k(\cdot,Y_i)-\mu(N(\widehat{m},\widehat{\Sigma}))
\right\|^2_{H(k)}
-\left\|\mu(P)-\mu(N(m_0,\Sigma_0))\right\|^2_{H(k)}
\right\}\\
&=\sqrt{n}\left<
\frac{1}{n}\sum_{i=1}^{n}k(\cdot,Y_i)-\mu(N(\widehat{m},\widehat{\Sigma}))-(\mu(P)-\mu(N(m_0,\Sigma_0))),\right.\\
&\hspace{50mm}\left.\frac{1}{n}\sum_{i=1}^{n}k(\cdot,Y_i)-\mu(N(\widehat{m},\widehat{\Sigma}))+(\mu(P)-\mu(N(m_0,\Sigma_0)))
\right>_{H(k)}\\
&=2\left<
\mu(P)-\mu(N(m_0,\Sigma_0)),\sqrt{n}\left\{
\frac{1}{n}\sum_{i=1}^{n}k(\cdot,Y_i)-\mu(N(\widehat{m},\widehat{\Sigma}))-(\mu(P)-\mu(N(m_0,\Sigma_0)))
\right\}
\right>_{H(k)}\\
&~~~~~+\frac{1}{\sqrt{n}}\left\|
\sqrt{n}\left\{
\frac{1}{n}\sum_{i=1}^{n}k(\cdot,Y_i)-\mu(N(\widehat{m},\widehat{\Sigma}))-(\mu(P)-\mu(N(m_0,\Sigma_0)))
\right\}
\right\|^2_{H(k)}.
\end{align*}
It follows from direct calculations as given in (\ref{E5}) that
\begin{align*}
&\sqrt{n}\left\{
\frac{1}{n}\sum_{i=1}^{n}k(\cdot,Y_i)-\mu(N(\widehat{m},\widehat{\Sigma}))-(\mu(P)-\mu(N(m_0,\Sigma_0)))
\right\}\\
&=\frac{1}{\sqrt{n}}\sum_{i=1}^{n}k(\cdot,Y_i)-\sqrt{n}\mu(N(m_0,\Sigma_0))\\
&~~~~~-\sqrt{n}D_{(m_0,\Sigma_0)}\mu\left(N\left(\frac{1}{n}\sum_{i=1}^{n}(Y_i-m_0),\frac{1}{n}\sum_{i=1}^{n}((Y_i-m_0)^{\otimes 2} -\Sigma_0)\right)\right)\\
&~~~~~-D_{(m_0,\Sigma_0)}\mu((0,\sqrt{n}(\widehat{m}-m_0)^{\otimes 2})) -\sqrt{n}R_2((m_0,\Sigma_0),(\widehat{m}-m_0,\widehat{\Sigma}-\Sigma_0)\\
&~~~~~-\sqrt{n}(\mu(P)-\mu(N(m_0,\Sigma_0)))\\
&=\frac{1}{\sqrt{n}}\sum_{i=1}^{n}\left\{f(Y_i)-(\mu(P)-\mu(N(m_0,\Sigma_0)))\right\}\\
&~~~~~-D_{(m_0,\Sigma_0)}\mu((0,\sqrt{n}(\widehat{m}-m_0)^{\otimes 2})) -\sqrt{n}R_2((m_0,\Sigma_0),(\widehat{m}-m_0,\widehat{\Sigma}-\Sigma_0).
\end{align*}
Furthermore, since $\mathbb{E}\left[f(Y_1)\right]
	=\mu(P) -\mu(N(m_0,\Sigma_0))$, 
so that 
\[
\sqrt{n}(\widehat{\Delta}^2-\Delta^2) \xrightarrow{\mathcal{D}} 2\left<\mu(P)-\mu(N(m_0,\Sigma_0)),X\right>_{H(k)}
\]
by the central limit theorem (see \cite{Hoffmann and Pisier}), where $X \sim N(0,V[f(Y_1)])$.
This finally gives
\begin{align*}
2\left<\mu(P)-\mu(N(m_0,\Sigma_0)),X\right>_{H(k)}\sim N\left(0,v^2\right). \qed
\end{align*}
\begin{Rem}\label{Rem3}
	Define that
	\begin{align*}
	\Sigma_k&=V[k(\cdot,Y_1)];~~H(k) \to H(k),\\
	\Sigma(m_0,\Sigma_0)&=V[(Y_1-m_0,(Y_1-m_0)^{\otimes 2} -\Sigma_0)];~~\mathcal{H} \oplus HS(\mathcal{H}) \to \mathcal{H} \oplus HS(\mathcal{H}),\\
	\Sigma_k(m_0,\Sigma_0) &=\mathbb{E}\left[
	(Y_1-m_0,(Y_1-m_0)^{\otimes 2} -\Sigma_0) \otimes k(\cdot,X)\right];~~H(k) \to \mathcal{H} \oplus HS(\mathcal{H})
	\end{align*}
	and let $D_{(m_0,\Sigma_0)}\mu(N)^*(\cdot)$ be  the adjoint operator of $D_{(m_0,\Sigma_0)}\mu(N(\cdot,\cdot))$ i.e. it is a linear operator (see Section I.1 of \cite{Reed and Simon} for details) from $H(k)$ to $\mathcal{H} \oplus HS(\mathcal{H})$ such that for all $g \in H(k),~(h,A) \in \mathcal{H} \oplus HS(\mathcal{H})$,
	\[
	\left<g,D_{(m_0,\Sigma_0)}\mu(N(m_0,\Sigma_0))\right>_{H(k)}=\left<D_{(m_0,\Sigma_0)}\mu(N)^*(g),(h,A)\right>_{\mathcal{H}\oplus HS(\mathcal{H})}.
	\]	
	Using above operators, $V[f(Y_1)]$ can be calculated as follows.
	For all $h,h' \in H(k)$,
	\begin{align*}
	&\left<V[f(Y_1)]h,h'\right>_{H(k)}\\
	&=\mathbb{E}\left[\left<f(Y_1)-\mathbb{E}\left[f(Y_1)\right],h\right>_{H(k)}\left<f(Y_1)-\mathbb{E}\left[f(Y_1)\right],h'\right>_{H(k)}\right]\\
	&=\left<(\Sigma_k-D_{(m_0,\Sigma_0)}\mu(N(\Sigma_k(m_0,\Sigma_0)))-D_{(m_0,\Sigma_0)}\mu(N(\Sigma_k(m_0,\Sigma_0)))^*\right.\\
	&~~~~~\left.+D_{(m_0,\Sigma_0)}\mu(N(\Sigma(m_0,\Sigma_0)))D_{(m_0,\Sigma_0)}\mu(N)^*)h,h'\right>_{H(k)},
	\end{align*}
	where 
	\[
	D_{(m_0,\Sigma_0)}\mu(N(\Sigma_k(m_0,\Sigma_0)))^*=(\Sigma_k(m_0,\Sigma_0))^*D_{(m_0,\Sigma_0)}\mu(N)^*.
	\]
	Therefore, 
	\begin{align*}
	V[f(Y_1)]&=\Sigma_k-D_{(m_0,\Sigma_0)}\mu(N(\Sigma_k(m_0,\Sigma_0)))-D_{(m_0,\Sigma_0)}\mu(N(\Sigma_k(m_0,\Sigma_0)))^*\\
	&~~~~~+D_{(m_0,\Sigma_0)}\mu(N(\Sigma(m_0,\Sigma_0)))D_{(m_0,\Sigma_0)}\mu(N)^*
	\end{align*}
	and $V[f(Y_1)]$ exists by assumption.
\end{Rem}
\subsection{Proof of Theorem \ref{T3}}
It follows from direct calculations as given in (\ref{E5}) that
\begin{align*}
n\widehat{\Delta}^2
&=\left\|
\frac{1}{\sqrt{n}}\sum_{i=1}^{n}k(\cdot,Y_i) -\sqrt{n}\mu(N(\widehat{m},\widehat{\Sigma}))
\right\|^2_{H(k)}\\
&=\left\|
\frac{1}{\sqrt{n}}\sum_{i=1}^{n}(k(\cdot,Y_i)
-D_{(m_0,\Sigma_0)} \mu(N(Y_i-m_0,(Y_i-m_0)^{\otimes 2} -\Sigma_0))
-\mu(N(m_0,\Sigma_0))
)
\right.\\
&\left.~~~~~-\sqrt{n}D_{(m_0,\Sigma_0)}\mu((0,(\widehat{m}-m_0)^{\otimes 2}))-\sqrt{n}R_2((m_0,\Sigma_0),(\widehat{m}-m_0,\widehat{\Sigma}-\Sigma_0))
\right\|^2_{H(k)}\\
&=\left\|\frac{1}{\sqrt{n}}\sum_{i=1}^{n}f(Y_i)-\sqrt{n}D_{(m_0,\Sigma_0)}\mu((0,(\widehat{m}-m_0)^{\otimes 2}))-\sqrt{n}R_2((m_0,\Sigma_0),(\widehat{m}-m_0,\widehat{\Sigma}-\Sigma_0))\right\|^2_{H(k)}.
\end{align*}
We see from Lemma \ref{L5} that $ \left\|D_{(m_0,\Sigma_0)}\mu((0,(\widehat{m}-m_0)^{\otimes 2}))\right\|_{H(k)}=O_p(n^{-1})$
and it also follows by Lemma \ref{L6} that $\displaystyle \left\|R_2((m_0,\Sigma_0),(\widehat{m}-m_0,\widehat{\Sigma}-\Sigma_0))\right\|_{H(k)}=O_p(n^{-1}).$
Hence
\begin{align*}
n\widehat{\Delta}^2
&=\frac{1}{n}\sum_{i,j=1}^{n}\left<f(Y_i),f(Y_j)\right>_{H(k)} +O_p\left(\frac{1}{\sqrt{n}}\right).
\end{align*}
We here utilize an operator $S_k$ in (\ref{E6}), which is a Hilbert-Schmidt operator, see Theorem VI.22 in \cite{Reed and Simon}.
Furthermore, $S_k$ is a self adjoint Hilbert-Schmidt operator.
Therefore 
\[
\left<f(x),f(y)\right>_{H(k)}=\sum_{\ell=1}^{\infty} \lambda_{\ell}\Psi_{\ell} (x) \Psi_{\ell}(y) 
\]
by Theorem 1 in \cite{Minh}, where for all $\ell \in \mathbb{N}$, $\lambda_{\ell}$ is eigenvalue of $S_k$ and $\Psi_{\ell}$ is eigenfunction  corresponding to $\lambda_{\ell}$, each satisfies (\ref{E30}) and (\ref{E31}).
Also, $\displaystyle \sum_{\ell=1}^{\infty}\lambda^2_{\ell} < \infty$ since $S_k$ is Hilbert-Schmidt.
The expected value of $f$ is $\displaystyle \mathbb{E}[f(Y_1)]=\mu(P)-\mu(N(m_0,\Sigma_0))=(1/\sqrt{n})\eta(Q).$
Therefore,
\begin{align*}
\lambda_{\ell} \mathbb{E}[\Psi_{\ell}(Y_1)]
=\int_{\mathcal{H}} \left<\mathbb{E}[f(Y_1)] ,f(y)\right>_{H(k)} \Psi_{\ell}(y) dN(m_0,\Sigma_0) (y)
=\frac{1}{\sqrt{n}}\eta_{\ell}(Q),
\end{align*}
where $\displaystyle \eta_{\ell}(Q) =\int_{\mathcal{H}} \left<\eta(Q), f(y)\right>_{H(k)} \Psi_{\ell} (y) dN(m_0,\Sigma_0) (y)$.
Hence $\displaystyle \mathbb{E}[\Psi_{\ell}(Y_1)] =\eta_{\ell}(Q)/(\sqrt{n} \lambda_{\ell}).$
Further, since $Y_i$ and $Y_j~(i \neq j)$ are independent and $\Psi_{\ell} \in L^2(\mathcal{H} , N(m_0,\Sigma_0))$, $\Psi_{\ell}(Y_i)$ and $\Psi_{\ell}(Y_j)$ are independent.
Therefore
\begin{align*}
V[\Psi_{\ell}(Y_1)] &=\mathbb{E}[\Psi_{\ell}(Y_1)^2] -\left\{
\mathbb{E}[\Psi_{\ell}(Y_1)] 
\right\}^2\\
&=\int_{\mathcal{H}} \Psi_{\ell} (y)^2dN(m_0,\Sigma_0) (y) +\frac{1}{\sqrt{n}}\int_{\mathcal{H}} \Psi_{\ell} (y) ^2d(Q-N(m_0,\Sigma_0)) (y) -\frac{1}{n} \cdot \frac{\eta^2_{\ell}(Q)}{\lambda^2_{\ell}}\\
&=1+\frac{1}{\sqrt{n}}\tau_{\ell \ell} -\frac{1}{n}\cdot \frac{\eta^2_{\ell}(Q)}{\lambda^2_{\ell}},
\end{align*}
where $\displaystyle \tau_{\ell s} =\int_{\mathcal{H}}\Psi_{\ell} (y) \Psi_s(y)d(Q-N(m_0,\Sigma_0))(y)$.  
Since the distribution to be obtained is 
\begin{align*}
n\widehat{\Delta}^2
&=\frac{1}{n}\sum_{i,j=1}^{n} \left<f(Y_i),f(Y_j)\right>_{H(k)} +O_p\left(\frac{1}{\sqrt{n}}\right)\\
&=\sum_{\ell=1}^{\infty} \lambda_{\ell} \left(\frac{1}{\sqrt{n}}\sum_{i=1}^{n}\Psi_{\ell} (Y_i)\right)^2+O_p\left(\frac{1}{\sqrt{n}}\right),
\end{align*}
 first we aim to obtain asymptotic distribution of $ (1/\sqrt{n})\sum_{i=1}^{n}\Psi_{\ell}(Y_i)$.
 Let $\psi_n$ be the characteristic function $ (1/\sqrt{n})\sum_{i=1}^{n} \Psi_{\ell}(Y_i) $ and $\varphi$ be the characteristic function $\displaystyle \Psi_{\ell} (Y_1)$.
 Then 
 \begin{align*}
 \psi_n(t)
 &=\mathbb{E}\left[
 \exp\left(it\frac{1}{\sqrt{n}}\sum_{i=1}^{n}\Psi_{\ell} (Y_i) \right)
 \right]\\
 &=\varphi\left(\frac{t}{\sqrt{n}}\right)^n\\
 &=\left\{ 
 1+i\left(\frac{t}{\sqrt{n}}\right)\mathbb{E}[\Psi_{\ell}(Y_1)]+\frac{1}{2}i^2\left(\frac{t}{\sqrt{n}}\right)^2\mathbb{E}\left[\Psi_{\ell}(Y_1)^2\right]+o\left(\frac{1}{n}\right)\right\}^n\\
 &=\left\{
 1+\frac{1}{n}\left\{
 it\cdot\frac{\eta_{\ell}(Q)}{\lambda_{\ell}}-\frac{t^2}{n}\left\{1+\frac{1}{\sqrt{n}}\tau_{\ell\ell}\right\}
 \right\}+o\left(\frac{1}{n}\right)
 \right\}^n
 \to \exp\left(it\cdot\frac{\eta_{\ell}(Q)}{\lambda_{\ell}}-\frac{t^2}{2}\right)
 \end{align*}
 as $n \to \infty$.
 Hence $(1/\sqrt{n})\sum_{i=1}^{n}\Psi_{\ell}(Y_i) \xrightarrow{\mathcal{D}}N\left((\eta_{\ell}(Q)/\lambda_{\ell}),1\right)$.
 Next, we obtain asymptotic distribution of 
 $
 (1/\sqrt{n})\sum_{i=1}^{n} 
 \begin{bmatrix}
 \Psi_{\ell}(Y_i)\\ 
 \Psi_{s} (Y_i)
 \end{bmatrix}
 ~(\ell \neq s) 
  $.
 Since 
 \[
 \mathbb{E}
 \begin{bmatrix}
 \Psi_{\ell}(Y_i)\\ 
 \Psi_{s}(Y_i)
 \end{bmatrix}
 =\frac{1}{\sqrt{n}}
 \begin{bmatrix}
 \displaystyle  \frac{\eta_{\ell}(Q)}{\lambda_{\ell}}\\ 
 \displaystyle \frac{\eta_{s}(Q)}{\lambda_{s}}
 \end{bmatrix},
 \]
 \begin{align*}
 V
 \begin{bmatrix}
 \Psi_{\ell}(Y_i)\\ 
 \Psi_{s}(Y_i)
 \end{bmatrix} 
 &=
 \begin{bmatrix}
 \displaystyle 1+\frac{1}{\sqrt{n}}\tau_{\ell\ell}-\frac{1}{n}\cdot\frac{\eta^2_{\ell}(Q)}{\lambda^2_{\ell}}& \displaystyle -\frac{1}{\sqrt{n}}\tau_{\ell s}-\frac{1}{n}\cdot\frac{\eta_{\ell}(Q)\eta_{s}(Q)}{\lambda_{\ell}\lambda_s}  \\ 
 \displaystyle -\frac{1}{\sqrt{n}}\tau_{\ell s}-\frac{1}{n}\cdot\frac{\eta_{\ell}(Q)\eta_{s}(Q)}{\lambda_{\ell}\lambda_s}& \displaystyle 1+\frac{1}{\sqrt{n}}\tau_{ss}-\frac{1}{n}\cdot\frac{\eta^2_{s}(Q)}{\lambda^2_{s}}
 \end{bmatrix}\\
 &=I+\frac{1}{\sqrt{n}}
 \begin{bmatrix}
 \tau_{\ell\ell}& \tau_{\ell s} \\ 
 \tau_{\ell s}& \tau_{ss} 
 \end{bmatrix} 
 -\frac{1}{n}
 \begin{bmatrix}
 \displaystyle \frac{\eta^2_{\ell}(Q)}{\lambda_{\ell}^2}& \displaystyle \frac{\eta_{\ell}(Q) \eta_{s}(Q)}{\lambda_{\ell s}} \\ 
 \displaystyle \frac{\eta_{\ell}(Q) \eta_{s}(Q)}{\lambda_{\ell} \lambda_{s}}& \displaystyle \frac{\eta^2_{s}(Q)}{\lambda_{s}^2} 
 \end{bmatrix} ,
 \end{align*}
 the characteristic function of $(1/\sqrt{n})\sum_{i=1}^{n}
 \begin{bmatrix}
 \Psi_{\ell}(Y_i)\\ 
 \Psi_{s}(Y_i)
 \end{bmatrix}$ can be evaluated as
 \begin{align*}
 &\mathbb{E}\left[
 \exp\left(i\underline{t}^T\frac{1}{\sqrt{n}}\sum_{i=1}^{n}
 \begin{bmatrix}
 \Psi_{\ell}(Y_i)\\ 
 \Psi_s(Y_i)
 \end{bmatrix} 
 \right)
 \right]\\
 &=\left\{
 1+i\frac{\underline{t}^T}{n}
 \begin{bmatrix}
 \displaystyle \frac{\eta_{\ell}(Q)}{\lambda_{\ell}}\\ 
 \displaystyle \frac{\eta_{s}(Q)}{\lambda_{s}}
 \end{bmatrix} 
 -\frac{1}{2}\frac{1}{n}\underline{t}^T
 \left\{
 I+\frac{1}{\sqrt{n}}
 \begin{bmatrix}
 \tau_{\ell \ell}& \tau_{\ell s} \\ 
 \tau_{\ell s}& \tau_{ss}
 \end{bmatrix} 
 \right\}\underline{t}+o\left(\frac{1}{n}\right)
 \right\}^n\\
 &\to \exp\left(i\underline{t}^T 
 \begin{bmatrix}
 \displaystyle \frac{\eta_{\ell}(Q)}{\lambda_{\ell}}\\ 
 \displaystyle \frac{\eta_{s}(Q)}{\lambda_{s}}
 \end{bmatrix} 
 -\frac{1}{2}\underline{t}^T\underline{t}
 \right)
 \end{align*}
 as $n \to \infty$.
 Hence
 \[
 \frac{1}{\sqrt{n}}\sum_{i=1}^{n}
 \begin{bmatrix}
 \Psi_{\ell}(Y_i)\\ 
 \Psi_{s}(Y_i)
 \end{bmatrix} \xrightarrow{\mathcal{D}} N\left(
 \begin{bmatrix}
 \displaystyle \frac{\eta_{\ell}(Q)}{\lambda_{\ell}}\\ 
 \displaystyle \frac{\eta_{s}(Q)}{\lambda_{s}}
 \end{bmatrix} ,I
 \right).
 \]
 Therefore, according to the same argument, for any $L \in \mathbb{N}$,
\begin{equation}\label{E25}
 \frac{1}{\sqrt{n}}\sum_{i=1}^{n}
\begin{bmatrix}
\Psi_{1}(Y_i)\\ 
\vdots\\
\Psi_{L}(Y_i)
\end{bmatrix} 
\xrightarrow{\mathcal{D}}
N\left(
\begin{bmatrix}
\displaystyle
\frac{\eta_1(Q)}{\lambda_1}\\ 
\vdots\\
\displaystyle
\frac{\eta_L(Q)}{\lambda_L}
\end{bmatrix} 
,I
\right).
\end{equation}
 Let 
 \begin{align*}
 &X_{L_{n}}=\sum_{\ell=1}^{L}\lambda_{\ell} \left(\dfrac{1}{\sqrt{n}}\Psi_{\ell}(Y_i)\right)^2,&&&&X_{L}=\sum_{\ell=1}^{L}\lambda_{\ell} W^2_{\ell},\\
 &Y_n=\sum_{\ell=1}^{\infty} \lambda_{\ell} \left(\dfrac{1}{\sqrt{n}}\sum_{i=1}^{n}\Psi_{\ell} (Y_i)\right)^2,&&\text{and}&&X=\sum_{\ell=1}^{\infty} \lambda_{\ell} W^2_{\ell}
 \end{align*}
 where $W_{\ell}~(\ell=1,2,\dots)$ are independent and $\displaystyle W_{\ell} \sim N\left((\eta_{\ell}(Q)/\lambda_{\ell}),1\right)$.
 Then $X_{L_{n}} \xrightarrow{ \mathcal{D}} X_{L}$ by (\ref{E25}).
 Also,
 \begin{align*}
 &\mathbb{E}\left[
 \left|
 \sum_{\ell=1}^{L}\lambda_{\ell} W^2_{\ell} -\sum_{\ell=1}^{\infty} \lambda_{\ell} W^2_{\ell}
 \right|^2
 \right]\\
 &=\sum_{\ell=L+1}^{\infty} \lambda^2_{\ell} \mathbb{E}\left[W^4_{\ell}\right]
 +\sum_{\ell \neq s} \lambda_{\ell} \lambda_{s} \mathbb{E}\left[W^2_{\ell}\right]\mathbb{E}\left[W^2_{s}\right]\\
&=2\sum_{\ell=L+1}^{\infty} \lambda^2_{\ell} +4\sum_{\ell=L+1}^{\infty} \eta^2_{\ell}(Q) +\left(\sum_{\ell=L+1}^{\infty} \lambda_{\ell}\right)^2
+2\left(
\sum_{\ell=L+1}^{\infty}\lambda_{\ell}
\right)
\left(\sum_{\ell=L+1}^{\infty}\frac{\eta^2_{\ell}(Q)}{\lambda_{\ell}}\right)
+\left(
\sum_{\ell=L+1}^{\infty} \frac{\eta^2_{\ell}(Q)}{\lambda_{\ell}}
\right)^2.
\end{align*}
Furthermore,
\[
\sum_{\ell=1}^{\infty} \frac{\eta^2_{\ell}(Q)}{\lambda_{\ell}}< \infty,
\]
by the definition of $\mathcal{A}_k$ and $\lim_{\ell \to 0}\lambda_{\ell}=0$ by $\sum_{\ell=1}^{\infty}\lambda_{\ell} < \infty$.
Hence, there exists $s \in \mathbb{N}$ such that for any $n > s$, $\lambda_n <1$.
Therefore, since
\[
\sum_{\ell=1}^{\infty} \eta^2_{\ell}(Q) =\sum_{\ell=1}^{s}\eta^2_{\ell}(Q) +\sum_{\ell=s+1}^{\infty} \eta^2_{\ell}(Q) =\sum_{\ell=1}^{s} \eta^2_{\ell}(Q) +\sum_{\ell=s+1}^{\infty} \frac{\eta^2_{\ell}(Q)}{\lambda_{\ell}} < \infty,
\]
we have 
\begin{equation}\label{E26}
\mathbb{E}\left[
\left|
\sum_{\ell=1}^{L} \lambda_{\ell} W^2_{\ell} -\sum_{\ell=1}^{\infty}\lambda_{\ell} W^2_{\ell} 
\right|^2
\right] \to 0
\end{equation}
as $L \to \infty$.
From (\ref{E26}), we get $X_{L} \xrightarrow{\mathcal{D}} X$.
Next, for any $\varepsilon > 0$,
\begin{align}
&\lim_{L \to \infty} \limsup_{n \to \infty} \mathbb{P}(|X_{L_{n}} -Y_n| > \varepsilon) \nonumber\\
& \leq \frac{1}{\varepsilon} \lim_{L \to \infty} \limsup_{n \to \infty} \mathbb{E}[|X_{L_{n}}-Y_n|]\nonumber\\
&=\frac{1}{\varepsilon}\lim_{L \to \infty} \limsup_{n \to \infty} \sum_{\ell=L+1}^{\infty} \lambda_{\ell} \mathbb{E}\left[
\frac{1}{n}\sum_{i,j=1}^{n}\Psi_{\ell}(Y_i) \Psi_{\ell}(Y_j)
\right]\label{E18}\\
&=\frac{1}{\varepsilon}\lim_{L \to \infty} \limsup_{n \to \infty} \sum_{\ell=L+1}^{\infty} \lambda_{\ell} \frac{1}{n} \left\{
\sum_{i=1}^{n}\mathbb{E}[\Psi_{\ell}(Y_i)^2]+\sum_{i \neq j} \mathbb{E}[\Psi_{\ell}(Y_i)]\mathbb{E}[\Psi_{\ell}(Y_j)]
\right\}\nonumber\\
&=\frac{1}{\varepsilon}\lim_{L \to \infty} \limsup_{n \to \infty} \sum_{\ell=L+1}^{\infty} \lambda_{\ell} \{
\mathbb{E}[\Psi_{\ell}(Y_1)^2]+(n-1)(\mathbb{E}[\Psi_{\ell}(Y_1)])^2
\}\nonumber\\
&=\frac{1}{\varepsilon}\lim_{L \to \infty} \limsup_{n \to \infty} \sum_{\ell=L+1}^{\infty} \lambda_{\ell}\left\{
1+\frac{\eta^2_{\ell}(Q)}{\lambda^2_{\ell}}+\frac{1}{\sqrt{n}}\tau_{\ell\ell} -\frac{1}{n}\frac{\eta^2_{\ell}(Q)}{\lambda^2_{\ell}}
\right\}\nonumber\\
&=\frac{1}{\varepsilon}\lim_{L \to \infty} \sum_{\ell=L+1}^{\infty} \left\{
\lambda_{\ell}+\frac{\eta^2_{\ell}(Q)}{\lambda_{\ell}}
\right\}\nonumber\\
&=0\nonumber
\end{align}
from Markov's inequality.
In (\ref{E18}), we have used the dominated convergence theorem because 
\[
\sum_{\ell=L+1}^{\infty} \lambda_{\ell} \left(\frac{1}{\sqrt{n}}\sum_{i=1}^{\infty}\Psi_{\ell}(Y_i)\right)^2 \leq \sum_{\ell=1}^{\infty}\left(
\frac{1}{\sqrt{n}}\sum_{i=1}^{\infty} \Psi_{\ell}(Y_i)
\right)^2=\frac{1}{n}\sum_{i=1}^{n} \sum_{j=1}^{n} \left<f(Y_i),f(Y_j)\right>_{H(k)}
\]
and
\begin{align*}
&\mathbb{E}\left[
\frac{1}{n}\sum_{i=1}^{n} \sum_{j=1}^{n} \left<f(Y_i),f(Y_j)\right>_{H(k)}
\right]\\
&=\frac{1}{n}\sum_{i=1}^{n}\mathbb{E}\left[
\left\|f(Y_i)\right\|^2_{H(k)}
\right]+\frac{1}{n}\sum_{i \neq j} \mathbb{E}\left[
\left<f(Y_i),f(Y_j)\right>_{H(k)}
\right]\\
&=\mathbb{E}\left[
\left\|f(Y_1)\right\|^2_{H(k)}
\right]
+\left(1-\frac{1}{n}\right)\left\|\eta(Q)\right\|^2_{H(k)}\\
&< \infty.
\end{align*}
Therefore, $Y_n \xrightarrow{\mathcal{D}} X$ by Theorem 4.2 in \cite{Billingsley}. \qed
\subsection{Proof of Corollary \ref{C1}}
For the case $Q=N(m_0,\Sigma_0)$ in Thereom \ref{T3}, we  can obtain asymptotic null distribution of $n\widehat{\Delta}^2$ under $H_0: P=N(m_0,\Sigma_0)$, since $\eta_{\ell}(Q)=0$ for all $\ell \in \mathbb{N}$.
Therefore, under $H_0$ as $n \to \infty$,
\[
n\widehat{\Delta}^2 \overset{\mathcal{D}}{\to} \sum_{\ell=1}^{\infty} \lambda_{\ell}Z^2_{\ell},
\]
whre $Z_{\ell} \overset{i.i.d.}{\sim} N(0,1)$.
\qed
\subsection{Proof of Lemma \ref{L5}}
For any $\delta >0$, there exsits $N \in \mathbb{N}$ such that for all $n > N$,
\begin{align*}
&\left(1-\frac{1}{n}\right)\left\{\left\|A(0,\Sigma_0)\right\|^2_{H(k)}+\mathbb{E}\left[\left\|A(0,(Y_j-m_0) \otimes (Y_i-m_0))\right\|^2_{H(k)}\right]\right.\\
&~~~~~\left.+\mathbb{E}\left[
\left<A(0,(Y_j-m_0) \otimes (Y_i-m_0)), A(0,(Y_i-m_0) \otimes (Y_j-m_0))\right>_{H(k)}
\right]
\right\}\\
&~~~~~+\frac{1}{n}\mathbb{E}\left[
\left\|
A(0,(Y_1-m_0)^{\otimes 2})
\right\|^2_{H(k)}
\right]\\
&< \left\|A(0,\Sigma_0)\right\|^2_{H(k)}+\mathbb{E}\left[\left\|A(0,(Y_j-m_0) \otimes (Y_i-m_0))\right\|^2_{H(k)}\right]\\
&~~~~~+\mathbb{E}\left[
\left<A(0,(Y_j-m_0) \otimes (Y_i-m_0)), A(0,(Y_i-m_0) \otimes (Y_j-m_0))\right>_{H(k)}
\right]+\delta.
\end{align*}
Let $M_{\delta} \in \mathbb{R}$ be such that 
\begin{align*}
M_{\delta} &> \frac{1}{\sqrt{\delta}}\left(
\left\|A(0,\Sigma_0)\right\|^2_{H(k)}+\mathbb{E}\left[\left\|A(0,(Y_j-m_0) \otimes (Y_i-m_0))\right\|^2_{H(k)}\right] \right.\\
&~~~~~\left.
+\mathbb{E}\left[
\left<A(0,(Y_j-m_0) \otimes (Y_i-m_0)), A(0,(Y_i-m_0) \otimes (Y_j-m_0))\right>_{H(k)}
\right]+\delta
\right)^{1/2}.
\end{align*}
By Theorem 2.1 in \cite{Rao}, for all $n > N$,
\begin{align*}
&\mathbb{P}\left(\left\|A(0,n(\widehat{m}-m_0)^{\otimes 2})\right\|_{H(k)}>M_{\delta}\right)\\
& \leq \frac{\mathbb{E}\left[\left\|A(0,n(\widehat{m}-m_0)^{\otimes 2})\right\|^2_{H(k)}\right]}{M_{\delta}^2}\\
&=\frac{\mathbb{E}\left[\left<A(0,n(\widehat{m}-m_0)^{\otimes 2}),A(0,n(\widehat{m}-m_0)^{\otimes 2})\right>_{H(k)}\right]}{M_{\delta}^2}\\
&=\frac{\sum_{i,j,s,\ell=1}^{n}\mathbb{E}\left[\left<A(0,\left<Y_i-m_0,\cdot\right>_{H(k)}(Y_j-m_0)),A(0,\left<Y_s-m_0,\cdot\right>_{H(k)}(Y_\ell-m_0))\right>_{H(k)}\right]}{n^2M_{\delta}^2}.
\end{align*}
We put $\mu_{ijs\ell} =\mathbb{E}\left[\left<A(0,\left<Y_i-m_0,\cdot\right>_{H(k)}(Y_j-m_0)),A(0,\left<Y_s-m_0,\cdot\right>_{H(k)}(Y_\ell-m_0))\right>_{H(k)}\right]$, and we aim to evaluate $\sum_{i,j,s,\ell=1}^{\infty} \mu_{ijs\ell}$.

If $i \neq j ,s, \ell$, we have 
\begin{align}
\mu_{ijs\ell}
=\mathbb{E}\left[\left<A(0,\left<\mathbb{E}[Y_i-m_0],\cdot\right>_{H(k)}(Y_j-m_0)),A(0,\left<Y_s-m_0,\cdot\right>_{H(k)}(Y_\ell-m_0))\right>_{H(k)}\right] 
=0 \label{E19}.
\end{align}
Similarly,
\begin{equation}\label{E20}
\mu_{ijs\ell}=0 
\end{equation}
if case $j \neq i, s, \ell$, case $s \neq i, j, \ell$ or case $ \ell \neq i,j,s$.
For the cases if $i=j,~s=\ell,~i \neq s$, we have
\begin{align}
\mu_{iiss}
=\left<A(0,\Sigma_0),A(0,\Sigma_0)\right>_{H(k)}
=\left\|A(0,\Sigma_0)\right\|^2_{H(k)} \label{E21}.
\end{align}
Also, if $i=s,~j=\ell,~i \neq j$, we obtain
\begin{align}
\mu_{ijij} 
=\mathbb{E}\left[
\left\|
A(0,\left<Y_i-m_0,\cdot\right>_{H(k)}(Y_j-m_0))\right\|^2_{H(k)}
\right]\label{E22}.
\end{align}
Similarly, if $i=\ell,~j=s,~i \neq j$, we get
\begin{align}
\mu_{ijji}
=\mathbb{E}\left[
\left<A(0,(Y_j-m_0) \otimes (Y_i-m_0)), A(0,(Y_i-m_0)  \otimes (Y_j-m_0))\right>_{H(k)}
\right]. \label{E23}
\end{align}
Finally, for the cases of $i=j=s=\ell$, we see that 
\begin{align}
\mu_{iiii} 
=\mathbb{E}\left[
\left\|A(0,(Y_i-m_0)^{\otimes 2})\right\|^2_{H(k)}
\right]\label{E24}.
\end{align}
All calculations (\ref{E19}), (\ref{E20}), (\ref{E21}), (\ref{E22}), (\ref{E22}), (\ref{E23}) and (\ref{E24}) furnish to reach
\begin{align*}
&\frac{\sum_{i,j,s,\ell=1}^{n}\mathbb{E}\left[\left<A(0,\left<Y_i-m_0,\cdot\right>_{H(k)}(Y_j-m_0)),A(0,\left<Y_s-m_0,\cdot\right>_{H(k)}(Y_\ell-m_0))\right>_{H(k)}\right]}{n^2M_{\delta}^2}\\
&=\frac{1}{M^2_{\delta}}\left(
\left(1-\frac{1}{n}\right)\left(\left\|A(0,\Sigma_0)\right\|^2_{H(k)} +\mu_{1212}+\mu_{1221}\right)+\frac{1}{n}\mathbb{E}\left[
\left\|A(0,(Y_i-m_0)^{\otimes 2})\right\|^2_{H(k)}
\right]
\right)\\
&< \frac{1}{M^2_{\delta}}\left(\left\|A(0,\Sigma_0)\right\|^2_{H(k)} +\mu_{1212}+\mu_{1221}+\delta\right)\\
&<\delta,
\end{align*}
which means $\displaystyle \left\|A(0,n(\widehat{m}-m_0)^{\otimes 2})\right\|_{H(k)}=O_p\left(n^{-1}\right)$.\qed
\subsection{Proof of Lemma \ref{L6}}
Since $\mu(N(\cdot,\cdot)) \in C^2(B((m_0,\Sigma_0),\varepsilon) ,H(k))$, for all $(m',\Sigma') \in B((m_0,\Sigma_0) ,\varepsilon)$, 
\[
\lim_{(m,\Sigma) \to (m',\Sigma')} \left\|D^2_{(m,\Sigma)} \mu(N(\cdot,\cdot))-D^2_{(m',\Sigma')} \mu(N(\cdot,\cdot))\right\|=0,
\] 
where $\|\cdot\|$ is the operator norm (see Section I.1 of \cite{Reed and Simon} for details).
Hence
\[
\lim_{(m,\Sigma) \to (m',\Sigma') } \sup_{x \neq 0} \frac{\left\|(D^2_{(m,\Sigma)} \mu(N(x)) ) -(D^2_{(m',\Sigma')}\mu(N(x)))\right\|_{H(k)}}{\left\|x\right\|_{(\mathcal{H} \oplus HS(\mathcal{H}))^2}}=0,
\]
i.e. for all $x \in (\mathcal{H} \oplus HS(\mathcal{H}))^2 \backslash \{0\}$,
\[
\lim_{(m,\Sigma) \to (m',\Sigma') }  \frac{\left\|(D^2_{(m,\Sigma)} \mu(N(x))) -(D^2_{(m',\Sigma')}\mu(N(x)))\right\|_{H(k)}}{\left\|x\right\|_{(\mathcal{H} \oplus HS(\mathcal{H}))^2}}=0.
\]
Therefore, for all $x \in (\mathcal{H} \oplus HS(\mathcal{H}))^2$,
\[
\lim_{(m,\Sigma) \to (m',\Sigma')} \left\|(D^2_{(m,\Sigma)}\mu(N(x))) -(D^2_{(m',\Sigma')}\mu(N(x))) \right\|_{H(k)}=0.
\]
This leads that
\begin{align*}
&\lim_{s \to s'} (D^2_{(m_0,\Sigma_0) +s(\widehat{m}-m_0,\widehat{\Sigma}-\Sigma_0)}\mu(N(\sqrt{n}(\widehat{m}-m_0),\sqrt{n}(\widehat{\Sigma}-\Sigma_0))^2))\\
&=(D^2_{(m_0,\Sigma_0) +s'(\widehat{m}-m_0,\widehat{\Sigma}-\Sigma_0)}\mu(N(\sqrt{n}(\widehat{m}-m_0),\sqrt{n}(\widehat{\Sigma}-\Sigma_0))^2)).
\end{align*}
Here we put
\[
K(n) =\max_{s \in [0,1]} \left\|
(D^2_{(m_0,\Sigma_0) +s(\widehat{m}-m_0,\widehat{\Sigma}-\Sigma_0)}\mu(N(\sqrt{n}(\widehat{m}-m_0),\sqrt{n}(\widehat{\Sigma}-\Sigma_0))^2))
\right\|_{H(k)}.
\]
Since
\begin{align*}
K(n) 
&=\max_{s \in [0,1]} \left\|
(D^2_{(m_0,\Sigma_0) +s(\widehat{m}-m_0,\widehat{\Sigma}-\Sigma_0)}\mu(N(\sqrt{n}(\widehat{m}-m_0),\sqrt{n}(\widehat{\Sigma}-\Sigma_0))^2))
\right\|_{H(k)}\\
&=2\left\|(\sqrt{n}(\widehat{m}-m_0),\sqrt{n}(\widehat{\Sigma}-\Sigma_0))\right\|_{\mathcal{H} \oplus HS(\mathcal{H})} \\
&~~~~~\times \max_{s \in [0,1]} 
\frac{\left\|(D^2_{(m_0,\Sigma_0) +s(\widehat{m}-m_0,\widehat{\Sigma}-\Sigma_0)}\mu(N(\sqrt{n}(\widehat{m}-m_0),\sqrt{n}(\widehat{\Sigma}-\Sigma_0))^2))\right\|_{H(k)}}{\left\|(\sqrt{n}(\widehat{m}-m_0),\sqrt{n}(\widehat{\Sigma}-\Sigma_0))^2\right\|_{(\mathcal{H} \oplus HS(\mathcal{H}))^2}}\\
& \leq  2\left\|(\sqrt{n}(\widehat{m}-m_0),\sqrt{n}(\widehat{\Sigma}-\Sigma_0))\right\|_{\mathcal{H} \oplus HS(\mathcal{H})} \\
&~~~~~\times \max_{s \in [0,1]} \sup_{((x,y),(x',y')) \neq 0}
\frac{\left\|(D^2_{(m_0,\Sigma_0) +s(\widehat{m}-m_0,\widehat{\Sigma}-\Sigma_0)}\mu(N((x,y),(x',y'))\right\|_{H(k)}}{\left\|((x,y),(x',y')\right\|_{(\mathcal{H} \oplus HS(\mathcal{H}))^2}}\\
&=2\left\|(\sqrt{n}(\widehat{m}-m_0),\sqrt{n}(\widehat{\Sigma}-\Sigma_0))\right\|_{\mathcal{H} \oplus HS(\mathcal{H})} \max_{s \in [0,1]} \left\|(D^2_{(m_0,\Sigma_0) +s(\widehat{m}-m_0,\widehat{\Sigma}-\Sigma_0)}\mu(N(\cdot,\cdot))\right\|\\
&=O_p(1),
\end{align*}
it follows that
\begin{align}
&\left\|
\int^{1}_{0} (1-s) D^2_{(m_0,\Sigma_0) +s(\widehat{m}-m_0,\widehat{\Sigma}-\Sigma_0)}\mu(N(\widehat{m}-m_0,\widehat{\Sigma}-\Sigma_0)^2) ds
\right\|^2_{H(k)}\nonumber\\
& \leq \frac{1}{n^2} \left(
\int^{1}_{0} (1-s)\left\|
 D^2_{(m_0,\Sigma_0) +s(\widehat{m}-m_0,\widehat{\Sigma}-\Sigma_0)}\mu(N(\sqrt{n}(\widehat{m}-m_0),\sqrt{n}(\widehat{\Sigma}-\Sigma_0))^2) 
\right\|_{H(k)}ds
\right)^2\nonumber\\
&\leq \frac{1}{n^2} \left(
\int^{1}_{0} (1-s) \max_{s \in [0,1]}\left\|
D^2_{(m_0,\Sigma_0) +s(\widehat{m}-m_0,\widehat{\Sigma}-\Sigma_0)}\mu(N(\sqrt{n}(\widehat{m}-m_0),\sqrt{n}(\widehat{\Sigma}-\Sigma_0))^2) 
\right\|_{H(k)}ds
\right)^2\nonumber\\
&=\frac{1}{n^2}(K(n))^2 \left(\int^{1}_{0} (1-s) ds\right)^2\nonumber\\
&=O_p\left(\frac{1}{n^2}\right). \label{E27}
\end{align}
From (\ref{E27}), we finally have
\[
\left\|
\int^{1}_{0} (1-s) D^2_{(m_0,\Sigma_0) +s(\widehat{m}-m_0,\widehat{\Sigma}-\Sigma_0)}\mu(N(\widehat{m}-m_0,\widehat{\Sigma}-\Sigma_0)^2) ds
\right\|_{H(k)}=O_p\left(\frac{1}{n}\right).
\]\qed
\subsection{Proof of (\ref{E32})}
A direct calculation gives
\begin{align*}
&\frac{\partial}{\partial \underline{m}} \mu(N(\underline{m},\Sigma)) \\
&=|I_d+2\sigma\Sigma|^{-1/2}
\frac{\partial}{\partial \underline{m}} \exp\left(-\sigma(\cdot-\underline{m})^T(I_d+2\sigma\Sigma)^{-1} (\cdot-\underline{m})\right)\\
&=2\sigma|I_d+2\sigma\Sigma|^{-1/2}\exp\left(-\sigma(\cdot-\underline{m})^T(I_d+2\sigma\Sigma)^{-1} (\cdot-\underline{m})\right)(I_d+2\sigma\Sigma)^{-1}(\cdot-\underline{m})\\
&=2\sigma\mu(N(\underline{m},\Sigma))(I_d+2\sigma\Sigma)^{-1}(\cdot-\underline{m}).
\end{align*}
\qed
\subsection{Proof of (\ref{E33})}
A straightforward but lengthy computation yields that
\begin{align*}
&\frac{\partial}{\partial \sigma_{ij}}\mu(N(\underline{m},\Sigma))\\
&=\frac{\partial}{\partial \sigma_{ij}}
|I_d+2\sigma\Sigma|^{-1/2} \exp\left(-\sigma(\cdot-\underline{m})^T(I_d+2\sigma\Sigma)^{-1} (\cdot-\underline{m})\right)\\
&=-\frac{1}{2}|I_d+2\sigma\Sigma|^{-3/2} \exp\left(-\sigma(\cdot-\underline{m})^T(I_d+2\sigma\Sigma)^{-1}(\cdot-\underline{m})\right)\frac{\partial}{\partial \sigma_{ij}}|I_d+2\sigma\Sigma|\\
&~~~~~-\sigma|I_d+2\sigma\Sigma|^{-1/2} \exp\left(-\sigma(\cdot-\underline{m})^T(I_d+2\sigma\Sigma)^{-1}(\cdot-\underline{m})\right)\\
&~~~~~\times
\frac{\partial}{\partial\sigma_{ij}}(\cdot-\underline{m})^T(I_d+2\sigma\Sigma)^{-1}(\cdot-\underline{m})\\
&=\sigma \mu(N(\underline{m},\Sigma)) \text{tr}\left(
\frac{\partial \Sigma}{\partial \sigma_{ij}}  (I_d+2\sigma\Sigma)^{-1} \left(
2\sigma (\cdot-\underline{m})(\cdot-\underline{m})^T  (I_d+2\sigma\Sigma)^{-1}-I_d
\right)
\right).
\end{align*}
\qed

\subsection{Proof of (\ref{E34})}
We shall try to obtain the Fourier transform of each term of (\ref{E40}).
Note that
\begin{align}
&\frac{1}{(2\pi)^d}\int_{\mathbb{R}^d} f(\underline{x})(\underline{z})\exp(-i\underline{t}^T\underline{z})d\underline{z}\nonumber\\
&=\frac{1}{(2\pi)^d}\int_{\mathbb{R}^d} k(\underline{x},\underline{z}) \exp(-i\underline{t}^T\underline{z}) d\underline{z}\nonumber\\
&~~~~~-\frac{1}{(2\pi)^d}\int_{\mathbb{R}^d} \mu(N(\underline{m}_0,\Sigma_0))(\underline{z}) \exp(-i\underline{t}^T\underline{z}) d\underline{z}\nonumber\\
&~~~~~-\frac{2\sigma}{(2\pi)^d} \int_{\mathbb{R}^d} \mu(N(\underline{m}_0,\Sigma_0))(\underline{z})(\underline{z}-\underline{m}_0)^TV^{-1} (\underline{x}-\underline{m}_0)\exp(-i\underline{t}^T\underline{z}) d\underline{z} \nonumber \\
&~~~~~-\frac{\sigma}{(2\pi)^d}\int_{\mathbb{R}^d} \mu(N(\underline{m}_0,\Sigma_0))(\underline{z}) A(\underline{m}_0,\Sigma_0)(\underline{z})^T\text{vech}\left((\underline{x}-\underline{m}_0)(\underline{x}-\underline{m}_0)^T-\Sigma_0\right) \exp(-i\underline{t}^T\underline{z}) d\underline{z}\nonumber\\
&\equiv I_1-I_2-I_3-I_4. \label{E48}
\end{align}
The first term $I_1$ is easily obtained as
\begin{align}
I_1
=\frac{1}{(2\pi)^d}\int_{\mathbb{R}^d} \exp\left(-\sigma\|\underline{z}-\underline{x}\|^2_{\mathbb{R}^d}\right) \exp(-i\underline{t}^T\underline{z}) d\underline{z}=\frac{1}{\sqrt{(4\pi\sigma)^d}} \exp\left(-i\underline{x}^T\underline{t}-\frac{1}{4\sigma}\underline{t}^T\underline{t}\right) \label{E41}.
\end{align}
The term $I_2$ can be verified as
\begin{align}
I_2
&=\frac{1}{(2\pi)^d}|V|^{-1/2} \int_{\mathbb{R}^d} \exp(-\sigma(\underline{z}-\underline{m}_0)^TV^{-1}(\underline{z}-\underline{m}_0))\exp(-i\underline{t}^T\underline{z}) d\underline{z}\nonumber\\
&=\frac{1}{\sqrt{(4\pi\sigma)^d}}\exp\left(-i\underline{m}_0^T\underline{t}-\frac{1}{4\sigma}\underline{t}^TV\underline{t}\right).\label{E42}
\end{align}
The term $I_3$ can be calculated that 
\begin{align}
I_3
&=\frac{2\sigma}{(2\pi)^d} |V|^{-1/2} \int_{\mathbb{R}^d} \exp\left(-\sigma(\underline{z}-\underline{m}_0)^T V^{-1} (\underline{z}-\underline{m}_0)\right) (\underline{z}-\underline{m}_0)^T V^{-1} (\underline{x}-\underline{m}_0) \exp(-i\underline{t}^T\underline{z}) d\underline{z} \nonumber\\
&=\frac{1}{(2\pi)^d}\cdot\frac{1}{\sqrt{(2\sigma)^{d-1}}} \int_{\mathbb{R}^d} \exp\left(-\frac{1}{2}\underline{u}^T\underline{u}\right) \underline{u}^TV^{-1/2} (\underline{x}-\underline{m}_0) 
\exp\left(-i\underline{t}^T\left(\frac{1}{\sqrt{2\sigma}}V^{1/2}\underline{u}+\underline{m}_0\right)\right)d\underline{u}\nonumber\\
&=- \frac{i}{\sqrt{(4\pi\sigma)^d}} \exp\left(-i\underline{t}^T\underline{m}_0-\frac{1}{4\sigma}\underline{t}^TV\underline{t}\right)\underline{t}^T(\underline{x}-\underline{m}_0),\label{E43}
\end{align}
where we have used the change of variables $\underline{u}=\sqrt{2\sigma} V^{-1/2}(\underline{z}-\underline{m}_0)$.
Necessary computations to get $I_4$ are
\begin{align}
I_4
&=\frac{\sigma}{(2\pi)^d} |V|^{-1/2} \int_{\mathbb{R}^d} \exp(-\sigma(\underline{z}-\underline{m}_0)^TV^{-1}(\underline{z}-\underline{m}_0))\nonumber\\
&~~~~~\times\text{tr}\left[
V^{-1} \{2\sigma(\underline{z}-\underline{m}_0)(\underline{z}-\underline{m}_0)^TV^{-1}-I_d\} B(\underline{x})
\right]
\exp(-i\underline{t}^T\underline{z}) d\underline{z}\nonumber\\
&=\frac{2\sigma^2}{(2\pi)^d}|V|^{-1/2}\nonumber \\
&~~~~~\times\text{tr}\left[
V^{-1}\int_{\mathbb{R}^d} \exp(-\sigma(\underline{z}-\underline{m}_0)^TV^{-1}(\underline{z}-\underline{m}_0)) (\underline{z}-\underline{m}_0)(\underline{z}-\underline{m}_0)^T \exp(-i\underline{t}^T\underline{z})d\underline{z}V^{-1}B(\underline{x})
\right] \nonumber\\
&~~~~~-\frac{\sigma}{(2\pi)^d}|V|^{-1/2} \text{tr}[V^{-1}B(\underline{x})] \int_{\mathbb{R}^d} \exp(-\sigma(\underline{z}-\underline{m}_0)^TV^{-1}(\underline{z}-\underline{m}_0)) \exp(-i\underline{t}^T\underline{z}) d\underline{z}\label{E44}.
\end{align}
Here the integral inside of trace is
\begin{align}
&\int_{\mathbb{R}^d} \exp(-\sigma(\underline{z}-\underline{m}_0)^TV^{-1}(\underline{z}-\underline{m}_0)) (\underline{z}-\underline{m}_0)(\underline{z}-\underline{m}_0)^T \exp(-i\underline{t}^T\underline{z})d\underline{z}\nonumber\\
&=\frac{1}{\sqrt{(2\sigma)^{d+2}}} |V|^{1/2} \int_{\mathbb{R}^d} \exp\left(
-\frac{1}{2}\underline{u}^T\underline{u} 
\right)
V^{1/2} \underline{u}\underline{u}^T V^{1/2} \exp\left(-i\underline{t}^T
\left(
\frac{1}{\sqrt{2\sigma}}V^{1/2} \underline{u}+\underline{m}_0
\right)
\right) d\underline{u}\nonumber\\
&=\frac{1}{(2\sigma)^2} \sqrt{\frac{\pi^d}{\sigma^d}} |V|^{1/2} V^{1/2} \exp\left(-i\underline{t}^T\underline{m}_0-\frac{1}{4\sigma}\underline{t}^TV\underline{t}\right) \left(2\sigma I_d-
V^{1/2} \underline{t}\underline{t}^TV^{1/2}
\right) V^{1/2}, \label{E45}
\end{align}
where we have used the change of variables $\underline{u}=\sqrt{2\sigma} V^{-1/2} (\underline{z}-\underline{m}_0)$.
Combining (\ref{E44}) and (\ref{E45}) yields
\begin{align}
I_4=-\frac{1}{2} \cdot
\frac{1}{\sqrt{(4\pi\sigma)^d}} \exp\left(
-i\underline{t}^T\underline{m}_0-\frac{1}{4\sigma}\underline{t}^TV\underline{t}
\right) \underline{t}^TB(\underline{x})\underline{t}\label{E46}.
\end{align}
Hence (\ref{E41}), (\ref{E42}), (\ref{E43}), (\ref{E46}) together with (\ref{E48}), we reach
\begin{align*}
\widehat{(f(\underline{x}))}(\underline{t}) 
&=\frac{1}{\sqrt{(4\pi\sigma)^d}}
\left\{
\exp\left(-i\underline{x}^T\underline{t}-\frac{1}{4\sigma}\underline{t}^T\underline{t}\right)+\exp\left(-i\underline{m}_0^T\underline{t}-\frac{1}{4\sigma}\underline{t}^TV\underline{t}\right)
\right.\\
&~~~~~\left.
\times\left(1
-i\underline{t}^T(\underline{x}-\underline{m}_0)-\frac{1}{2}\underline{t}^T B(\underline{x}) \underline{t}
\right)
\right\}.
\end{align*}
\qed

\subsection{Proof of (\ref{E47})}
We calculate $\left<f(\underline{x})(\cdot), f(\underline{y})(\cdot) \right>_{H(k)}$ using (\ref{E34}) .
Note that
\begin{align}
&\left<f(\underline{x})(\cdot), f(\underline{y})(\cdot) \right>_{H(k)} \nonumber\\
&=\left<k(\cdot,\underline{x}),k(\cdot,\underline{y})\right>_{H(k)}\nonumber\\
&~~~~~-\left<k(\cdot,\underline{x}),\mu(N(\underline{m}_0,\Sigma_0))(\cdot) \{
1+2\sigma(\cdot-\underline{m}_0)^TV^{-1} (\underline{y}-\underline{m}_0)
+\sigma A(\underline{m}_0,\Sigma_0)(\cdot)^T \text{vech}(B(\underline{y}))\} \right>_{H(k)}\nonumber\\
&~~~~~-\left<k(\cdot,\underline{y}),\mu(N(\underline{m}_0,\Sigma_0))(\cdot) \{
1+2\sigma(\cdot-\underline{m}_0)^TV^{-1} (\underline{x}-\underline{m}_0)
+\sigma A(\underline{m}_0,\Sigma_0)(\cdot)^T \text{vech}(B(\underline{x}))\} \right>_{H(k)}\nonumber\\
&~~~~~+\frac{1}{\sqrt{(4\pi\sigma)^d}} \int_{\mathbb{R}^d} \exp\left(-\frac{1}{2\sigma}\underline{t}^T\left(V-\frac{1}{2}I_d\right)\underline{t}\right)
\left\{1 -i\underline{t}^T(\underline{x}-\underline{m}_0)-\frac{1}{2}\underline{t}^TB(\underline{x}) \underline{t} \right\} \nonumber\\
&~~~~~~~~~~\times
\left\{1 +i\underline{t}^T(\underline{y}-\underline{m}_0)-\frac{1}{2}\underline{t}^TB(\underline{y}) \underline{t} \right\} d\underline{t}\nonumber\\
&\equiv J_1-J_2-J_3+J_4. \label{E49}
\end{align}
Our focus goes to $J_4$, since it includes a bit messy calculations.
We see that
\begin{align}
J_4
&=\frac{1}{\sqrt{(4\pi\sigma)^d}} \int_{\mathbb{R}^d} \exp\left(-\frac{1}{4\sigma}\underline{t}^T\left(2V- I_d\right)\underline{t}\right)\biggl\{
1-i\underline{t}^T(\underline{x}-\underline{y})\nonumber\\
&~~~~~+\underline{t}^T\left((\underline{x}-\underline{m}_0)(\underline{y}-\underline{m}_0)^T-\frac{1}{2}B(\underline{x})-\frac{1}{2}B(\underline{y})\right)\underline{t}\nonumber\\
&~~~~~+\frac{i}{2}\underline{t}^T\left(
(\underline{x}-\underline{m}_0)\underline{t}^TB(\underline{y})
-(\underline{y}-\underline{m}_0)\underline{t}^TB(\underline{x})
\right)\underline{t}+\frac{1}{4}\underline{t}^TB(\underline{x})\underline{t}\underline{t}^TB(\underline{y})\underline{t}\biggr\}d\underline{t}\nonumber\\
&\equiv J_{41} -J_{42} +J_{43} +J_{44} +J_{45}. \label{E50}
\end{align}
The term $J_{41}$ is easily obtained as
\begin{align}
J_{41}
&=\frac{1}{\sqrt{(4\pi\sigma)^d}} \int_{\mathbb{R}^d} \exp\left(-\frac{1}{4\sigma}\underline{t}^T\left(2V- I_d\right)\underline{t}\right) d\underline{t}= |2V- I_d|^{-1/2}  \label{E35}.
\end{align}
By the form of density function of normal distribution with mean $0$, the term $J_{42}$ is in fact
\begin{align}
J_{42}
&=\frac{i}{\sqrt{(4\pi\sigma)^d}}\int_{\mathbb{R}^d} \underline{t}^T((\underline{x}-\underline{m}_0)-(\underline{y}-\underline{m}_0))\exp\left(-\frac{1}{4\sigma}\underline{t}^T\left(2V- I_d\right)\underline{t}\right) d\underline{t}=0\label{E36}.
\end{align}
The term $J_{43}$ can be calculated by using Theorem 9.18 in \cite{Schott} as
\begin{align}
J_{43}&=\frac{1}{\sqrt{(4\pi\sigma)^d}} \int_{\mathbb{R}^d}
\underline{t}^T\left((\underline{x}-\underline{m}_0)(\underline{y}-\underline{m}_0)^T-\frac{1}{2}B(\underline{x})-\frac{1}{2}B(\underline{y})\right)\underline{t}\nonumber\exp\left(-\frac{1}{4\sigma}\underline{t}^T\left(2V- I_d\right)\underline{t}\right) d\underline{t}\nonumber\\
&=2\sigma |2V- I_d|^{-1/2} \text{tr}\left[
\left((\underline{x}-\underline{m}_0)(\underline{y}-\underline{m}_0)^T-\frac{1}{2}B(\underline{x})-\frac{1}{2}B(\underline{y})\right)(2V- I_d)^{-1}
\right].\label{E37}
\end{align}
From Section 2.6.2 in \cite{Anderson}, we have
\begin{align}
J_{44}
&=\frac{1}{\sqrt{(4\pi\sigma)^d}}\int_{\mathbb{R}^d} \underline{t}^T\left(
(\underline{x}-\underline{m}_0)\underline{t}^TB(\underline{y})
-(\underline{y}-\underline{m}_0)\underline{t}^TB(\underline{x})
\right)\underline{t}\exp\left(-\frac{1}{4\sigma}\underline{t}^T\left(2V- I_d\right)\underline{t}\right) d\underline{t}=0.\label{E38}
\end{align}
Theorem 9.21 in \cite{Schott} yields that
\begin{align}
J_{45}
&=\frac{1}{4\sqrt{(4\pi\sigma)^d}} \int_{\mathbb{R}^d} \underline{t}^TB(\underline{x}) \underline{t} \underline{t}^TB(\underline{y}) \underline{t}\exp\left(-\frac{1}{4\sigma}\underline{t}^T\left(2V- I_d\right)\underline{t}\right) d\underline{t} \nonumber\\
&=\sigma^2 |2V- I_d|^{-1/2} \Bigl\{
\text{tr}\left[B(\underline{x})(2V- I_d)^{-1}\right]
\text{tr}\left[B(\underline{y})(2V- I_d)^{-1}\right] \nonumber\\
&~~~~~ +2\text{tr}\left[ B(\underline{x})(2V- I_d)^{-1} B(\underline{y}) (2V- I_d)^{-1}\right]
\Bigr\}.\label{E39}
\end{align}
By combining (\ref{E50}), (\ref{E35}), (\ref{E36}), (\ref{E37}), (\ref{E38}) and (\ref{E39}), we finally have
\begin{align}
J_{4}&=|2V- I_d|^{-1/2} \Bigl\{
1+\sigma \text{tr}\left[
\Bigl(2(\underline{x}-\underline{m}_0)(\underline{y}-\underline{m}_0)^T-B(\underline{x})-B(\underline{y})\Bigr)(2V- I_d)^{-1}
\right]\nonumber
\\
&~~~~~
+\sigma^2\Bigl(\text{tr}\left[B(\underline{x}) (2V- I_d)^{-1}\right] \text{tr}\left[B(\underline{y}) (2V- I_d)^{-1}\right]\nonumber\\
&~~~~~
+2\text{tr}\left[B(\underline{x}) (2V- I_d)^{-1} B(\underline{y}) (2V- I_d)^{-1}\right]
\Bigr)
\Bigr\}. \label{E51}
\end{align}
Also, $J_1$ is obtained from the property of the reproducing kernel as
\begin{equation}\label{E52}
J_1=\left<k(\cdot,\underline{x}),k(\cdot,\underline{y})\right>_{H(k)} =\exp\left(-\sigma \|\underline{x}-\underline{y}\|^2_{\mathbb{R}^d}\right).
\end{equation}
The term $J_2$ is having the form
\begin{align}
J_2
&=\mu(N(\underline{m}_0,\Sigma_0))(\underline{x}) \Bigl\{
1+2\sigma(\underline{x}-\underline{m}_0)^TV^{-1} (\underline{y}-\underline{m}_0)
+\sigma A(\underline{m}_0,\Sigma_0)(\underline{x})^T \text{vech}\left(B(\underline{y})\right)\Bigr\}\nonumber\\
&=|V|^{-1/2} \exp\left(-\sigma(\underline{x}-\underline{m}_0)^TV^{-1}(\underline{x}-\underline{m}_0) \right)
\Bigl\{
1+2\sigma(\underline{x}-\underline{m}_0)^TV^{-1} (\underline{y}-\underline{m}_0) \nonumber\\
&~~~~~+\sigma \text{tr}\left[
V^{-1}\{2\sigma(B(\underline{x}) +\Sigma_0) V^{-1} -I_d \}B(\underline{y})\right] \Bigr\}, \label{E53}
\end{align}
and symmetric calculation gives $J_3$ as
\begin{align}
J_{3}
&=|V|^{-1/2} \exp\left(-\sigma(\underline{y}-\underline{m}_0)^TV^{-1}(\underline{y}-\underline{m}_0) \right)
\Bigl\{
1+2\sigma(\underline{y}-\underline{m}_0)^TV^{-1} (\underline{x}-\underline{m}_0)\nonumber \\
&~~~~~+\sigma \text{tr}\left[
V^{-1}\{2\sigma(B(\underline{y}) +\Sigma_0) V^{-1} -I_d \}B(\underline{x})\right] \Bigr\}.\label{E54}
\end{align}
(\ref{E51}), (\ref{E52}), (\ref{E53}) and (\ref{E54}) furnish to reach (\ref{E47}).
\qed

\subsection{Proof of Proposition \ref{P1} }
In the sequent discussions, we repeat use the results for expectation of multitple quadratic forms summarized in Section 9.6 of \cite{Schott} and the formula
\begin{eqnarray}
\phi_{\Sigma_{1}}(\underline{x}-\underline{m}_{1})\phi_{\Sigma_{2}}(\underline{x}-\underline{m}_{2})=\phi_{\Sigma_{1}+\Sigma_{2}}(\underline{m}_{1}-\underline{m}_{2})\phi_{(\Sigma_{1}^{-1}+\Sigma_{2}^{-1})^{-1}}(\underline{x}-\underline{m}^{*}), \label{Prop_GaussDensity}
\end{eqnarray}
where 
$$
\underline{m}^{*}=(\Sigma_{1}^{-1}+\Sigma_{2}^{-1})^{-1}(\Sigma_{1}^{-1}\underline{m}_{2}+\Sigma_{2}^{-1}\underline{m}_{1})
$$
and $\phi_{\Sigma}(\cdot-\underline{m})$ designates the density of $N_{d}(\underline{m},\Sigma)$, see e.g. Appendix C in \cite{Wand and Jones}.

First we devide 
\begin{eqnarray} \label{Int_Expectation}
\left<f(\underline{x})(\cdot),f(\underline{x})(\cdot)\right>_{H(k)}=\mathcal{C}-I_1-I_2,
\end{eqnarray}
where
\begin{align*}
\mathcal{C}&=1+|V+2\sigma \Sigma_{0}|^{-1/2}\Big\{
1+2\sigma \text{tr}[(V+2\sigma \Sigma_{0})^{-1}\Sigma_0]+\sigma^2\{\text{tr}[(V+2\sigma \Sigma_{0})^{-1}\Sigma_0]\}^2\\
&\hspace{4cm}+2\sigma^2\text{tr}[\{(V+2\sigma \Sigma_{0})^{-1} \Sigma_{0}\}^2]
\Big\}, \nonumber \\
I_{1}&=2|V|^{-1/2} \exp(-\sigma(\underline{x}-\underline{m}_0)^TV^{-1}(\underline{x}-\underline{m}_0)\Big\{
1+\sigma\text{tr}[V^{-1}\Sigma_0] \nonumber \\
&\hspace{2cm}+\sigma(\underline{x}-\underline{m}_0)^TV^{-1}(I_d-2\sigma \Sigma_0 V^{-1}) (\underline{x}-\underline{m}_0)+2\sigma^2((\underline{x}-\underline{m}_0)^TV^{-1}(\underline{x}-\underline{m}_0))^2\Big\} \\
\ &\equiv I_{11}+I_{12}+I_{13}, \nonumber \\
I_{2}&=\sigma^2|V+2\sigma \Sigma_{0}|^{-1/2}\\
&~~~~~\times\Big\{
2(\underline{x}-\underline{m}_0)^T(V+2\sigma \Sigma_{0})^{-1}(\text{tr}[(V+2\sigma \Sigma_{0})^{-1}\Sigma_0] I_d+2\Sigma_0(V+2\sigma \Sigma_{0})^{-1})(\underline{x}-\underline{m}_0) \nonumber \\
&\hspace{1.5cm}
-3\{
(\underline{x}-\underline{m}_0)^T(V+2\sigma \Sigma_{0})^{-1}(\underline{x}-\underline{m}_0)
\}^2
\Big\} \nonumber \\ 
\ &\equiv I_{21}-I_{22}. 
\end{align*}
And we note the notation (\ref{DEF_MT_V}).

We have by direct computations using results in Section 9.6 of \cite{Schott} that
\begin{align*}
&\int_{\mathbb{R}^d} I_{21} dN(\underline{m}_0,\Sigma_0)(\underline{x})\\
&=2\sigma^2|V+2\sigma \Sigma_{0}|^{-1/2} \\
&~~~~~\times\int_{\mathbb{R}^d} 
(\underline{x}-\underline{m}_0)^T(V+2\sigma \Sigma_{0})^{-1}(\text{tr}[(V+2\sigma \Sigma_{0})^{-1}\Sigma_0] I_d+2\Sigma_0(V+2\sigma \Sigma_{0})^{-1})(\underline{x}-\underline{m}_0) dN(\underline{m}_0,\Sigma_0)(\underline{x})\\
&=2\sigma^2|V+2\sigma \Sigma_{0}|^{-1/2} \Big\{
\{\text{tr}[(V+2\sigma \Sigma_{0})^{-1}\Sigma_0]\}^2+2\text{tr}\left[
\{(V+2\sigma \Sigma_{0})^{-1}\Sigma_{0}\}^2
\right]
\Big\}
\end{align*}
and
\begin{align*}
&\int_{\mathbb{R}^d} I_{22} dN(\underline{m}_0,\Sigma_0)(\underline{x})\\
&=3\sigma^2|V+2\sigma \Sigma_{0}|^{-1/2}\int_{\mathbb{R}^d} \{(\underline{x}-\underline{m}_0)^T(V+2\sigma \Sigma_{0})^{-1}(\underline{x}-\underline{m}_0)\}^2 dN(\underline{m}_0,\Sigma_0)(\underline{x})\\
&=3\sigma^2|V+2\sigma \Sigma_{0}|^{-1/2}\Big\{
\{\text{tr}[(V+2\sigma \Sigma_{0})^{-1}\Sigma_0]\}^2+2\text{tr}\left[
\{(V+2\sigma \Sigma_{0})^{-1}\Sigma_0\}^2
\right]
\Big\},
\end{align*}
from which it follows that
\begin{eqnarray} \label{Int_I2}
\int_{\mathbb{R}^d} I_2 dN(\underline{m}_0,\Sigma_0)(\underline{x})=-\sigma^2|V+2\sigma \Sigma_{0}|^{-1/2}\Big\{
\{\text{tr}[(V+2\sigma \Sigma_{0})^{-1}\Sigma_0]\}^2+2\text{tr}\left[
\{(V+2\sigma \Sigma_{0})^{-1}\Sigma_0\}^2
\right]
\Big\}.
\end{eqnarray}
Next our focus goes to $I_{1}$.
Note that 
\begin{align*}
\exp(-\sigma(\underline{x}-\underline{m}_0)^T V^{-1} (\underline{x}-\underline{m}_0))=\left(\frac{\pi}{\sigma}\right)^{d/2} |V|^{1/2} \phi_{\frac{1}{2\sigma}V}(\underline{x}-\underline{m}_0)
\end{align*}
and (\ref{Prop_GaussDensity}) yields that
\begin{align*}
\phi_{\frac{1}{2\sigma}V}(\underline{x}-\underline{m}_0) \phi_{\Sigma_0}(\underline{x}-\underline{m}_0)=\phi_{\frac{1}{2\sigma}V+\Sigma_0}(\underline{0}) \phi_{(2\sigma V^{-1} +\Sigma_0^{-1})^{-1}}(\underline{x}-\underline{m}^*)
\end{align*}
with
\[
m^*=(2\sigma V^{-1} +\Sigma^{-1}_0)^{-1} (2\sigma V^{-1} \underline{m}_0+\Sigma_0^{-1}\underline{m}_0)=\underline{m}_0.
\]
Further it is easy to check that 
\begin{align*}
(2\sigma V^{-1}+\Sigma^{-1}_0)^{-1}
&=\Sigma_0(2V-I_d)^{-1}V,
\end{align*}
by which we have
\[
\phi_{\frac{1}{2\sigma}V}(\underline{x}-\underline{m}_0) \phi_{\Sigma_0}(\underline{x}-\underline{m}_0)
=\phi_{\frac{1}{2\sigma}(2V-I_d)} (\underline{0}) \phi_{\Sigma_0(2V-I_d)^{-1}V} (\underline{x}-\underline{m}_0).
\]
Repeat use of above equalities gives that
\begin{align*}
&\int_{\mathbb{R}^d} I_{11} dN(\underline{m}_0,\Sigma_0)(\underline{x})\\
&=2|V|^{-1/2} (1+\sigma\text{tr}[V^{-1}\Sigma_0])\int_{\mathbb{R}^d} \exp(-\sigma(\underline{x}-\underline{m}_0)^T V^{-1} (\underline{x}-\underline{m}_0)) dN(\underline{m}_0,\Sigma_0)(\underline{x})\\
&=2|V+2\sigma \Sigma_{0}|^{-1/2} (1+\sigma \text{tr}[V^{-1}\Sigma_0]),
\end{align*}
\begin{align*}
&\int_{\mathbb{R}^d} I_{12} dN(\underline{m}_0,\Sigma_0)(\underline{x})\\
&=2\sigma|V|^{-1/2} \int_{\mathbb{R}^d} (\underline{x}-\underline{m}_0)^TV^{-1}(I_d-2\sigma\Sigma_0V^{-1})(\underline{x}-\underline{m}_0)
\exp\left(
-\sigma(\underline{x}-\underline{m}_0)^TV^{-1}(\underline{x}-\underline{m}_0)
\right)dN(\underline{m}_0,\Sigma_0)(\underline{x})\\
&=2\sigma|V+2\sigma \Sigma_{0}|^{-1/2} \text{tr}\left[
\Sigma_0(V+2\sigma \Sigma_{0})^{-1}V^{-1}
\right]
\end{align*}
and
\begin{align*}
&\int_{\mathbb{R}^d} I_{13} dN(\underline{m}_0,\Sigma_0) (\underline{x})\\
&=4\sigma^2|V|^{-1/2} 
\int_{\mathbb{R}^d} \{
(\underline{x}-\underline{m}_0)^TV^{-1}(\underline{x}-\underline{m}_0)
\}^2
\exp\left(
-\sigma(\underline{x}-\underline{m}_0)^TV^{-1}(\underline{x}-\underline{m}_0)
\right)dN(\underline{m}_0,\Sigma_0)(\underline{x})\\
&=4\sigma^2|V+2\sigma \Sigma_{0}|^{-1/2} \Big\{
\{\text{tr}[(V+2\sigma \Sigma_{0})^{-1}\Sigma_0]\}^2+2\text{tr}[\{(V+2\sigma \Sigma_{0})^{-1}\Sigma_0\}^2]
\Big\}.
\end{align*}
These are combined into
\begin{align}
&\int_{\mathbb{R}^d} I_1dN(\underline{m}_0,\Sigma_0)(\underline{x}) \nonumber \\
&=2|V+2\sigma \Sigma_{0}|^{-1/2} \Big\{
1+\sigma\text{tr}[V^{-1}\Sigma_0] +\sigma \text{tr}[\Sigma_0(V+2\sigma \Sigma_{0})^{-1}V^{-1}]\nonumber \\ 
&~~~~~+2\sigma^2
\{\text{tr}[(V+2\sigma \Sigma_{0})^{-1}\Sigma_0]\}^2+4\sigma^2\text{tr}[\{(V+2\sigma \Sigma_{0})^{-1}\Sigma_0\}^2]
\Big\} \label{Int_I1}
\end{align}
We finally obtain Proposition \ref{P1} by (\ref{Int_Expectation}), (\ref{Int_I2}) and (\ref{Int_I1}) and the fact taht
\begin{align*}
\text{tr}[V^{-1}\Sigma_0]+\text{tr}[\Sigma_0( V+2\sigma \Sigma_{0})^{-1}V^{-1}]=2\text{tr}[(V+2\sigma \Sigma_{0})^{-1}\Sigma_0]. \qed
\end{align*}

\subsection{Proof of Proposition \ref{P2} }

We now need to introduce a certain function $\mathcal{F}$ determined by 12 $d \times d$ matrices $T_{k}(k=1,...,6)$, $T'_{k}(k=1,...,6)$, 14 real values $C_{k}(k=0,1,...,6)$ and $C'_{k}(k=0,1,...,6)$ and 6 binary variables $k_{1},k_{2},k_{3}$ and $k'_{1},k'_{2},k'_{3}$, defined by
\begin{align}
&\mathcal{F}(k_1,k_2,k_3|C_0,\dots,C_6|T_1,\dots,T_6\ ||\ k'_1,k'_2,k'_3|C'_0,\dots,C'_6|T'_1,\dots,T'_6) \nonumber\\
&=|R|^{-1/2} C_0C'_0 \Big\{
C_1C'_1
+C_1 Q_1(R^{-1/2} S'_1 R^{-1/2})
+C'_1 Q_1(R^{-1/2} S_1 R^{-1/2}) \nonumber \\
&~~~~~
+C_1C'_3 Q_2(R^{-1/2} S'_2 R^{-1/2},R^{-1/2} S'_2 R^{-1/2})
+C'_1C_3 Q_2(R^{-1/2} S_2 R^{-1/2},R^{-1/2} S_2 R^{-1/2}) \nonumber \\
&~~~~~
+C_1C'_6 Q_2(R^{-1/2} S'_3 R^{-1/2},R^{-1/2} S'_4 R^{-1/2})
+C'_1C_6 Q_2(R^{-1/2} S_3 R^{-1/2},R^{-1/2} S_4 R^{-1/2})\nonumber \\
&~~~~~
+Q_2(R^{-1/2} S_1 R^{-1/2},R^{-1/2} S'_1 R^{-1/2})\nonumber \\
&~~~~~
+C'_3 Q_3(R^{-1/2} S_1 R^{-1/2},R^{-1/2} S'_2 R^{-1/2},R^{-1/2} S'_2 R^{-1/2}) \nonumber \\
&~~~~~
+C_3 Q_3(R^{-1/2} S'_1 R^{-1/2},R^{-1/2} S_2 R^{-1/2},R^{-1/2} S_2 R^{-1/2}) \nonumber \\
&~~~~~
+C'_6Q_3(R^{-1/2} S_1 R^{-1/2},R^{-1/2} S'_3 R^{-1/2},R^{-1/2} S'_4 R^{-1/2}) \nonumber \\
&~~~~~
+C_6Q_3(R^{-1/2} S'_1 R^{-1/2},R^{-1/2} S_3 R^{-1/2},R^{-1/2} S_4 R^{-1/2}) \nonumber \\
&~~~~~
+C_3C'_3 Q_4(R^{-1/2} S_2 R^{-1/2},R^{-1/2} S_2 R^{-1/2},R^{-1/2} S'_2 R^{-1/2},R^{-1/2} S'_2 R^{-1/2}) \nonumber \\
&~~~~~
+C_3C'_6  Q_4(R^{-1/2} S_2 R^{-1/2},R^{-1/2} S_2 R^{-1/2},R^{-1/2} S'_3 R^{-1/2},R^{-1/2} S'_4 R^{-1/2}) \nonumber \\
&~~~~~
+C'_3C_6  Q_4(R^{-1/2} S'_2 R^{-1/2},R^{-1/2} S'_2 R^{-1/2},R^{-1/2} S_3 R^{-1/2},R^{-1/2} S_4 R^{-1/2}) \nonumber \\
&~~~~~+C_6C'_6 Q_4(R^{-1/2} S_3 R^{-1/2},R^{-1/2} S_4 R^{-1/2},R^{-1/2} S'_3 R^{-1/2},R^{-1/2} S'_4 R^{-1/2})
\Big\}, \label{MATHCAL_F_EXP}
\end{align}
where
\begin{align}
&R \nonumber \\
&=I_{2d} +2(k_1+k'_1) 
\begin{bmatrix}
\Sigma_0& -\Sigma_0 \\ 
-\Sigma_0& \Sigma_0 
\end{bmatrix} 
+
\begin{bmatrix}
2(k_2+k'_2) \sigma \Sigma_0^{1/2} V^{-1} \Sigma_0^{1/2}& O \\ 
O& 2(k_3+k'_3) \sigma \Sigma_0^{1/2} V^{-1} \Sigma_0^{1/2} 
\end{bmatrix} \nonumber \\
&=I_{2d} +2(k_1+k'_1) \begin{bmatrix}
1& -1 \\ 
-1& 1
\end{bmatrix} \otimes \Sigma_0
+ 
\begin{bmatrix}
k_2+k'_2& 0  \\ 
0& k_3+k'_3 
\end{bmatrix} \otimes 2\sigma \Sigma_0^{1/2} V^{-1} \Sigma_0^{1/2}	\label{DEF_MAT_R}
\end{align}	
and
\begin{eqnarray}
S_1&=&\begin{bmatrix}
C_4 \Sigma_0^{1/2} T_3 \Sigma_0^{1/2}& \frac{1}{2} C_2 \Sigma_0^{1/2} T_1 \Sigma_0^{1/2}  \\ 
\frac{1}{2} C_2 \Sigma_0^{1/2} T_1 \Sigma_0^{1/2} & C_5 \Sigma_0^{1/2} T_4 \Sigma_0^{1/2}
\end{bmatrix}, \label{DEF_MAT_S1}\\
S_2&=&\begin{bmatrix}
O& \frac{1}{2} \Sigma_0^{1/2} T_2 \Sigma_0^{1/2}  \\ 
\frac{1}{2}\Sigma_0^{1/2} T_2 \Sigma_0^{1/2} & O
\end{bmatrix}, \label{DEF_MAT_S2}\\
S_3&=&\begin{bmatrix}
\Sigma_0^{1/2} T_5 \Sigma_0^{1/2}& O  \\ 
O & O
\end{bmatrix}, \label{DEF_MAT_S3}\\ 
S_4&=&
\begin{bmatrix}
O& O  \\ 
O & \Sigma_0^{1/2} T_6 \Sigma_0^{1/2} \label{DEF_MAT_S4}
\end{bmatrix}
\end{eqnarray}
and $S'_{k}(k=1,2,3,4)$ are the same as $S_{k}(k=1,2,3,4)$ but with $T'_{k}(k=1,...,6)$ and $C'_k(k=0,...,6)$ instead of $T_{k}(k=1,...,6)$ and $C_k(k=0,...,6)$. 
Using these quantities, direct but long calculations furnish to reach the following expression:
\begin{lemm}\label{L7}
	\begin{align}
	&\frac{1}{2}V[Z]=\mathcal{F}_{11}+\mathcal{F}_{22}+\mathcal{F}_{33}+\mathcal{F}_{44}-2\mathcal{F}_{12}-2\mathcal{F}_{13}+2\mathcal{F}_{14}+2\mathcal{F}_{23}-2\mathcal{F}_{24}-2\mathcal{F}_{34} \label{VZ_FULL}
	\end{align}
	where
	\begin{align}
	&\mathcal{F}_{ij} \nonumber\\
	&=\mathcal{F}(k^{(i)}_1,k^{(i)}_2,k^{(i)}_3|C^{(i)}_0,\dots,C^{(i)}_6|T^{(i)}_1,\dots,T^{(i)}_6\ ||\ k^{(j)}_1,k^{(j)}_2,k^{(j)}_3|C^{(j)}_0,\dots,C^{(j)}_6|T^{(j)}_1,\dots,T^{(j)}_6) \label{MATHCAL_F_ij}
	\end{align} 
	for $i,j=1,2,3,4$ and
	\begin{align*}
		(k^{(1)}_{1},k^{(1)}_{2},k^{(1)}_{3})&=(1,0,0),\ \ (k^{(2)}_{1},k^{(2)}_{2},k^{(2)}_{3})=(0 ,1,0), \\
		(k^{(3)}_{1},k^{(3)}_{2},k^{(3)}_{3})&=(0,0,1) ,\ \ (k^{(4)}_{1},k^{(4)}_{2},k^{(4)}_{3})=(0,0,0), 
	\end{align*}
	\begin{align*}
		(C^{(1)}_{0},C^{(1)}_{1},C^{(1)}_{2},C^{(1)}_{3},C^{(1)}_{4},C^{(1)}_{5},C^{(1)}_{6})&=(1,1,0,0,0,0,0), \\
		(C^{(2)}_{0},C^{(2)}_{1},C^{(2)}_{2},C^{(2)}_{3},C^{(2)}_{4},C^{(2)}_{5},C^{(2)}_{6})&=(|V|^{-1/2},1+\sigma \text{tr}[V^{-1} \Sigma_0] ,2\sigma,2\sigma^2,-2\sigma^2,-\sigma,0), \\
		(C^{(3)}_{0},C^{(3)}_{1},C^{(3)}_{2},C^{(3)}_{3},C^{(3)}_{4},C^{(3)}_{5},C^{(3)}_{6})&=(|V|^{-1/2},1+\sigma \text{tr}[V^{-1} \Sigma_0] ,2\sigma,2\sigma^2,-\sigma,-2\sigma^2,0),\\
		(C^{(4)}_{0},C^{(4)}_{2},C^{(4)}_{3},C^{(4)}_{4},C^{(4)}_{5},C^{(4)}_{6})&=(|2V-I_d|^{-1/2},2\sigma,2\sigma^2,-\sigma,-\sigma,\sigma^2),
	\end{align*}
	\[
	C^{(4)}_{1}=1+2\sigma \text{tr}[\Sigma_0(2V-I_d)^{-1}] +\sigma^2 \{\text{tr}[\Sigma_0(2V-I_d)^{-1} ]\}^2 +2\sigma^2 \text{tr}[\{\Sigma_0(2V-I_d)^{-1}\}^2],
	\]
	$T^{(1)}_{1}=T^{(1)}_{2}=T^{(1)}_{3}=T^{(1)}_{4}=T^{(1)}_{5}=T^{(1)}_{6}=O$,
	\begin{align*}
		(T^{(2)}_{1},T^{(2)}_{2},T^{(2)}_{3},T^{(2)}_{4},T^{(2)}_{5},T^{(2)}_{6})&=(V^{-1},V^{-1} ,V^{-1}\Sigma_0V^{-1} ,V^{-1} ,O,O),\\
		(T^{(3)}_{1},T^{(3)}_{2},T^{(3)}_{3},T^{(3)}_{4},T^{(3)}_{5},T^{(3)}_{6})&=(V^{-1},V^{-1} ,V^{-1},V^{-1}\Sigma_0 V^{-1} ,O,O),
	\end{align*}
	$T^{(4)}_{1}=T^{(4)}_{2}=T^{(4)}_{5}=T^{(4)}_{6}=(2V-I_d)^{-1}$ and
	\[
	T^{(4)}_{3}=T^{(4)}_{4}=\Big(
	(1+\sigma \text{tr}[\Sigma_0(2V-I_d)^{-1}])I_d+2\sigma(2V-I_d)^{-1} \Sigma_0
	\Big)(2V-I_d)^{-1}.
	\]
\end{lemm}

The expressions in Lemma \ref{L7} look complicated, but it can be drastically reduced as follows:

 \begin{lemm}\label{L8}
The following equalities hold:
\begin{align}
&\mathcal{F}_{22}=\mathcal{F}_{33} \label{J22=J33}\\
&\mathcal{F}_{12}=\mathcal{F}_{13} \label{J12=J13}\\
&\mathcal{F}_{34}=\mathcal{F}_{24}\label{J34=J24}\\
&\mathcal{F}_{22}=\mathcal{F}_{12}\label{J22=J12}\\
&\mathcal{F}_{44}=\mathcal{F}_{14}\label{J44=J14}\\
&\mathcal{F}_{44}=\mathcal{F}_{23}\label{J44=J23}\\
&\mathcal{F}_{14}=\mathcal{F}_{24}\label{J14=J24}.
\end{align}
Therefore, $V[Z]$ has the expression
 	\begin{align}
 	&\frac{1}{2}V[Z]=\mathcal{F}_{11}-2\mathcal{F}_{22}+\mathcal{F}_{44}. \label{VZ_SIMPLE}
\end{align}
\end{lemm}

Though it needs a bit long calculations of matrices, we can obtain $\mathcal{F}_{11}$, $\mathcal{F}_{22}$ and $\mathcal{F}_{44}$ by almost same manner as addressed in the proof of Lemma \ref{L8}.
Especially we get
\begin{eqnarray}
\mathcal{F}_{11}=|R|^{-1/2}=|I_d+8\sigma \Sigma_0|^{-1/2},\label{MATHCAL_F11}
\end{eqnarray}
\begin{align}
&\mathcal{F}_{22} \nonumber\\
&=|V|^{-1/2} |V+4\sigma\Sigma_0|^{-1/2}  \Big\{
1
+\frac{1}{2}\sigma^2\{\text{tr}[V^{-1} \Sigma_0]\}^2
+\sigma^2 \text{tr}[\{V^{-1} \Sigma_0\}^2]
+\frac{1}{2}\sigma^2\{\text{tr}[(V+4\sigma \Sigma_0)^{-1} \Sigma_0]\}^2 \nonumber\\
&~~~~~~~~~~~~~~~~~~~~~~~~~~~~~~~~
+\sigma^2 \text{tr}[\{(V+4\sigma \Sigma_0)^{-1} \Sigma_0\}^2]
+\sigma\text{tr}[V^{-1} \Sigma_0]
-\sigma \text{tr}[(V+4\sigma \Sigma_0)^{-1} \Sigma_0] \nonumber \\
&~~~~~~~~~~~~~~~~~~~~~~~~~~~~~~~~
-\sigma^2 \text{tr}[V^{-1} \Sigma_0] \text{tr}[(V+4\sigma \Sigma_0)^{-1} \Sigma_0]
\Big\} \label{MATHCAL_F22}
\end{align}
and
\begin{align}
&\mathcal{F}_{44} \nonumber \\
&=|V+2\sigma \Sigma_0|^{-1} \Big\{
1+8\sigma^2 \text{tr}[\{(V+2\sigma \Sigma_0)^{-1} \Sigma_0\}^2]
+12\sigma^4 \{\text{tr}[\{(V+2\sigma \Sigma_0)^{-1} \Sigma_0\}^2]\}^2 \nonumber\\
&~~~~~~~~~~~~~~~~~~~~~~~~~~~~
+24\sigma^4 \text{tr}[\{(V+2\sigma \Sigma_0)^{-1} \Sigma_0\}^4]
\Big\}. \label{MATHCAL_F44}
\end{align}
(\ref{MATHCAL_F11}), (\ref{MATHCAL_F22}) and (\ref{MATHCAL_F44}) are combined into (\ref{VZ_SIMPLE}), which completes the proof of Proposition \ref{P2}. \qed

\subsection{Proof of Lemma \ref{L7} }

The essential point to obtain the expression of $V[Z]$ is that it finally consists of expectations of multiple for quadratic forms of Gaussian variable.
To see this, first we aim to find the expression of $\left<f(\underline{x})(\cdot),f(\underline{y})(\cdot)\right>_{H(k)}^2$ with terms of constant, linear form, bilinear form, quadratic form and multilple of quadratic forms.
It is easily confirmed from the definition of $B(\underline{x})$ in (\ref{DEF_Bx}) that
 	\begin{align*}
 	&\text{tr}[V^{-1} \{2\sigma (B(\underline{x})+\Sigma_0)V^{-1} -I_d\} B(\underline{y})]\\
 	&=2\sigma \Big((\underline{x}-\underline{m}_0)^T V^{-1} (\underline{y}-\underline{m}_0)\Big)^2
 	-2\sigma (\underline{x}-\underline{m}_0)V^{-1}\Sigma_0 V^{-1} (\underline{x}-\underline{m}_0) \\
 	&~~~~~
 	-(\underline{y}-\underline{m}_0)^T V^{-1} (\underline{y}-\underline{m}_0)
 	+\text{tr}[V^{-1} \Sigma_0],\\
 	&\text{tr}[(2(\underline{x}-\underline{m}_0)(\underline{y}-\underline{m}_0)^T-B(\underline{x}) -B(\underline{y}) )(2V-I_d)^{-1}]\\
 	&=2(\underline{x}-\underline{m}_0)^T (2V-I_d)^{-1} (\underline{y}-\underline{m}_0) 
 	-(\underline{x}-\underline{m}_0)^T (2V-I_d)^{-1} (\underline{x}-\underline{m}_0) \\
 	&~~~~~
 	-(\underline{y}-\underline{m}_0)^T (2V-I_d)^{-1} (\underline{y}-\underline{m}_0)
 	+2\text{tr}[\Sigma_0(2V-I_d)^{-1}],\\
 	&\text{tr}[B(\underline{x}) (2V-I_d)^{-1}]
 	\text{tr}[B(\underline{y}) (2V-I_d)^{-1}]\\
 	&=(\underline{x}-\underline{m}_0)^T (2V-I_d)^{-1} (\underline{x}-\underline{m}_0) 
 	(\underline{y}-\underline{m}_0)^T (2V-I_d)^{-1} (\underline{y}-\underline{m}_0) \\
 	&~~~~~-\text{tr}[\Sigma_0(2V-I_d)^{-1}] (\underline{x}-\underline{m}_0)^T (2V-I_d)^{-1} (\underline{x}-\underline{m}_0)\\
 	&~~~~~-\text{tr}[\Sigma_0(2V-I_d)^{-1}] (\underline{y}-\underline{m}_0)^T (2V-I_d)^{-1} (\underline{y}-\underline{m}_0)
 	+\{\text{tr}[\Sigma_0(2V-I_d)^{-1}]\}^2,
 	\end{align*}
 	and
 	\begin{align*}
 	&\text{tr}[B(\underline{x}) (2V-I_d)^{-1} B(\underline{y}) (2V-I_d)^{-1}]\\
 	&=\Big(
 	(\underline{x}-\underline{m}_0)^T (2V-I_d)^{-1} (\underline{y}-\underline{m}_0)
 	\Big)^2
 	-(\underline{x}-\underline{m}_0)^T (2V-I_d)^{-1} \Sigma_0(2V-I_d)^{-1} (\underline{x}-\underline{m}_0)\\
 	&~~~~~
 	-(\underline{y}-\underline{m}_0)^T (2V-I_d)^{-1} \Sigma_0(2V-I_d)^{-1} (\underline{y}-\underline{m}_0)
 	+\text{tr}[\{\Sigma_0(2V-I_d)^{-1}\}^2].
 	\end{align*}
These equalities give the another expression of (\ref{E47}) as
\begin{align}
&\left<f(\underline{x})(\cdot),f(\underline{y})(\cdot)\right>_{H(k)} \nonumber \\
&=\exp\Big(-\sigma(\underline{x}-\underline{m}_0)^T(\underline{x}-\underline{m}_0)-\sigma(\underline{y}-\underline{m}_0)^T(\underline{y}-\underline{m}_0)
+2\sigma (\underline{x}-\underline{m}_0)^T (\underline{y}-\underline{m}_0) \Big) \nonumber \\
&~~~~~-|V|^{-1/2} \exp\Big(
-\sigma(\underline{x}-\underline{m}_0)^T V^{-1} (\underline{x}-\underline{m}_0)
\Big)
\Big\{
1+\sigma \text{tr}[V^{-1} \Sigma_0]
+2\sigma (\underline{x}-\underline{m}_0)^T V^{-1} (\underline{y}-\underline{m}_0) \nonumber \\
&~~~~~~~~~~
+2\sigma^2 \Big(
(\underline{x}-\underline{m}_0)^T V^{-1} (\underline{y}-\underline{m}_0)
\Big)^2
-2\sigma^2 (\underline{x}-\underline{m}_0)^T V^{-1} \Sigma_0 V^{-1} (\underline{x}-\underline{m}_0)
-\sigma (\underline{y}-\underline{m}_0)^T V^{-1} (\underline{y}-\underline{m}_0)
\Big\} \nonumber \\
&~~~~~-|V|^{-1/2} \exp\Big(
-\sigma(\underline{y}-\underline{m}_0)^T V^{-1} (\underline{y}-\underline{m}_0)
\Big)
\Big\{
1+\sigma \text{tr}[V^{-1} \Sigma_0]
+2\sigma (\underline{x}-\underline{m}_0)^T V^{-1} (\underline{y}-\underline{m}_0) \nonumber \\
&~~~~~~~~~~
+2\sigma^2 \Big(
(\underline{x}-\underline{m}_0)^T V^{-1} (\underline{y}-\underline{m}_0)
\Big)^2
-\sigma (\underline{x}-\underline{m}_0)^T V^{-1} (\underline{x}-\underline{m}_0)
-2\sigma^2 (\underline{y}-\underline{m}_0)^T V^{-1} \Sigma_0 V^{-1} (\underline{y}-\underline{m}_0)  	
\Big\} \nonumber \\
&~~~~~+|2V-I_d|^{-1/2} \Big\{
1+2\sigma \text{tr}[\Sigma_0(2V-I_d)^{-1}] 
+\sigma^2 \{\text{tr}[\Sigma_0(2V-I_d)^{-1}]\}^2
+2\sigma^2 \text{tr}[\{\Sigma_0(2V-I_d)^{-1}\}^2] \nonumber \\
&~~~~~~~~~~+2\sigma (\underline{x}-\underline{m}_0)^T(2V-I_d)^{-1} (\underline{y}-\underline{m}_0)
+2\sigma^2 \Big(
(\underline{x}-\underline{m}_0)^T (2V-I_d)^{-1} (\underline{y}-\underline{m}_0)\Big)^2 \nonumber \\
&~~~~~~~~~~-\sigma(\underline{x}-\underline{m}_0)^T \Big(
(1+\sigma\text{tr}[\Sigma_0(2V-I_d)^{-1}] )I_d+2\sigma(2V-I_d)^{-1} \Sigma_0
\Big)(2V-I_d)^{-1}(\underline{x}-\underline{m}_0) \nonumber \\
&~~~~~~~~~~-\sigma(\underline{y}-\underline{m}_0)^T \Big(
(1+\sigma\text{tr}[\Sigma_0(2V-I_d)^{-1}] )I_d+2\sigma(2V-I_d)^{-1} \Sigma_0
\Big)(2V-I_d)^{-1}(\underline{y}-\underline{m}_0) \nonumber \\
&~~~~~~~~~~+\sigma^2 (\underline{x}-\underline{m}_0)^T(2V-I_d)^{-1}(\underline{x}-\underline{m}_0)
(\underline{y}-\underline{m}_0)^T(2V-I_d)^{-1}(\underline{y}-\underline{m}_0)
\Big\}. \label{ALT_EXP_1}
\end{align}
Now we introduce
  	\begin{align}
  	&\mathcal{T}[k_1,k_2,k_3|C_0,\dots,C_6|T_1,\dots,T_6](\underline{x},\underline{y}) \nonumber \\
  	&=\exp\Big(
  	-k_1\{\sigma \norm{\underline{x}-\underline{m}_0}^2_{\mathbb{R}^d} +\sigma \norm{\underline{y}-\underline{m}_0}^2_{\mathbb{R}^d} -2\sigma (\underline{x}-\underline{m}_0)^T (\underline{y}-\underline{m}_0) \} \nonumber \\
  	&~~~~~~~~~~
  	-k_2\sigma (\underline{x}-\underline{m}_0)^T V^{-1} (\underline{x}-\underline{m}_0)
  	-k_3\sigma (\underline{y}-\underline{m}_0)^T V^{-1} (\underline{y}-\underline{m}_0)
  	\Big) \nonumber \\
  	&~~~~~\times
  	C_{0}\left\{
  	C_1+C_2 (\underline{x}-\underline{m}_0)^T T_1 (\underline{y}-\underline{m}_0)
  	+C_3\{(\underline{x}-\underline{m}_0)^T T_2 (\underline{y}-\underline{m}_0)
  	\}^2
  	+C_4 (\underline{x}-\underline{m}_0)^T T_3 (\underline{x}-\underline{m}_0) \right. \nonumber \\
  	&~~~~~~~~~~
  	\left.+C_5 (\underline{y}-\underline{m}_0)^T T_4 (\underline{y}-\underline{m}_0)
  	+C_6 (\underline{x}-\underline{m}_0)^T T_5 (\underline{x}-\underline{m}_0) (\underline{y}-\underline{m}_0)^T T_6 (\underline{y}-\underline{m}_0)\right\} \label{MATHCAL_T_EXP}
  	\end{align}
	for binary variables  $k_{1},k_{2},k_{3}$, real values $C_{k}(k=0,1,...,6)$ and $d \times d$ matrices $T_{k}(k=1,...,6)$.
Then, by a careful check of the structure in (\ref{ALT_EXP_1}), we see that
\begin{align*}
\left<f(\underline{x})(\cdot),f(\underline{y})(\cdot)\right>_{H(k)}&
=\mathcal{T}[k^{(1)}_1,k^{(1)}_2,k^{(1)}_3|C^{(1)}_0,\dots,C^{(1)}_6|T^{(1)}_1,\dots,T^{(1)}_6](\underline{x},\underline{y}) \\
&\ \  -\mathcal{T}[k^{(2)}_1,k^{(2)}_2,k^{(2)}_3|C^{(2)}_0,\dots,C^{(2)}_6|T^{(2)}_1,\dots,T^{(2)}_6](\underline{x},\underline{y}) \\
&\ \  -\mathcal{T}[k^{(3)}_1,k^{(3)}_2,k^{(3)}_3|C^{(3)}_0,\dots,C^{(3)}_6|T^{(3)}_1,\dots,T^{(3)}_6)](\underline{x},\underline{y})\\
&\ \  +\mathcal{T}[k^{(4)}_1,k^{(4)}_2,k^{(4)}_3|C^{(4)}_0,\dots,C^{(4)}_6|T^{(4)}_1,\dots,T^{(4)}_6](\underline{x},\underline{y}) ,
\end{align*}
where $k^{(t)}_{i}(i=1,2,3)$, $C^{(t)}_{i}(i=0,1,...,6)$ and $T^{(t)}_{i}(i=1,...,6)$ for $t=1,2,3,4$ are those given in the proof of Proposition \ref{P2}.
Defining
\begin{align}
&\mathcal{F}(k_1,k_2,k_3|C_0,\dots,C_6|T_1,\dots,T_6\ ||\ k'_1,k'_2,k'_3|C'_0,\dots,C'_6|T'_1,\dots,T'_6) \nonumber\\
&=\int_{\mathbb{R}^d} \int_{\mathbb{R}^d} \mathcal{T}[k_1,k_2,k_3|C_0,\dots,C_6|T_1,\dots,T_6](\underline{x},\underline{y})  \nonumber \\
&~~~~~~~~~~~\times\mathcal{T}[k'_1,k'_2,k'_3|C'_0,\dots,C'_6|T'_1,\dots,T'_6](\underline{x},\underline{y}) dN(\underline{m}_0,\Sigma_0)(\underline{x}) dN(\underline{m}_0,\Sigma_0)(\underline{y}) \label{MATHCAL_F_EXP2}
\end{align}
derives the expression in Lemma \ref{L7}.

Finally we have to show that the integral (\ref{MATHCAL_F_EXP2}) certainly leads to (\ref{MATHCAL_F_EXP}).
Starting from change of variables $\underline{z}=\Sigma_0^{-1/2} (\underline{x}-\underline{m}_0),~~\underline{w}=\Sigma^{-1/2}_0(\underline{y}-\underline{m}_0)$, and we try to obtain the integral expression with the stacked variable $\underline{u}^{T}=[\underline{z}^{T}~~\underline{w}^{T}]$ of $2d$ dimension.
Now we see that change of variables above yields another expression of (\ref{MATHCAL_T_EXP}) as
\begin{align}
&\mathcal{T}[k_1,k_2,k_3|C_0,\dots,C_6|T_1,\dots,T_6](\underline{z},\underline{w})\nonumber \\
&=\exp\Big(
-k_1\{
\sigma \underline{z}^T \Sigma_0 \underline{z} +\sigma \underline{w}^T \Sigma_0 \underline{w}
-2\sigma \underline{z}^T \Sigma_0 \underline{w}
\}
-k_2 \sigma \underline{z}^T \Sigma_0^{1/2} V^{-1} \Sigma_0^{1/2} \underline{z} 
-k_3 \sigma \underline{w}^T \Sigma_0^{1/2} V^{-1} \Sigma_0^{1/2} \underline{w}
\Big) \nonumber \\
&~~~~~\times C_0\{
C_1+C_2 \underline{z}^T \Sigma_0^{1/2} T_1 \Sigma_0^{1/2} \underline{w} 
+C_3 \{\underline{z}^T \Sigma_0^{1/2} T_2 \Sigma_0^{1/2} \underline{w}\}^2
+C_4 \underline{z}^T \Sigma^{1/2}_0 T_3 \Sigma^{1/2}_0 \underline{z} \nonumber \\
&~~~~~~~~~~+C_5 \underline{w}^T \Sigma_0^{1/2} T_4 \Sigma_0^{1/2} \underline{w}
+C_6 \underline{z}^T \Sigma_0^{1/2} T_5 \Sigma_0^{1/2} \underline{z} \underline{w}^T \Sigma_0^{1/2} T_6 \Sigma_0^{1/2}\underline{w}
\}. \label{MATHCAL_T_EXP2}
\end{align}
Consider the integral of the product
$$
\mathcal{T}[k_1,k_2,k_3|C_0,\dots,C_6|T_1,\dots,T_6](\underline{z},\underline{w}) \times\mathcal{T}[k'_1,k'_2,k'_3|C'_0,\dots,C'_6|T'_1,\dots,T'_6](\underline{z},\underline{w}).
$$
The product of exponential parts including the Gaussian densities appeared from $\underline{Z}\sim N_{d}(\underline{0},I_{d})$ and $\underline{W}\sim N_{d}(\underline{0},I_{d})$ as seen in (\ref{MATHCAL_F_EXP2}) can be combined into
 \begin{align}
 &\exp\Big(
 -k_1\{
 \sigma \underline{z}^T \Sigma_0 \underline{z} +\sigma \underline{w}^T \Sigma_0 \underline{w}
 -2\sigma \underline{z}^T \Sigma_0 \underline{w}
 \}
 -k_2 \sigma \underline{z}^T \Sigma_0^{1/2} V^{-1} \Sigma_0^{1/2} \underline{z} 
 -k_3 \sigma \underline{w}^T \Sigma_0^{1/2} V^{-1} \Sigma_0^{1/2} \underline{w}
 \Big)
 \nonumber \\
 &~~~~~\times 
 \exp\Big(
 -k'_1\{
 \sigma \underline{z}^T \Sigma_0 \underline{z} +\sigma \underline{w}^T \Sigma_0 \underline{w}
 -2\sigma \underline{z}^T \Sigma_0 \underline{w}
 \}
 -k'_2 \sigma \underline{z}^T \Sigma_0^{1/2} V^{-1} \Sigma_0^{1/2} \underline{z} 
 -k'_3 \sigma \underline{w}^T \Sigma_0^{1/2} V^{-1} \Sigma_0^{1/2} \underline{w}
 \Big) \nonumber \\
 &~~~~~\times |2\pi I_d|^{-1} \exp \left(
 -\frac{1}{2}\underline{z}^T \underline{z}
 \right)
 \exp \left(
 -\frac{1}{2}\underline{w}^T \underline{w}
 \right) \nonumber \\
 &=(2\pi)^{-d} \exp\Bigg(
 -\frac{1}{2}\underline{z}^T \left\{I_d+2(k_1+k'_1) \sigma \Sigma_0+2(k_2+k'_2) \sigma \Sigma_0^{1/2} V^{-1} \Sigma^{1/2}_0 \right\} \underline{z} \nonumber \\
 &~~~~~~~~~~~~~~~~~~-\frac{1}{2}\underline{w}^T \left\{
 I_d+2(k_1+k'_1) \sigma \Sigma_0 +2(k_3+k'_3) \sigma \Sigma_0^{1/2} V^{-1} \Sigma_0^{1/2} 
 \right\}\underline{w} \nonumber \\
 &~~~~~~~~~~~~~~~~~~+\frac{1}{2} \underline{z}^T \left\{
 4(k_1+k'_1)\sigma \Sigma_0
 \right\} \underline{w}
 \Bigg) \nonumber \\
 &=(2\pi)^{-d} \exp\Bigg( -\frac{1}{2}\underline{u}^T R \underline{u}\Bigg), \label{EXP_PART}
 \end{align}
where $R$ is in (\ref{DEF_MAT_R}).
Further the product of non-exponential parts become to be
\begin{align}
&C_0\{
C_1+C_2 \underline{z}^T \Sigma_0^{1/2} T_1 \Sigma_0^{1/2} \underline{w} 
+C_3 \{\underline{z}^T \Sigma_0^{1/2} T_2 \Sigma_0^{1/2} \underline{w}\}^2
+C_4 \underline{z}^T \Sigma^{1/2}_0 T_3 \Sigma^{1/2}_0 \underline{z} \nonumber \\
&~~~~~+C_5 \underline{w}^T \Sigma_0^{1/2} T_4 \Sigma_0^{1/2} \underline{w}
+C_6 \underline{z}^T \Sigma_0^{1/2} T_5 \Sigma_0^{1/2} \underline{z} \underline{w}^T \Sigma_0^{1/2} T_6 \Sigma_0^{1/2}\underline{w}
\} \nonumber \\
&~~~~~\times C'_0\{
C'_1+C'_2 \underline{z}^T \Sigma_0^{1/2} T'_1 \Sigma_0^{1/2} \underline{w} 
+C'_3 \{\underline{z}^T \Sigma_0^{1/2} T'_2 \Sigma_0^{1/2} \underline{w}\}^2
+C'_4 \underline{z}^T \Sigma^{1/2}_0 T'_3 \Sigma^{1/2}_0 \underline{z} \nonumber \\
&~~~~~~~~~~+C'_5 \underline{w}^T \Sigma_0^{1/2} T'_4 \Sigma_0^{1/2} \underline{w}
+C'_6 \underline{z}^T \Sigma_0^{1/2} T'_5 \Sigma_0^{1/2} \underline{z} \underline{w}^T \Sigma_0^{1/2} T'_6 \Sigma_0^{1/2}\underline{w}
\} \nonumber \\
&=C_0\left\{
C_1+
\begin{bmatrix}
\underline{z}\\ 
\underline{w}
\end{bmatrix}^T
\begin{bmatrix}
C_4 \Sigma_0^{1/2} T_3 \Sigma_0^{1/2}& \frac{1}{2} C_2 \Sigma_0^{1/2} T_1 \Sigma_0^{1/2}  \\ 
\frac{1}{2} C_2 \Sigma_0^{1/2} T_1 \Sigma_0^{1/2} & C_5 \Sigma_0^{1/2} T_4 \Sigma_0^{1/2}
\end{bmatrix} 
\begin{bmatrix}
\underline{z}\\ 
\underline{w}
\end{bmatrix} \right.
 \nonumber\\
&~~~~~~~~~~+C_3 \left\{
\begin{bmatrix}
\underline{z}\\ 
\underline{w}
\end{bmatrix}^T
\begin{bmatrix}
O& \frac{1}{2} \Sigma_0^{1/2} T_2 \Sigma_0^{1/2}  \\ 
\frac{1}{2}\Sigma_0^{1/2} T_2 \Sigma_0^{1/2} & O
\end{bmatrix} 
\begin{bmatrix}
\underline{z}\\ 
\underline{w}
\end{bmatrix}
\right\}^2 \nonumber \\
&~~~~~~~~~~\left. +C_6 
\begin{bmatrix}
\underline{z}\\ 
\underline{w}
\end{bmatrix}^T
\begin{bmatrix}
\Sigma_0^{1/2} T_5 \Sigma_0^{1/2}& O  \\ 
O & O
\end{bmatrix} 
\begin{bmatrix}
\underline{z}\\ 
\underline{w}
\end{bmatrix}
\begin{bmatrix}
\underline{z}\\ 
\underline{w}
\end{bmatrix}^T
\begin{bmatrix}
O& O  \\ 
O & \Sigma_0^{1/2} T_6 \Sigma_0^{1/2}
\end{bmatrix} 
\begin{bmatrix}
\underline{z}\\ 
\underline{w}
\end{bmatrix}
\right\} \nonumber \\
&~~~~~\times C'_0\left\{
C'_1+
\begin{bmatrix}
\underline{z}\\ 
\underline{w}
\end{bmatrix}^T
\begin{bmatrix}
C'_4 \Sigma_0^{1/2} T'_3 \Sigma_0^{1/2}& \frac{1}{2} C'_2 \Sigma_0^{1/2} T'_1 \Sigma_0^{1/2}  \\ 
\frac{1}{2} C'_2 \Sigma_0^{1/2} T'_1 \Sigma_0^{1/2} & C'_5 \Sigma_0^{1/2} T'_4 \Sigma_0^{1/2}
\end{bmatrix} 
\begin{bmatrix}
\underline{z}\\ 
\underline{w}
\end{bmatrix} \right.
\nonumber \\
&~~~~~~~~~~+C'_3 \left\{
\begin{bmatrix}
\underline{z}\\ 
\underline{w}
\end{bmatrix}^T
\begin{bmatrix}
O& \frac{1}{2} \Sigma_0^{1/2} T'_2 \Sigma_0^{1/2}  \\ 
\frac{1}{2}\Sigma_0^{1/2} T'_2 \Sigma_0^{1/2} & O
\end{bmatrix} 
\begin{bmatrix}
\underline{z}\\ 
\underline{w}
\end{bmatrix}
\right\}^2 \nonumber \\
&~~~~~~~~~~\left. +C'_6 
\begin{bmatrix}
\underline{z}\\ 
\underline{w}
\end{bmatrix}^T
\begin{bmatrix}
\Sigma_0^{1/2} T'_5 \Sigma_0^{1/2}& O  \\ 
O & O
\end{bmatrix} 
\begin{bmatrix}
\underline{z}\\ 
\underline{w}
\end{bmatrix}
\begin{bmatrix}
\underline{z}\\ 
\underline{w}
\end{bmatrix}^T
\begin{bmatrix}
O& O  \\ 
O & \Sigma_0^{1/2} T'_6 \Sigma_0^{1/2}
\end{bmatrix} 
\begin{bmatrix}
\underline{z}\\ 
\underline{w}
\end{bmatrix}
\right\} \nonumber \\
&=C_0 \left\{
C_1+\underline{u}^T S_1\underline{u} +C_3 \{\underline{u}^T S_2 \underline{u}\}^2
+C_6 \underline{u}^T S_3 \underline{u} \underline{u}^T S_4 \underline{u}
\right\} \nonumber\\
&~~~~~\times C'_0 \left\{
C'_1+\underline{u}^T S'_1\underline{u} +C'_3 \{\underline{u}^T S'_2 \underline{u}\}^2
+C'_6 \underline{u}^T S'_3 \underline{u} \underline{u}^T S'_4 \underline{u}
\right\} \nonumber \\
&=C_0 C'_0 \Big\{
C_1 C'_1
+C_1 \underline{u}^T S'_1 \underline{u} +C'_1 \underline{u}^T S_1 \underline{u}
+C_1 C'_3 \{\underline{u}^T S_2' \underline{u}\}^2
+C'_1 C_3 \{\underline{u}^T S_2 \underline{u}\}^2 \nonumber\\
&~~~~~
+C_1C'_6 \underline{u}^T S'_3 \underline{u} \underline{u}^T S'_4 \underline{u}
+C'_1C_6 \underline{u}^T S_3 \underline{u} \underline{u}^T S_4 \underline{u}
+\underline{u}^T S_1 \underline{u} \underline{u}^T S'_1 \underline{u} \nonumber \\
&~~~~~
+C'_3 \underline{u}^T S_1 \underline{u} \{\underline{u}^T S'_2 \underline{u}\}^2
+C_3 \underline{u}^T S'_1 \underline{u} \{\underline{u}^T S_2 \underline{u}\}^2
+C'_6\underline{u}^T S_1 \underline{u} \underline{u}^T S'_3 \underline{u} \underline{u}^T S'_4 \underline{u}
+C_6\underline{u}^T S'_1 \underline{u} \underline{u}^T S_3 \underline{u} \underline{u}^T S_4 \underline{u} \nonumber\\
&~~~~~
+C_3C'_3 \{\underline{u}^T S_2 \underline{u}\}^2 \{\underline{u}^T S'_2 \underline{u}\}^2
+C_3 C'_6 \{\underline{u}^T S_2 \underline{u}\}^2 \underline{u}^T S'_3 \underline{u} \underline{u}^T S'_4 \underline{u}
+C'_3 C_6 \{\underline{u}^T S'_2 \underline{u}\}^2 \underline{u}^T S_3 \underline{u} \underline{u}^T S_4 \underline{u} \nonumber \\
&~~~~~+C_6C'_6 \underline{u}^T S_3 \underline{u} \underline{u}^T S_4 \underline{u} \underline{u}^T S'_3 \underline{u} \underline{u}^T S'_4 \underline{u}
\Big\}, \label{NON_EXP_PART}
\end{align}
where $S_{1}$, $S_{2}$, $S_{3}$ and $S_{4}$ are those in (\ref{DEF_MAT_S1}), (\ref{DEF_MAT_S2}), (\ref{DEF_MAT_S3}) and (\ref{DEF_MAT_S4}), respectively, and the dashed versions are corresponding to including $T'_{k}$s and $C'_{k}$s instead of $T_{k}$s and $C_{k}$s.
We then obtain by combining (\ref{MATHCAL_F_EXP2}), (\ref{MATHCAL_T_EXP2}), (\ref{EXP_PART}) and (\ref{NON_EXP_PART}) that
\begin{align*}
&\mathcal{F}(k_1,k_2,k_3|C_0,\dots,C_6|T_1,\dots,T_6||k'_1,k'_2,k'_3|C'_0,\dots,C'_6|T'_1,\dots,T'_6)\\
&=(2\pi)^{-d} C_0 C'_0\int_{\mathbb{R}^{2d}} 
\exp\left(-\frac{1}{2} \underline{u}^T R \underline{u}\right) 
\Big\{
C_1 C'_1
+C_1 \underline{u}^T S'_1 \underline{u} +C'_1 \underline{u}^T S_1 \underline{u}
+C_1 C'_3 \{\underline{u}^T S_2' \underline{u}\}^2
+C'_1 C_3 \{\underline{u}^T S_2 \underline{u}\}^2\\
&~~~~~
+C_1C'_6 \underline{u}^T S'_3 \underline{u} \underline{u}^T S'_4 \underline{u}
+C'_1C_6 \underline{u}^T S_3 \underline{u} \underline{u}^T S_4 \underline{u}
+\underline{u}^T S_1 \underline{u} \underline{u}^T S'_1 \underline{u}\\
&~~~~~
+C'_3 \underline{u}^T S_1 \underline{u} \{\underline{u}^T S'_2 \underline{u}\}^2
+C_3 \underline{u}^T S'_1 \underline{u} \{\underline{u}^T S_2 \underline{u}\}^2
+C'_6\underline{u}^T S_1 \underline{u} \underline{u}^T S'_3 \underline{u} \underline{u}^T S'_4 \underline{u}
+C_6\underline{u}^T S'_1 \underline{u} \underline{u}^T S_3 \underline{u} \underline{u}^T S_4 \underline{u}\\
&~~~~~
+C_3C'_3 \{\underline{u}^T S_2 \underline{u}\}^2 \{\underline{u}^T S'_2 \underline{u}\}^2
+C_3 C'_6 \{\underline{u}^T S_2 \underline{u}\}^2 \underline{u}^T S'_3 \underline{u} \underline{u}^T S'_4 \underline{u}
+C'_3 C_6 \{\underline{u}^T S'_2 \underline{u}\}^2 \underline{u}^T S_3 \underline{u} \underline{u}^T S_4 \underline{u}\\
&~~~~~+C_6C'_6 \underline{u}^T S_3 \underline{u} \underline{u}^T S_4 \underline{u} \underline{u}^T S'_3 \underline{u} \underline{u}^T S'_4 \underline{u}
\Big\} d\underline{u}\\
&=|R|^{-1/2} C_0C'_0 \int_{\mathbb{R}^d} \Big\{
C_1C'_1
+C_1 \underline{v}^T R^{-1/2} S'_1 R^{-1/2} \underline{v}
+C'_1 \underline{v}^T R^{-1/2} S_1 R^{-1/2} \underline{v}\\
&~~~~~+C_1C'_3 \{\underline{v}^T R^{-1/2} S'_2 R^{-1/2} \underline{v}\}^2
+C'_1C_3 \{\underline{v}^T R^{-1/2} S_2 R^{-1/2} \underline{v}\}^2\\
&~~~~~+C_1C'_6 \underline{v}^T R^{-1/2} S'_3 R^{-1/2} \underline{v} \underline{v}^T R^{-1/2} S'_4 R^{-1/2} \underline{v}
+C'_1C_6 \underline{v}^T R^{-1/2} S_3 R^{-1/2} \underline{v} \underline{v}^T R^{-1/2} S_4 R^{-1/2} \underline{v}\\
&~~~~~+\underline{v}^T R^{-1/2} S_1 R^{-1/2} \underline{v} \underline{v}^T R^{-1/2} S'_1 R^{-1/2} \underline{v}\\
&~~~~~+C'_3 \underline{v}^T R^{-1/2} S_1 R^{-1/2} \underline{v} \{\underline{v}^T R^{-1/2} S'_2 R^{-1/2} \underline{v}\}^2
+C_3 \underline{v}^T R^{-1/2} S'_1 R^{-1/2} \underline{v} \{\underline{v}^T R^{-1/2} S_2 R^{-1/2} \underline{v}\}^2\\
&~~~~~+C'_6 \underline{v}^T R^{-1/2} S_1 R^{-1/2} \underline{v} \underline{v}^T R^{-1/2} S'_3 R^{-1/2} \underline{v}\underline{v}^T R^{-1/2} S'_4 R^{-1/2} \underline{v}\\
&~~~~~+C_6 \underline{v}^T R^{-1/2} S'_1 R^{-1/2} \underline{v} \underline{v}^T R^{-1/2} S_3 R^{-1/2} \underline{v}\underline{v}^T R^{-1/2} S_4 R^{-1/2} \underline{v}\\
&~~~~~+C_3C'_3 \{\underline{v}^T R^{-1/2} S_2 R^{-1/2} \underline{v}\}^2
\{\underline{v}^T R^{-1/2} S'_2 R^{-1/2} \underline{v}\}^2\\
&~~~~~+C_3C'_6 \{\underline{v}^T R^{-1/2} S_2 R^{-1/2} \underline{v}\}^2 \underline{v}^T R^{-1/2} S'_3 R^{-1/2} \underline{v} \underline{v}^T R^{-1/2} S'_4 R^{-1/2} \underline{v} \\
&~~~~~+C'_3C_6 \{\underline{v}^T R^{-1/2} S'_2 R^{-1/2} \underline{v}\}^2 \underline{v}^T R^{-1/2} S_3 R^{-1/2} \underline{v} \underline{v}^T R^{-1/2} S_4 R^{-1/2} \underline{v} \\
&~~~~~+C_6C'_6 \underline{v}^T R^{-1/2} S_3 R^{-1/2} \underline{v}\underline{v}^T R^{-1/2} S_4 R^{-1/2} \underline{v}\underline{v}^T R^{-1/2} S'_3 R^{-1/2} \underline{v}\underline{v}^T R^{-1/2} S'_4 R^{-1/2} \underline{v}
\Big\}dN(\underline{0},I_{2d})(\underline{v}),
\end{align*}
where a trivial change of variable $\underline{v}=R^{1/2}\underline{u}$ has been used, and this expression implies (\ref{MATHCAL_F_EXP}) using the definitions (\ref{Q_1}), (\ref{Q_2}), (\ref{Q_3}) and (\ref{Q_4}) in Section 4.3.1.\qed

\subsection{Proof of Lemma \ref{L8} }
Equalities (\ref{J22=J33}), (\ref{J12=J13}) and (\ref{J34=J24}) are easily confirmed from the structures of $(k^{(2)}_{i},k^{(3)}_{i}),i=1,2,3$, $(C^{(2)}_{i},C^{(3)}_{i}),i=0,1,...,6$ and $(T^{(2)}_{i},T^{(3)}_{i}),i=1,...,6$.
In fact the (1,1)- and (2,2)-blocks of matrices appeared in the calculations for the left hand side are just permutated into (2,2)- and (1,1)-blocks of matrices  in those of the right hand side, which gives equalities.

Among (\ref{J22=J12}), (\ref{J44=J14}), (\ref{J44=J23}) and (\ref{J14=J24}), we give the proof for (\ref{J14=J24}) since it is the most complicated.
(\ref{J22=J12}), (\ref{J44=J14}), (\ref{J44=J23}) can be obtained in the same manner.

So we shall start to address the calculation of $\mathcal{F}_{14}$.
Remember that, as listed in Lemma \ref{L7}, 
\begin{align*}
	(k_{1},k_{2},k_{3})&=(k_{1}^{(1)},k_{2}^{(1)},k_{3}^{(1)})=(1,0,0), \\
	(k'_{1},k'_{2},k'_{3})&=(k_{1}^{(4)},k_{2}^{(4)},k_{3}^{(4)})=(0,0,0), \\
	(C_{0},C_{1},C_{2},C_{3},C_{4},C_{5},C_{6})&=(C^{(1)}_{0},C^{(1)}_{1},C^{(1)}_{2},C^{(1)}_{3},C^{(1)}_{4},C^{(1)}_{5},C^{(1)}_{6})=(1,1,0,0,0,0,0), \\
	(C'_{0},C'_{1},C'_{2},C'_{3},C'_{4},C'_{5},C'_{6})&=(C^{(4)}_{0},C^{(4)}_{2},C^{(4)}_{3},C^{(4)}_{4},C^{(4)}_{5},C^{(4)}_{6})=(|2V-I_d|^{-1/2},2\sigma,2\sigma^2,-\sigma,-\sigma,\sigma^2),
\end{align*}
which implies that
\begin{align}
&\mathcal{F}_{14} \nonumber \\
&=\mathcal{F}(k^{(1)}_1,k^{(1)}_2,k^{(1)}_3|C^{(1)}_0,\dots,C^{(1)}_6|T^{(1)}_1,\dots,T^{(1)}_6\ ||\ k^{(4)}_1,k^{(4)}_2,k^{(4)}_3|C^{(4)}_0,\dots,C^{(4)}_6|T^{(4)}_1,\dots,T^{(4)}_6) \nonumber \\
&=|R|^{-1/2} |V+2\sigma \Sigma_0|^{-1/2} \Big\{
C_1'+Q_1(S_1'R^{-1}) +2\sigma^2 Q_2(S_2'R^{-1},S_2'R^{-1}) +\sigma^2 Q_2(S_3'R^{-1} ,S_4'R^{-1})
\Big\} \label{MATHCAL_F14}
\end{align}
with referring the definition of $\mathcal{F}$.
Furthermore, in this case, matrices $R$, $S_{1}$, $S_{2}$ $S_{3}$ and $S_{4}$ defined respectively in (\ref{DEF_MAT_R}), (\ref{DEF_MAT_S1}), (\ref{DEF_MAT_S2}), (\ref{DEF_MAT_S3}) and (\ref{DEF_MAT_S4}) and those dashed version $S_{i}$s as well as matrix $T_{3}'$ are obtained as
 	\begin{itemize}
	\item $R=
	\begin{bmatrix}
	I_d+2\sigma \Sigma_0&-2\sigma \Sigma_0  \\ 
	-2\sigma \Sigma_0& I_d+2\sigma \Sigma_0 
	\end{bmatrix}=\begin{bmatrix}
	V &-2\sigma \Sigma_0  \\ 
	-2\sigma \Sigma_0& V 
	\end{bmatrix}, 
	$
	\item  $S_1=S_2=S_3=S_4=O$,
	\item $S_1'=
	\begin{bmatrix}
	-\sigma \Sigma_0^{1/2} T_3' \Sigma_0^{1/2}& \sigma \Sigma_0^{1/2} (V+2\sigma \Sigma_0)^{-1} \Sigma_0^{1/2}  \\ 
	\sigma \Sigma_0^{1/2} (V+2\sigma \Sigma_0)^{-1} \Sigma_0^{1/2}& -\sigma \Sigma_0^{1/2} T_3' \Sigma_0^{1/2} 
	\end{bmatrix}, 
	$
	\item  $S_2'=
	\begin{bmatrix}
	O&\dfrac{1}{2}\Sigma_0^{1/2} (V+2\sigma \Sigma_0)^{-1} \Sigma_0^{1/2}  \\ 
	\dfrac{1}{2}\Sigma_0^{1/2} (V+2\sigma \Sigma_0)^{-1} \Sigma_0^{1/2}&O 
	\end{bmatrix} $
	\item  $S_3'=
	\begin{bmatrix}
	\Sigma_0^{1/2} (V+2\sigma \Sigma_0)^{-1} \Sigma_0^{1/2}& O \\ 
	O&O 
	\end{bmatrix}, 
	$
	\item $S_4'=
	\begin{bmatrix}
	O&O  \\ 
	O&\Sigma_0^{1/2} (V+2\sigma \Sigma_0)^{-1} \Sigma_0^{1/2} 
	\end{bmatrix}, 
	$
\end{itemize}
where we note again $C_{1}'=C_{1}^{(4)}$ and $T_{3}'=T_{3}^{(4)}$ given in Lemma \ref{L7}.
We need to clarify each term appeared in $\mathcal{F}_{14}$.
First we see that
\begin{align*}
\{T_3'\Sigma_0\}^2
&=\{
(V+2\sigma \Sigma_0)^{-1} \Sigma_0+\sigma \text{tr}[(V+2\sigma \Sigma_0)^{-1} \Sigma_0] (V+2\sigma \Sigma_0)^{-1} \Sigma_0+2\sigma\{(V+2\sigma \Sigma_0)^{-1} \Sigma_0\}^2
\}^2\\
&=\{(V+2\sigma \Sigma_0)^{-1} \Sigma_0\}^2
+2\sigma \text{tr}[(V+2\sigma \Sigma_0)^{-1} \Sigma_0] \{(V+2\sigma \Sigma_0)\}^2
+4\sigma \{(V+2\sigma \Sigma_0)^{-1} \Sigma_0\}^3\\
&~~~~~+\sigma^2 \{\text{tr}[(V+2\sigma \Sigma_0)^{-1} \Sigma_0]\}^2 \{(V+2\sigma \Sigma_0)^{-1} \Sigma_0\}^2
\\
&~~~~~
+4\sigma^2 \text{tr}[(V+2\sigma \Sigma_0)^{-1} \Sigma_0]\{(V+2\sigma \Sigma_0)^{-1} \Sigma_0\}^3
+4\sigma^2 \{(V+2\sigma \Sigma_0)^{-1} \Sigma_0\}^4.
\end{align*}
And we have the inverse
\begin{align}
R^{-1} 
&=
\begin{bmatrix}
V^{-1}+4\sigma^2V^{-1}\Sigma_0 S^{-1} \Sigma_0 V^{-1}& 2\sigma V^{-1} \Sigma_0 S^{-1}  \\ 
2\sigma S^{-1}\Sigma_0 V^{-1}&S^{-1} 
\end{bmatrix},\label{R_MAT_INV} 
\end{align}
where
\begin{align*}
S
&=V-4\sigma^2 \Sigma_0 V^{-1} \Sigma_0\\
&=(\Sigma_0^{-1} +2\sigma I_d-4\sigma^2 \Sigma_0 V^{-1}) \Sigma_0\\
&=(\Sigma_0^{-1}V+2\sigma V-4\sigma^2 \Sigma_0) V^{-1}\Sigma_0\\
&=(\Sigma_0^{-1}+2\sigma I_d +2\sigma I_d +4\sigma^2 \Sigma_0 -4\sigma^2 \Sigma_0) V^{-1}\Sigma_0\\
&=(\Sigma_0^{-1}+4\sigma I_d) V^{-1}\Sigma_0\\
&=\Sigma_0^{-1/2}(I_d+4\sigma \Sigma_0) \Sigma_0^{-1/2} V^{-1}\Sigma_0\\
&=\Sigma_0^{-1/2} (V+2\sigma \Sigma_0) \Sigma_0^{-1/2} V^{-1} \Sigma_0. 
\end{align*}
Hence the (1,1)-block of (\ref{R_MAT_INV}) can be calculated as
\begin{align*}
&V^{-1}+4\sigma^2 V^{-1} \Sigma_0 \Sigma_0^{-1} V \Sigma_0^{1/2} (V+2\sigma \Sigma_0)^{-1} \Sigma_0^{3/2} V^{-1}\\
&=V^{-1}+4\sigma^2\Sigma_0^{1/2} (V+2\sigma \Sigma_0)^{-1} \Sigma_0^{3/2} V^{-1}\\
&=(I_d+4\sigma^2 \Sigma_0^{1/2} (V+2\sigma \Sigma_0)^{-1} \Sigma_0^{3/2} )V^{-1}\\
&=\Sigma_0^{1/2} (\Sigma_0^{-2}+4\sigma^2 (V+2\sigma \Sigma_0)^{-1}) \Sigma_0^{3/2} V^{-1}\\
&=\Sigma_0^{1/2} (V+2\sigma \Sigma_0)^{-1} \left(
(V+2\sigma \Sigma_0) \Sigma_0^{-2} +4\sigma^2 I_d
\right) \Sigma_0^{3/2} V^{-1}\\
&=\Sigma_0^{1/2} (V+2\sigma \Sigma_0)^{-1} (V\Sigma_0^{-2} +2\sigma \Sigma_0^{-1} +4\sigma^2 I_d) \Sigma_0^{3/2} V^{-1}\\
&=\Sigma_0^{1/2} (V+2\sigma \Sigma_0)^{-1} \left(
(I_d+2\sigma \Sigma_0)\Sigma_0^{-2} +2\sigma \Sigma_0^{1/2} (I_d+2\sigma \Sigma_0) \Sigma_0^{-3/2}
\right) \Sigma_0^{3/2}  V^{-1}\\
&=\Sigma_0^{1/2} (V+2\sigma \Sigma_0)^{-1} \left(
\Sigma_0^{-1/2} (I_d+2\sigma \Sigma_0) \Sigma_0^{-3/2}+2\sigma \Sigma_0^{1/2} V \Sigma_0^{1/2} 
\right) \Sigma_0^{3/2} V^{-1}\\
&=\Sigma_0^{1/2} (V+2\sigma \Sigma_0)^{-1} (\Sigma_0^{-1/2} +2\sigma \Sigma_0^{1/2})\\
&=\Sigma_0^{1/2} (V+2\sigma \Sigma_0)^{-1} (I_d+2\sigma \Sigma_0) \Sigma_0^{-1/2} \\
&= \Sigma_0^{1/2} (V+2\sigma \Sigma_0)^{-1} V \Sigma_0^{-1/2} \\
&= \Sigma_0^{1/2} (V+2\sigma \Sigma_0)^{-1} (V+2\sigma \Sigma_0-2\sigma \Sigma_0) \Sigma_0^{-1/2}\\
&=I_d-2\sigma \Sigma_0^{1/2} (V+2\sigma \Sigma_0)^{-1} \Sigma_0^{1/2}.
\end{align*}
The (2,2)-block of (\ref{R_MAT_INV}) is
\begin{align*}
S^{-1}
&=\Sigma_0^{-1} V \Sigma_0^{1/2} (V+2\sigma \Sigma_0)^{-1} \Sigma_0\\
&=\Sigma_0^{-1} (I_d+2\sigma \Sigma_0) \Sigma_0^{1/2} (V+2\sigma \Sigma_0)^{-1} \Sigma_0^{1/2}\\
&=\Sigma_0^{-1/2} (I_d+2\sigma \Sigma_0) (V+2\sigma \Sigma_0)^{-1} \Sigma_0^{1/2}\\
&=\Sigma_0^{-1/2} (V+2\sigma \Sigma_0-2\sigma \Sigma_0) (V+2\sigma \Sigma_0)^{-1} \Sigma_0^{1/2}\\
&=I_d-2\sigma \Sigma_0^{1/2} (V+2\sigma \Sigma_0)^{-1} \Sigma_0^{1/2}.
\end{align*}
And off-diagonal blocks of (\ref{R_MAT_INV}) are both calculated as
\begin{align*}
2\sigma\Sigma_0^{1/2} (V+2\sigma \Sigma_0)^{-1} \Sigma_0^{1/2}. 
\end{align*}
Above expressions are combined into
\begin{align*}
R^{-1}
&=
\begin{bmatrix}
I_d-2\sigma \Sigma_0^{1/2} (V+2\sigma \Sigma_0)^{-1} \Sigma_0^{1/2}& 2\sigma \Sigma_0^{1/2} (V+2\sigma \Sigma_0)^{-1} \Sigma_0 \\ 
2\sigma \Sigma_0^{1/2} (V+2\sigma \Sigma_0)^{-1} \Sigma_0&I_d-2\sigma \Sigma_0^{1/2} (V+2\sigma \Sigma_0)^{-1} \Sigma_0^{1/2} 
\end{bmatrix}.
\end{align*}
The determinant of $R$ is easily evaluated as
\begin{align*}
|R|
&=|V||V-4\sigma^2 \Sigma_0 V^{-1} \Sigma_0|\\
&=|V||S|\\
&=|V||\Sigma_0^{-1/2} (V+2\sigma \Sigma_0) \Sigma_0^{-1/2} V^{-1} \Sigma_0|\\
&=|V+2\sigma \Sigma_0|.
\end{align*}
Hereafter we aim to obtain the expressions of the matrices appeared in $Q_{1}$ and $Q_{2}$ in $\mathcal{F}_{14}$.
Because we will take the trace in the final calculations, it suffices to obtain the diagonal blocks concretely.
Straightforward calculations show that those are in fact obtained as  
\begin{align*}
&S_1'R^{-1}
=
\begin{bmatrix}
-\sigma \Sigma_0^{1/2} T_3' \Sigma_0^{1/2} +2\sigma^2 \Sigma_0^{1/2} T_3' \Sigma_0 (V+2\sigma \Sigma_0)^{-1} \Sigma_0^{1/2} +2\sigma^2\{\Sigma_0^{1/2} (V+2\sigma \Sigma_0)^{-1} \Sigma_0^{1/2}\}^2~~~~~  *\\
*~~~~~ -\sigma \Sigma_0^{1/2} T_3' \Sigma_0^{1/2} +2\sigma^2 \Sigma_0^{1/2} T_3' \Sigma_0 (V+2\sigma \Sigma_0)^{-1} \Sigma_0^{1/2} +2\sigma^2\{\Sigma_0^{1/2} (V+2\sigma \Sigma_0)^{-1} \Sigma_0^{1/2}\}^2
\end{bmatrix}, \\
&S_2'R^{-1}
=
\begin{bmatrix}
\sigma \{\Sigma_0^{1/2}  (V+2\sigma \Sigma_0)^{-1} \Sigma_0^{1/2}\}^2~~~~~ \dfrac{1}{2}\Sigma_0^{1/2}  (V+2\sigma \Sigma_0)^{-1} \Sigma_0^{1/2}-\sigma \{\Sigma_0^{1/2}  (V+2\sigma \Sigma_0)^{-1} \Sigma_0^{1/2}\}^2 \\ 
\dfrac{1}{2}\Sigma_0^{1/2}  (V+2\sigma \Sigma_0)^{-1} \Sigma_0^{1/2}-\sigma \{\Sigma_0^{1/2}  (V+2\sigma \Sigma_0)^{-1} \Sigma_0^{1/2}\}^2~~~~~ \sigma \{\Sigma_0^{1/2}  (V+2\sigma \Sigma_0)^{-1} \Sigma_0^{1/2}\}^2
\end{bmatrix}, \\
&\{S_2'R^{-1}\}^2
=
\begin{bmatrix}
2\sigma^2\{\Sigma_0^{1/2}  (V+2\sigma \Sigma_0)^{-1} \Sigma_0^{1/2}\}^4-\sigma \{\Sigma_0^{1/2}  (V+2\sigma \Sigma_0)^{-1} \Sigma_0^{1/2}\}^3~~~~~*\\
+\dfrac{1}{4} \{\Sigma_0^{1/2}  (V+2\sigma \Sigma_0)^{-1} \Sigma_0^{1/2}\}^2~~~~~~~~~~~~~~~~~~~~~~~~~ \\ 
~~~~~~2\sigma^2\{\Sigma_0^{1/2}  (V+2\sigma \Sigma_0)^{-1} \Sigma_0^{1/2}\}^4-\sigma \{\Sigma_0^{1/2}  (V+2\sigma \Sigma_0)^{-1} \Sigma_0^{1/2}\}^3\\
*~~~~~~~~~~~~~~~~~~~~~~~~~~~~~~~~~~~~~~~~~~~ +\dfrac{1}{4} \{\Sigma_0^{1/2}  (V+2\sigma \Sigma_0)^{-1} \Sigma_0^{1/2}\}^2\\ 
\end{bmatrix}, \\
&S_3'R^{-1}
=
\begin{bmatrix}
\Sigma_0^{1/2}(V+2\sigma \Sigma_0)^{-1}\Sigma_0^{1/2}-2\sigma\{\Sigma_0^{1/2}(V+2\sigma \Sigma_0)^{-1}\Sigma_0^{1/2}\}^2& 2\sigma \{\Sigma_0^{1/2}(V+2\sigma \Sigma_0)^{-1}\Sigma_0^{1/2}\}^2  \\ 
O&O 
\end{bmatrix}, \\
&S_4'R^{-1}
=
\begin{bmatrix}
O&O \\
2\sigma \{\Sigma_0^{1/2}(V+2\sigma \Sigma_0)^{-1}\Sigma_0^{1/2}\}^2   	&
\Sigma_0^{1/2}(V+2\sigma \Sigma_0)^{-1}\Sigma_0^{1/2}-2\sigma\{\Sigma_0^{1/2}(V+2\sigma \Sigma_0)^{-1}\Sigma_0^{1/2}\}^2
\end{bmatrix} 
\end{align*}
and
\begin{align*}
S_3'R^{-1}S_4'R^{-1}
&=
\begin{bmatrix}
4\sigma^2\{\Sigma_0^{1/2}(V+2\sigma \Sigma_0)^{-1}\Sigma_0^{1/2}\}^4&*  \\ 
*&O 
\end{bmatrix}. 
\end{align*}
By the use of above expressions, we see that
\begin{align*}
&Q_1(S_1'R^{-1})\\
&=\text{tr}[S_1'R^{-1}]\\
&=-2\sigma \text{tr}[T_3' \Sigma_0]
+4\sigma^2 \text{tr}[T_3' \Sigma_0 (V+2\sigma \Sigma_0)^{-1} \Sigma_0] 
+4\sigma^2 \text{tr}[\{(V+2\sigma \Sigma_0)^{-1} \Sigma_0\}^2]\\
&=-2\sigma \text{tr}[(V+2\sigma \Sigma_0)^{-1} \Sigma_0+\sigma \text{tr}[(V+2\sigma \Sigma_0)^{-1} \Sigma_0] (V+2\sigma \Sigma_0)^{-1} \Sigma_0+2\sigma\{(V+2\sigma \Sigma_0)^{-1} \Sigma_0\}^2]\\
&~~~~~+4\sigma^2 \text{tr}[\{(V+2\sigma \Sigma_0)^{-1} \Sigma_0+\sigma \text{tr}[(V+2\sigma \Sigma_0)^{-1} \Sigma_0] (V+2\sigma \Sigma_0)^{-1} \Sigma_0+2\sigma\{(V+2\sigma \Sigma_0)^{-1} \Sigma_0\}^2\}\\
&~~~~~~~~~~\times (V+2\sigma \Sigma_0)^{-1} \Sigma_0]\\
&~~~~~+4\sigma^2 \text{tr}[\{(V+2\sigma \Sigma_0)^{-1}\Sigma_0\}^2]\\
&=-2\sigma \text{tr}[ (V+2\sigma \Sigma_0)^{-1} \Sigma_0]
-2\sigma^2 \{\text{tr}[ (V+2\sigma \Sigma_0)^{-1} \Sigma_0]\}^2
-4\sigma^2 \text{tr}[\{ (V+2\sigma \Sigma_0)^{-1} \Sigma_0\}^2]\\
&~~~~~+4\sigma^2 \text{tr}[\{ (V+2\sigma \Sigma_0)^{-1} \Sigma_0\}^2]
+4\sigma^3 \text{tr}[ (V+2\sigma \Sigma_0)^{-1} \Sigma_0] \text{tr}[\{ (V+2\sigma \Sigma_0)^{-1} \Sigma_0\}^2]\\
&~~~~~
+8\sigma^3 \text{tr}[\{ (V+2\sigma \Sigma_0)^{-1} \Sigma_0\}^3]
+4\sigma^2 \text{tr}[\{ (V+2\sigma \Sigma_0)^{-1} \Sigma_0\}^2]\\
&=-2\sigma \text{tr}[ (V+2\sigma \Sigma_0)^{-1} \Sigma_0]
-2\sigma^2 \{\text{tr}[ (V+2\sigma \Sigma_0)^{-1} \Sigma_0]\}^2
+4\sigma^2 \text{tr}[\{ (V+2\sigma \Sigma_0)^{-1} \Sigma_0\}^2]\\
&~~~~~+4\sigma^3 \text{tr}[ (V+2\sigma \Sigma_0)^{-1} \Sigma_0] \text{tr}[\{ (V+2\sigma \Sigma_0)^{-1} \Sigma_0\}^2]
+8\sigma^3\text{tr}[\{ (V+2\sigma \Sigma_0)^{-1} \Sigma_0\}^3],\\
&Q_2(S_2'R^{-1},S_2'R^{-1})\\
&=\{\text{tr}[S_2'R^{-1}]\}^2+2\text{tr}[\{S_2'R^{-1}\}^2]\\
&=\{2\sigma\text{tr}[\{(V+2\sigma \Sigma_0)^{-1} \Sigma_0\}^2]\}^2
+8\sigma^2\text{tr}[\{(V+2\sigma \Sigma_0)^{-1} \Sigma_0\}^4]
-4\sigma \text{tr}[\{(V+2\sigma \Sigma_0)^{-1} \Sigma_0\}^3]\\
&~~~~~
+\text{tr}[\{(V+2\sigma \Sigma_0)^{-1} \Sigma_0\}^2]\\
&=4\sigma^2 \{\text{tr}[\{(V+2\sigma \Sigma_0)^{-1} \Sigma_0\}^2]\}^2
+8\sigma^2\text{tr}[\{(V+2\sigma \Sigma_0)^{-1} \Sigma_0\}^4]
-4\sigma \text{tr}[\{(V+2\sigma \Sigma_0)^{-1} \Sigma_0\}^3]\\
&~~~~~
+\text{tr}[\{(V+2\sigma \Sigma_0)^{-1} \Sigma_0\}^2]
\end{align*}
and
\begin{align*}
&Q_2(S_3'R^{-1},S_4'R^{-1})\\
&=\text{tr}[S_3'R^{-1}] \text{tr}[S_4'R^{-1}] +2\text{tr}[S_3'R^{-1} S_4'R^{-1}]\\
&=\{\text{tr}[(V+2\sigma\Sigma_0)^{-1}\Sigma_0]-2\sigma \text{tr}[\{(V+2\sigma\Sigma_0)^{-1}\Sigma_0\}^2]\}^2
+8\sigma^2 \text{tr}[\{(V+2\sigma\Sigma_0)^{-1}\Sigma_0\}^4]\\
&=\{\text{tr}[(V+2\sigma\Sigma_0)^{-1}\Sigma_0]\}^2
-4\sigma \text{tr}[(V+2\sigma\Sigma_0)^{-1}\Sigma_0] \text{tr}[\{(V+2\sigma\Sigma_0)^{-1}\Sigma_0\}^2]\\
&~~~~~
+4\sigma^2 \{\text{tr}[\{(V+2\sigma\Sigma_0)^{-1}\Sigma_0\}^2]\}^2
+8\sigma^2\text{tr}[\{(V+2\sigma\Sigma_0)^{-1}\Sigma_0\}^4].
\end{align*}
Substituting all expressions above into (\ref{MATHCAL_F14}) yields that
\begin{align}
&\mathcal{F}_{14}=|V+2\sigma \Sigma_0|^{-1} \Big\{
1+8\sigma^2 \text{tr}[\{(V+2\sigma \Sigma_0)^{-1} \Sigma_0\}^2]
+12\sigma^4 \{\text{tr}[\{(V+2\sigma \Sigma_0)^{-1} \Sigma_0\}^2]\}^2 \nonumber\\
&~~~~~~~~~~~~~~~~~~~~~~~~~~~~~~~~~~+24\sigma^4 \text{tr}[\{(V+2\sigma \Sigma_0)^{-1} \Sigma_0\}^4]
\Big\}. \label{equation_F14}
\end{align}
Next our focus goes to $\mathcal{F}_{24}$.
By noting the $k_{i}^{(2)}(i=1,2,3)$, $C_{i}^{(2)}(i=0,1,...,6)$ and $T_{i}^{(2)}(i=1,...,6)$ listed in Lemma \ref{L7}, we have 
\begin{align}
&\mathcal{F}_{24} \nonumber\\
&=\mathcal{F}(k^{(2)}_1,k^{(2)}_2,k^{(2)}_3|C^{(2)}_0,\dots,C^{(2)}_6|T^{(2)}_1,\dots,T^{(2)}_6\ ||\ k^{(4)}_1,k^{(4)}_2,k^{(4)}_3|C^{(4)}_0,\dots,C^{(4)}_6|T^{(4)}_1,\dots,T^{(4)}_6)\nonumber \\
&=|R|^{-1/2} |V|^{-1/2} |V+2\sigma \Sigma_0|^{-1/2} \Big\{
(1+\sigma \text{tr}[V^{-1} \Sigma_0]) C_1'
+(1+\sigma\text{tr}[V^{-1} \Sigma_0])Q_1(S_1'R^{-1})
+C_1' Q_1(S_1R^{-1}) \nonumber\\
&~~~~~
+2\sigma^2(1+\sigma \text{tr}[V^{-1} \Sigma_0])Q_2(S_2'R^{-1},S_2'R^{-1})
+2\sigma^2 C_1' Q_2(S_2R^{-1},S_2R^{-1})\nonumber \\
&~~~~~
+\sigma^2 (1+\sigma \text{tr}[V^{-1} \Sigma_0]) Q_2(S_3'R^{-1} ,S_4'R^{-1})
+Q_2(S_1R^{-1},S_1'R^{-1})
+2\sigma^2Q_3(S_1R^{-1},S_2'R^{-1},S_2'R^{-1})\nonumber \\
&~~~~~
+2\sigma^2 Q_3(S_1'R^{-1},S_2R^{-1},S_2R^{-1})
+\sigma^2 Q_3(S_1R^{-1},S_3'R^{-1},S_4'R^{-1})\nonumber \\
&~~~~~
+4\sigma^4 Q_4(S_2R^{-1},S_2R^{-1},S_2'R^{-1},S_2'R^{-1})
+2\sigma^4 Q_4(S_2R^{-1},S_2R^{-1},S_3'R^{-1},S_4'R^{-1})
\Big\} \label{MATHCAL_F24}
\end{align}
For the sequent calculations, we shall figure the necessary quantities as follows.
	\begin{itemize}
	\item $R=
	\begin{bmatrix}
	I_d+2\sigma \Sigma_0^{1/2} V^{-1} \Sigma_0^{1/2}& O \\ 
	O&I_d 
	\end{bmatrix} 
	$		
	\item $S_1=
	\begin{bmatrix}
	-2\sigma^2 \Sigma_0^{1/2} V^{-1} \Sigma_0 V^{-1} \Sigma_0^{1/2}&  \sigma \Sigma_0^{1/2} V^{-1} \Sigma_0^{1/2} \\ 
	\sigma \Sigma_0^{1/2} V^{-1} \Sigma_0^{1/2}&-\sigma \Sigma_0^{1/2} V^{-1} \Sigma_0^{1/2} 
	\end{bmatrix} 
	$
	\item $S_2=
	\begin{bmatrix}
	O&\dfrac{1}{2}\Sigma_0^{1/2} V^{-1} \Sigma_0^{1/2}  \\ 
	\dfrac{1}{2}\Sigma_0^{1/2} V^{-1} \Sigma_0^{1/2}&O 
	\end{bmatrix} 
	$
	\item  $S_3=S_4=O$
	\item $S_1'=
	\begin{bmatrix}
	-\sigma \Sigma_0^{1/2} T_3' \Sigma_0^{1/2}& \sigma \Sigma_0^{1/2} (V+2\sigma \Sigma_0)^{-1} \Sigma_0^{1/2}  \\ 
	\sigma \Sigma_0^{1/2} (V+2\sigma \Sigma_0)^{-1} \Sigma_0^{1/2}& -\sigma \Sigma_0^{1/2} T_3' \Sigma_0^{1/2} 
	\end{bmatrix} 
	$
	\item  $S_2'=
	\begin{bmatrix}
	O&\dfrac{1}{2}\Sigma_0^{1/2} (V+2\sigma \Sigma_0)^{-1} \Sigma_0^{1/2}  \\ 
	\dfrac{1}{2}\Sigma_0^{1/2} (V+2\sigma \Sigma_0)^{-1} \Sigma_0^{1/2}&O 
	\end{bmatrix} $
	\item  $S_3'=
	\begin{bmatrix}
	\Sigma_0^{1/2} (V+2\sigma \Sigma_0)^{-1} \Sigma_0^{1/2}& O \\ 
	O&O 
	\end{bmatrix} 
	$
	\item $S_4'=
	\begin{bmatrix}
	O&O  \\ 
	O&\Sigma_0^{1/2} (V+2\sigma \Sigma_0)^{-1} \Sigma_0^{1/2} 
	\end{bmatrix} 
	$
\end{itemize}
And we further notice that $C'_{1}=C_{1}^{(4)}$ and $T'_{3}=T_{3}^{(4)}$, both of which can be found in Lemma \ref{L7} as well.
As preliminary calculation, we have
\begin{align*}
\{T_3'\Sigma_0\}^2
&=\{
(V+2\sigma \Sigma_0)^{-1} \Sigma_0+\sigma \text{tr}[(V+2\sigma \Sigma_0)^{-1} \Sigma_0] (V+2\sigma \Sigma_0)^{-1} \Sigma_0+2\sigma\{(V+2\sigma \Sigma_0)^{-1} \Sigma_0\}^2
\}^2\\
&=\{(V+2\sigma \Sigma_0)^{-1} \Sigma_0\}^2
+2\sigma \text{tr}[(V+2\sigma \Sigma_0)^{-1} \Sigma_0] \{(V+2\sigma \Sigma_0)\}^2
+4\sigma \{(V+2\sigma \Sigma_0)^{-1} \Sigma_0\}^3\\
&~~~~~+\sigma^2 \{\text{tr}[(V+2\sigma \Sigma_0)^{-1} \Sigma_0]\}^2 \{(V+2\sigma \Sigma_0)^{-1} \Sigma_0\}^2\\
&~~~~~
+4\sigma^2 \text{tr}[(V+2\sigma \Sigma_0)^{-1} \Sigma_0]\{(V+2\sigma \Sigma_0)^{-1} \Sigma_0\}^3
+4\sigma^2 \{(V+2\sigma \Sigma_0)^{-1} \Sigma_0\}^4,\\
|R|
&=|I_d+2\sigma \Sigma_0^{1/2} V^{-1} \Sigma_0^{1/2}|\\
&=|\Sigma_0||\Sigma_0^{-1}+2\sigma V^{-1}|\\
&=|V|^{-1}|V+2\sigma\Sigma_0|
\end{align*}
and
\begin{align*}
R^{-1}
&=
\begin{bmatrix}
(I_d+2\sigma \Sigma_0^{1/2}V^{-1} \Sigma_0^{1/2})^{-1}&O  \\ 
O&I_d 
\end{bmatrix} \\
&=\begin{bmatrix}
\Sigma_0^{-1/2} V(V+2\sigma \Sigma_0)^{-1} \Sigma_0^{1/2}&O  \\ 
O&I_d 
\end{bmatrix} .
\end{align*}
Using these, it follows that
\begin{align*}
&S_1' R^{-1}
=
\begin{bmatrix}
-\sigma \Sigma_0^{1/2} T_3' V(V+2\sigma \Sigma_0)^{-1} \Sigma_0^{1/2}&\sigma \Sigma_0^{1/2} (V+2\sigma \Sigma_0)^{-1} \Sigma_0^{1/2}  \\ 
\sigma \Sigma_0^{1/2} (V+2\sigma \Sigma_0)^{-1} V(V+2\sigma \Sigma_0)^{-1} \Sigma_0^{1/2}&-\sigma \Sigma_0^{1/2} T_3' \Sigma_0^{1/2} 
\end{bmatrix} \\
&\hspace{1.1cm}=
\begin{bmatrix}
-\sigma \Sigma_0^{1/2} T_3' \Sigma_0^{1/2}+2\sigma^2 \Sigma_0^{1/2} T_3' \Sigma_0 (V+2\sigma \Sigma_0)^{-1} \Sigma_0^{1/2} &\sigma \Sigma_0^{1/2} (V+2\sigma \Sigma_0)^{-1} \Sigma_0^{1/2}  \\ 
\sigma \Sigma_0^{1/2} (V+2\sigma \Sigma_0)^{-1}\Sigma_0^{1/2}-2\sigma^2 \{\Sigma_0^{1/2} (V+2\sigma \Sigma_0)^{-1} \Sigma_0^{1/2} \}^2&-\sigma \Sigma_0^{1/2} T_3' \Sigma_0^{1/2} 
\end{bmatrix},\\
&S_1R^{-1}
=
\begin{bmatrix}
-2\sigma^2 \Sigma_0^{1/2} V^{-1} \Sigma_0 (V+2\sigma \Sigma_0)^{-1} \Sigma_0^{1/2}&\sigma \Sigma_0^{1/2} V^{-1} \Sigma_0^{1/2}  \\ 
\sigma \Sigma_0^{1/2} (V+2\sigma \Sigma_0)^{-1} \Sigma_0^{1/2}& -\sigma \Sigma_0^{1/2} V^{-1} \Sigma_0^{1/2}
\end{bmatrix}, \\
&S_2'R^{-1}
=
\begin{bmatrix}
O& \dfrac{1}{2}\Sigma_0^{1/2} (V+2\sigma \Sigma_0)^{-1} \Sigma_0^{1/2} \\ 
\dfrac{1}{2}\Sigma_0^{1/2} (V+2\sigma \Sigma_0)^{-1} \Sigma_0^{1/2} -\sigma \{\Sigma_0^{1/2} (V+2\sigma \Sigma_0)^{-1} \Sigma_0^{1/2}\}^2& O
\end{bmatrix},\\ 
&\{S_2'R^{-1}\}^2
=
\begin{bmatrix}
\dfrac{1}{4}\{\Sigma_0^{1/2} (V+2\sigma \Sigma_0)^{-1} \Sigma_0^{1/2}\}^2-\dfrac{1}{2}\sigma\{\Sigma_0^{1/2} (V+2\sigma \Sigma_0)^{-1} \Sigma_0^{1/2}\}^3~~~~~ O \\ 
O~~~~~ \dfrac{1}{4}\{\Sigma_0^{1/2} (V+2\sigma \Sigma_0)^{-1} \Sigma_0^{1/2}\}^2-\dfrac{1}{2}\sigma\{\Sigma_0^{1/2} (V+2\sigma \Sigma_0)^{-1} \Sigma_0^{1/2}\}^3
\end{bmatrix}, \\
&S_2R^{-1}
=
\begin{bmatrix}
O&\dfrac{1}{2}\Sigma_0^{1/2} V^{-1} \Sigma_0^{1/2}  \\ 
\dfrac{1}{2}\Sigma_0^{1/2} (V+2\sigma \Sigma_0)^{-1} \Sigma_0^{1/2}&O 
\end{bmatrix}, \\
&\{S_2R^{-1}\}^2
=
\begin{bmatrix}
\dfrac{1}{4} \Sigma_0^{1/2} V^{-1} \Sigma_0 (V+2\sigma \Sigma_0)^{-1} \Sigma_0^{1/2}&O  \\ 
O&\dfrac{1}{4}\Sigma_0^{1/2} (V+2\sigma \Sigma_0)^{-1} \Sigma_0 V^{-1} \Sigma_0^{1/2} 
\end{bmatrix}, \\
&S_3'R^{-1}
=
\begin{bmatrix}
\Sigma_0^{1/2} (V+2\sigma\Sigma_0)^{-1} \Sigma_0^{1/2} -2\sigma \{\Sigma_0^{1/2} (V+2\sigma\Sigma_0)^{-1} \Sigma_0^{1/2}\}^2&O  \\ 
O&O 
\end{bmatrix}, \\
&S_4'R^{-1}
=
\begin{bmatrix}
O&O  \\ 
O&\Sigma_0^{1/2} (V+2\sigma\Sigma_0)^{-1} \Sigma_0^{1/2} 
\end{bmatrix}, \\
&S_3'R^{-1} S_4'R^{-1}
=O, \\
&S_1R^{-1} S_1'R^{-1}
=
\begin{bmatrix}
2\sigma^3 \Sigma_0^{1/2} V^{-1} \Sigma_0 (V+2\sigma \Sigma_0)^{-1} \Sigma_0 T_3' \Sigma_0^{1/2} ~~~~~~~~~~~~~~~~~~~~~~~~~~~~~~~~~~~~~~~~~~~~~~~~~~~\\
~-4\sigma^4\Sigma_0^{1/2} V^{-1} \Sigma_0 (V+2\sigma \Sigma_0)^{-1} \Sigma_0 T_3' \Sigma_0 (V+2\sigma \Sigma_0)^{-1} \Sigma_0^{1/2} ~~~~~~~~~~~~~~~~~~~~~~~~~*
\\
~+\sigma^2\Sigma_0^{1/2} V^{-1} \Sigma_0 (V+2\sigma \Sigma_0)^{-1} \Sigma_0^{1/2} ~~~~~~~~~~~~~~~~~~~~~~~~~~~~~~~~~~~~~~~~~~~~~~~~~~~~~~~~\\
~~~~~
-2\sigma^3 \Sigma_0^{1/2} V^{-1} \Sigma_0^{1/2} \{\Sigma_0^{1/2} (V+2\sigma \Sigma_0)^{-1} \Sigma_0^{1/2}\}^2~~~~~~~~~~~~~~~~~~~~~~~~~~~~~~~~~~~~~~~~~~~~~~\\
*~~~~~~~~~~~~~~~~~~~~~~~~~~\sigma^2 \{\Sigma_0^{1/2}(V+2\sigma \Sigma_0)^{-1} \Sigma_0\}^2+\sigma^2 \Sigma_0^{1/2} V^{-1} \Sigma_0 T_3' \Sigma_0^{1/2}
\end{bmatrix}, \\
&S_1R^{-1} \{S_2'R^{-1}\}^2\\
&=
\begin{bmatrix}
-\dfrac{1}{2}\sigma^2 \Sigma_0^{1/2} V^{-1}\Sigma_0^{1/2}\{ \Sigma_0^{1/2}(V+2\sigma \Sigma_0)^{-1} \Sigma_0^{1/2}\}^3 +\sigma^3 \Sigma_0^{1/2} V^{-1}\Sigma_0^{1/2}\{ \Sigma_0^{1/2}(V+2\sigma \Sigma_0)^{-1} \Sigma_0^{1/2}\}^4~~~~~* \\ 
*~~~~~ -\dfrac{1}{4}\sigma \Sigma_0^{1/2} V^{-1} \Sigma_0^{1/2} \{\Sigma_0^{1/2} (V+2\sigma \Sigma_0)^{-1} \Sigma_0^{1/2}\}^2+\dfrac{1}{2}\sigma^2 \Sigma_0^{1/2} V^{-1} \Sigma_0^{1/2} \{\Sigma_0^{1/2} (V+2\sigma \Sigma_0)^{-1} \Sigma_0^{1/2}\}^3 
\end{bmatrix}, \\
&S_1'R^{-1}\{S_2R^{-1}\}^2\\
&=
\begin{bmatrix}
-\dfrac{1}{4}\sigma \Sigma_0^{1/2} T_3' \Sigma_0 V^{-1} \Sigma_0 (V+2\sigma \Sigma_0)^{-1} \Sigma_0^{1/2}+\dfrac{1}{2}\sigma^2 \Sigma_0^{1/2} T_3' \Sigma_0 (V+2\sigma \Sigma_0)^{-1} \Sigma_0 V^{-1} \Sigma_0 (V+2\sigma \Sigma_0)^{-1} \Sigma_0^{-1/2} ~~*  \\ 
*~~~~~~~~~~~~~~~~~~~~~~~~~~~~~~~~~~~~~~~~-\dfrac{1}{4}\sigma \Sigma_0^{1/2} T_3' \Sigma_0 (V+2\sigma \Sigma_0)^{-1} \Sigma_0 V^{-1} \Sigma_0^{1/2}
\end{bmatrix}, \\
&S_1R^{-1} S_4'R^{-1}
=
\begin{bmatrix}
O&\sigma \Sigma_0^{1/2} V^{-1} \Sigma_0(V+2\sigma \Sigma_0)^{-1} \Sigma_0^{1/2}  \\ 
O& -\sigma \Sigma_0^{1/2} V^{-1} \Sigma_0(V+2\sigma \Sigma_0)^{-1} \Sigma_0^{1/2}  
\end{bmatrix}, \\
&S_1R^{-1} S_3'R^{-1}\\
&=
\begin{bmatrix}
-2\sigma^2 \Sigma_0^{1/2} V^{-1} \Sigma_0^{1/2} \{\Sigma_0^{1/2} (V+2\sigma \Sigma_0)^{-1} \Sigma_0\}^2 +4\sigma^3 \Sigma_0^{1/2} V^{-1} \Sigma_0^{1/2} \{\Sigma_0^{1/2} (V+2\sigma \Sigma_0)^{-1} \Sigma_0\}^3&O  \\ 
\sigma \{\Sigma_0^{1/2} (V+2\sigma \Sigma_0)^{-1} \Sigma_0\}^2-2\sigma^2 \{\Sigma_0^{1/2} (V+2\sigma \Sigma_0)^{-1} \Sigma_0\}^3& O
\end{bmatrix}, \\
&S_2R^{-1} S_2'R^{-1}
=
\begin{bmatrix}
\dfrac{1}{4} \Sigma_0^{1/2} V^{-1} \Sigma_0 (V+2\sigma \Sigma_0)^{-1} \Sigma_0^{1/2}-\dfrac{1}{2}\sigma \Sigma_0^{1/2} V^{-1} \Sigma_0^{1/2} \{  \Sigma_0^{1/2} (V+2\sigma \Sigma_0)^{-1} \Sigma_0^{1/2}\}^2~~~~~O  \\ 
O~~~~~ \dfrac{1}{4} \{\Sigma_0^{1/2}  (V+2\sigma \Sigma_0)^{-1} \Sigma_0^{1/2}\}^2
\end{bmatrix}, \\
&\{S_2R^{-1}\}^2 \{S_2'R^{-1}\}^2\\
&=
\begin{bmatrix}
\dfrac{1}{16} \Sigma_0^{1/2}V^{-1} \Sigma_0^{1/2} \{\Sigma_0^{1/2} (V+2\sigma \Sigma_0)^{-1} \Sigma_0^{1/2}\}^3~~~~~~~~~~~~~~~~~~~~~~~~~~~~~~~~~~~~~~~~~~~~~O\\ ~~~~~-\dfrac{1}{8}\sigma\Sigma_0^{1/2}V^{-1} \Sigma_0^{1/2} \{\Sigma_0^{1/2} (V+2\sigma \Sigma_0)^{-1} \Sigma_0^{1/2}\}^4 ~~~~~~~~~~~~~~~~~~~~~~~~~~~~~~~~~~~~~~~~~~~~~\\ 
~~~~~~~~~~~~~~~~~~ \dfrac{1}{16} \Sigma_0^{1/2} (V+2\sigma \Sigma_0)^{-1} \Sigma_0 V^{-1} \Sigma_0^{1/2} \{\Sigma_0^{1/2} (V+2\sigma \Sigma_0)^{-1} \Sigma_0^{1/2}\}^2 \\
O~~~~~~~~~~~~~~~~~~-\dfrac{1}{8}\sigma \Sigma_0^{1/2} (V+2\sigma \Sigma_0)^{-1} \Sigma_0 V^{-1} \Sigma_0^{1/2} \{\Sigma_0^{1/2} (V+2\sigma \Sigma_0)^{-1} \Sigma_0^{1/2}\}^3
\end{bmatrix}, \\
&\{S_2R^{-1} S_2'R^{-1}\}^2\\
&=
\begin{bmatrix}
\left\{\dfrac{1}{4} \Sigma_0^{1/2} V^{-1} \Sigma_0 (V+2\sigma \Sigma_0)^{-1} \Sigma_0^{1/2}-\dfrac{1}{2}\sigma \Sigma_0^{1/2} V^{-1} \Sigma_0^{1/2} \{  \Sigma_0^{1/2} (V+2\sigma \Sigma_0)^{-1} \Sigma_0^{1/2}\}^2\right\}^2~~~~~O  \\ 
O~~~~~ \dfrac{1}{16} \{\Sigma_0^{1/2}  (V+2\sigma \Sigma_0)^{-1} \Sigma_0^{1/2}\}^4
\end{bmatrix}, \\
&\{S_2R^{-1}\}^2 S_4'R^{-1}
=
\begin{bmatrix}
O&O  \\ 
O&  \dfrac{1}{4}\Sigma_0^{1/2} (V+2\sigma \Sigma_0)^{-1} \Sigma_0 V^{-1} \Sigma_0 (V+2\sigma \Sigma_0)^{-1} \Sigma_0^{1/2}
\end{bmatrix}, \\
&\{S_2R^{-1}\}^2 S_3'R^{-1}\\
&=
\begin{bmatrix}
\dfrac{1}{4} \Sigma_0^{1/2} V^{-1} \Sigma_0^{1/2} \{\Sigma_0^{1/2} (V+2\sigma \Sigma_0)^{-1} \Sigma_0^{1/2}\}^2-\dfrac{1}{2}\sigma \Sigma_0^{1/2} V^{-1} \Sigma_0^{1/2} \{\Sigma_0^{1/2} (V+2\sigma \Sigma_0)^{-1} \Sigma_0^{1/2}\}^3 &O  \\ 
O& O
\end{bmatrix}, \\
&S_2R^{-1} S_3' R^{-1}
=
\begin{bmatrix}
O&O  \\ 
\dfrac{1}{2} \{\Sigma_0^{1/2} (V+2\sigma \Sigma_0)^{-1} \Sigma_0^{1/2}\}^2-\sigma \{\Sigma_0^{1/2} (V+2\sigma \Sigma_0)^{-1} \Sigma_0^{1/2}\}^3& O
\end{bmatrix}, \\
&S_2R^{-1} S_4'R^{-1}
=
\begin{bmatrix}
O&\dfrac{1}{2}\Sigma_0^{1/2} V^{-1} \Sigma_0 (V+2\sigma \Sigma_0)^{-1} \Sigma_0^{1/2}  \\ 
O&O 
\end{bmatrix}
\end{align*}
and
\begin{align*}
&S_2R^{-1} S_3'R^{-1} S_2R^{-1} S_4'R^{-1}\\
&=
\begin{bmatrix}
O&O  \\ 
O&\dfrac{1}{4}\{\Sigma_0^{1/2} (V+2\sigma \Sigma_0)^{-1} \Sigma_0^{1/2}\}^2 \Sigma_0^{1/2} V^{-1} \Sigma_0 (V+2\sigma \Sigma_0)^{-1} \Sigma_0^{1/2}\\
&~~~~~-\dfrac{1}{2}\sigma \{\Sigma_0^{1/2} (V+2\sigma \Sigma_0)^{-1} \Sigma_0^{1/2}\}^3 \Sigma_0^{1/2} V^{-1} \Sigma_0 (V+2\sigma \Sigma_0)^{-1} \Sigma_0^{1/2} 
\end{bmatrix},
\end{align*}
where we omit the calculations of the off-diagonal block for matrix in the left hand side such that it appears only in the trace in the final calculations $Q_{1}$, $Q_{2}$, $Q_{3}$ and $Q_{4}$. 

We now start to evaluate each term in (\ref{MATHCAL_F24}).

For a necessary scalar appeared firstly, it holds that 
\begin{align*}
&(1+\sigma \text{tr}[V^{-1} \Sigma_0]) C_1'\\
&=1+2\sigma \text{tr}[(V+2\sigma \Sigma_0)^{-1} \Sigma_0]
+\sigma \text{tr}[V^{-1} \Sigma_0]
+\sigma^2\{\text{tr}[(V+2\sigma \Sigma_0)^{-1} \Sigma_0]\}^2+2\sigma^2\text{tr}[\{(V+2\sigma \Sigma_0)^{-1} \Sigma_0\}^2]\\
&~~~~~
+2\sigma^2 \text{tr}[V^{-1}\Sigma_0]\text{tr}[(V+2\sigma \Sigma_0)^{-1} \Sigma_0]
+\sigma^3 \text{tr}[V^{-1} \Sigma_0] \{\text{tr}[(V+2\sigma \Sigma_0)^{-1} \Sigma_0]\}^2\\
&~~~~~
+2\sigma^3 \text{tr}[V^{-1}\Sigma_0] \text{tr}[\{(V+2\sigma \Sigma_0)^{-1} \Sigma_0\}^2].
\end{align*}
Using above expressions previously obtained, we see that terms related to $Q_{1}$ can be calculated as
\begin{align*}
&Q_1(S_1'R^{-1})\\
&=\text{tr}[S_1'R^{-1}]\\
&=-2\sigma \text{tr}[T_3' \Sigma_0] +2\sigma^2 \text{tr}[T_3' \Sigma_0 (V+2\sigma \Sigma_0)^{-1} \Sigma_0] \\
&=-2\sigma \text{tr}[
(V+2\sigma \Sigma_0)^{-1} \Sigma_0+\sigma \text{tr}[(V+2\sigma \Sigma_0)^{-1} \Sigma_0] (V+2\sigma \Sigma_0)^{-1} \Sigma_0+2\sigma\{(V+2\sigma \Sigma_0)^{-1} \Sigma_0\}^2
]\\
&~~~~~+2\sigma^2 \text{tr}[\{
(V+2\sigma \Sigma_0)^{-1} \Sigma_0+\sigma \text{tr}[(V+2\sigma \Sigma_0)^{-1} \Sigma_0] (V+2\sigma \Sigma_0)^{-1} \Sigma_0+2\sigma\{(V+2\sigma \Sigma_0)^{-1} \Sigma_0\}^2
\}\\
&~~~~~~~~~~\times(V+2\sigma \Sigma_0)^{-1} \Sigma_0]\\
&=-2\sigma \text{tr}[(V+2\sigma \Sigma_0)^{-1} \Sigma_0]
-2\sigma^2 \{\text{tr}[(V+2\sigma \Sigma_0)^{-1} \Sigma_0]\}^2
-4\sigma^2 \text{tr}[\{(V+2\sigma \Sigma_0)^{-1} \Sigma_0\}^2]\\
&~~~~~+2\sigma^2 \text{tr}[\{(V+2\sigma \Sigma_0)^{-1} \Sigma_0\}^2]
+2\sigma^3 \text{tr}[(V+2\sigma \Sigma_0)^{-1} \Sigma_0] \text{tr}[\{(V+2\sigma \Sigma_0)^{-1} \Sigma_0\}^2]\\
&~~~~~+4\sigma^3 \text{tr}[\{(V+2\sigma \Sigma_0)^{-1} \Sigma_0\}^3]\\
&=-2\sigma \text{tr}[(V+2\sigma \Sigma_0)^{-1} \Sigma_0]
-2\sigma^2\{\text{tr}[(V+2\sigma \Sigma_0)^{-1} \Sigma_0]\}^2
-2\sigma^2 \text{tr}[\{(V+2\sigma \Sigma_0)^{-1} \Sigma_0\}^2]\\
&~~~~~+2\sigma^3 \text{tr}[(V+2\sigma \Sigma_0)^{-1} \Sigma_0]\text{tr}[\{(V+2\sigma \Sigma_0)^{-1} \Sigma_0\}^2]
+4\sigma^3 \text{tr}[\{(V+2\sigma \Sigma_0)^{-1} \Sigma_0\}^3],\\
&(1+\sigma \text{tr}[V^{-1} \Sigma_0]) Q_1(S_1'R^{-1})\\
&=-2\sigma \text{tr}[(V+2\sigma \Sigma_0)^{-1} \Sigma_0]
-2\sigma^2 \{\text{tr}[(V+2\sigma \Sigma_0)^{-1} \Sigma_0]\}^2
-2\sigma^2 \text{tr}[\{(V+2\sigma \Sigma_0)^{-1} \Sigma_0\}^2]\\
&~~~~~
-2\sigma^2 \text{tr}[V^{-1} \Sigma_0] \text{tr}[(V+2\sigma \Sigma_0)^{-1} \Sigma_0]
+2\sigma^3 \text{tr}[(V+2\sigma \Sigma_0)^{-1} \Sigma_0]\text{tr}[\{(V+2\sigma \Sigma_0)^{-1} \Sigma_0\}^2]\\
&~~~~~
+4\sigma^3 \text{tr}[\{(V+2\sigma \Sigma_0)^{-1} \Sigma_0\}^3]
-2\sigma^3 \text{tr}[V^{-1} \Sigma_0] \{\text{tr}[(V+2\sigma \Sigma_0)^{-1} \Sigma_0]\}^2\\
&~~~~~
-2\sigma^3 \text{tr}[V^{-1} \Sigma_0] \text{tr}[\{(V+2\sigma \Sigma_0)^{-1} \Sigma_0\}^2]
+2\sigma^4 \text{tr}[V^{-1} \Sigma_0 ] \text{tr}[(V+2\sigma \Sigma_0)^{-1} \Sigma_0] \text{tr}[\{(V+2\sigma \Sigma_0)^{-1} \Sigma_0\}^2]\\
&~~~~~
+4\sigma^4 \text{tr}[V^{-1} \Sigma_0] \text{tr}[\{(V+2\sigma \Sigma_0)^{-1} \Sigma_0\}^3],\\
&Q_1(S_1R^{-1})\\
&=\text{tr}[S_1R^{-1}]\\
&=-2\sigma^2\text{tr}[V^{-1}\Sigma_0(V+2\sigma \Sigma_0)^{-1} \Sigma_0] 
-\sigma \text{tr}[V^{-1} \Sigma_0]\\
&=-2\sigma \text{tr}[V^{-1} \Sigma_0] 
+\sigma \text{tr}[(V+2\sigma \Sigma_0)^{-1} \Sigma_0]
\end{align*}
and
\begin{align*}
&C_1'Q_1(S_1R^{-1})\\
&=\{
1+2\sigma \text{tr}[(V+2\sigma \Sigma_0)^{-1} \Sigma_0]+\sigma^2 \{\text{tr}[(V+2\sigma \Sigma_0)^{-1} \Sigma_0]\}^2+2\sigma^2 \text{tr}[\{(V+2\sigma \Sigma_0)^{-1} \Sigma_0\}^2]
\}\\
&~~~~~\times\{
-2\sigma \text{tr}[V^{-1} \Sigma_0] 
+\sigma \text{tr}[(V+2\sigma \Sigma_0)^{-1} \Sigma_0]
\}\\
&=-2\sigma \text{tr}[V^{-1} \Sigma_0] 
+\sigma \text{tr}[(V+2\sigma \Sigma_0)^{-1} \Sigma_0]
-4\sigma^2 \text{tr}[V^{-1} \Sigma_0] \text{tr}[(V+2\sigma \Sigma_0)^{-1} \Sigma_0]
\\
&~~~~~
+2\sigma^2 \{\text{tr}[(V+2\sigma \Sigma_0)^{-1} \Sigma_0]\}^2
-2\sigma^3 \text{tr}[V^{-1} \Sigma_0] \{\text{tr}[(V+2\sigma \Sigma_0)^{-1} \Sigma_0]\}^2
+\sigma^3 \{\text{tr}[(V+2\sigma \Sigma_0)^{-1} \Sigma_0]\}^3\\
&~~~~~
-4\sigma^3 \text{tr}[V^{-1} \Sigma_0] \text{tr}[\{(V+2\sigma \Sigma_0)^{-1} \Sigma_0\}^2]+2\sigma^3 \text{tr}[(V+2\sigma \Sigma_0)^{-1} \Sigma_0] \text{tr}[\{(V+2\sigma \Sigma_0)^{-1} \Sigma_0\}^2].
\end{align*}
We have for the terms related to $Q_{2}$ as
\begin{align*}
&Q_2(S_2'R^{-1} ,S_2'R^{-1})\\
&=\{\text{tr}[S_2'R^{-1}]\}^2+2\text{tr}[\{S_2'R^{-1}\}^2]\\
&=\text{tr}[\{(V+2\sigma \Sigma_0)^{-1} \Sigma_0\}^2]-2\sigma \text{tr}[\{(V+2\sigma \Sigma_0)^{-1} \Sigma_0\}^3],\\
&2\sigma^2 (1+\sigma \text{tr}[V^{-1} \Sigma_0]) Q_2(S_2'R^{-1} ,S_2'R^{-1})\\
&=2\sigma^2 \text{tr}[\{(V+2\sigma \Sigma_0)^{-1} \Sigma_0\}^2] 
-4\sigma^3 \text{tr}[\{(V+2\sigma \Sigma_0)^{-1} \Sigma_0\}^3]\\
&~~~~~+2\sigma^3 \text{tr}[V^{-1} \Sigma_0] \text{tr}[\{(V+2\sigma \Sigma_0)^{-1} \Sigma_0\}^2]
-4\sigma^4 \text{tr}[V^{-1} \Sigma_0] \text{tr}[\{(V+2\sigma \Sigma_0)^{-1} \Sigma_0\}^3],\\
&Q_2(S_2R^{-1},S_2R^{-1})\\
&=\{\text{tr}[S_2R^{-1}]\}^2+2\text{tr}[\{S_2R^{-1}\}^2]\\
&=\text{tr}[V^{-1} \Sigma_0 (V+2\sigma \Sigma_0)^{-1} \Sigma_0]\\
&=\frac{1}{2\sigma} \text{tr}[V^{-1} \Sigma_0]-\frac{1}{2\sigma} \text{tr}[(V+2\sigma \Sigma_0)^{-1} \Sigma_0],\\
&2\sigma^2 C_1' Q_2(S_2R^{-1},S_2R^{-1})\\
&=\sigma \{
1+2\sigma \text{tr}[(V+2\sigma \Sigma_0)^{-1} \Sigma_0]+\sigma^2 \{\text{tr}[(V+2\sigma \Sigma_0)^{-1} \Sigma_0]\}^2+2\sigma^2 \text{tr}[\{(V+2\sigma \Sigma_0)^{-1} \Sigma_0\}^2]
\}\\
&~~~~~\times\{
\text{tr}[V^{-1} \Sigma_0] -\text{tr}[(V+2\sigma \Sigma_0)^{-1} \Sigma_0]
\}\\
&=\sigma \text{tr}[V^{-1} \Sigma_0] -\sigma \text{tr}[(V+2\sigma \Sigma_0)^{-1} \Sigma_0]
+2\sigma^2 \text{tr}[V^{-1}\Sigma_0] \text{tr}[(V+2\sigma \Sigma_0)^{-1} \Sigma_0]
-2\sigma^2 \{\text{tr}[(V+2\sigma \Sigma_0)^{-1} \Sigma_0]\}^2\\
&~~~~~
+\sigma^3 \text{tr}[V^{-1} \Sigma_0] \{\text{tr}[(V+2\sigma \Sigma_0)^{-1} \Sigma_0]\}^2
-\sigma^3 \{\text{tr}[(V+2\sigma \Sigma_0)^{-1} \Sigma_0]\}^3\\
&~~~~~
+2\sigma^3 \text{tr}[V^{-1} \Sigma_0] \text{tr}[\{(V+2\sigma \Sigma_0)^{-1} \Sigma_0\}^2]
-2\sigma^3 \text{tr}[(V+2\sigma \Sigma_0)^{-1} \Sigma_0] \text{tr}[\{(V+2\sigma \Sigma_0)^{-1} \Sigma_0\}^2],\\
&Q_2(S_3'R^{-1},S_4'R^{-1})\\
&=\text{tr}[S_3'R^{-1}] \text{tr}[S_4'R^{-1}] +2\text{tr}[S_3'R^{-1} S_4'R^{-1}]\\
&=\{\text{tr}[(V+2\sigma \Sigma_0)^{-1} \Sigma_0]-2\sigma \text{tr}[\{(V+2\sigma \Sigma_0)^{-1} \Sigma_0\}^2]\} \text{tr}[(V+2\sigma \Sigma_0)^{-1} \Sigma_0]\\
&=\{\text{tr}[(V+2\sigma \Sigma_0)^{-1} \Sigma_0]\}^2
-2\sigma \text{tr}[(V+2\sigma \Sigma_0)^{-1} \Sigma_0] \text{tr}[\{(V+2\sigma \Sigma_0)^{-1} \Sigma_0\}^2]
\end{align*}
and
\begin{align*}
&\sigma^2 (1+\sigma \text{tr}[V^{-1} \Sigma_0]) Q_2(S_3'R^{-1},S_4'R^{-1})\\
&=\sigma^2 \{\text{tr}[(V+2\sigma \Sigma_0)^{-1} \Sigma_0]\}^2
-2\sigma^3 \text{tr}[(V+2\sigma \Sigma_0)^{-1} \Sigma_0] \text{tr}[\{(V+2\sigma \Sigma_0)^{-1} \Sigma_0\}^2]\\
&~~~~~+\sigma^3 \text{tr}[V^{-1} \Sigma_0] \{\text{tr}[(V+2\sigma \Sigma_0)^{-1} \Sigma_0]\}^2
-2\sigma^4 \text{tr}[V^{-1} \Sigma_0] \text{tr}[(V+2\sigma \Sigma_0)^{-1} \Sigma_0] \text{tr}[\{(V+2\sigma \Sigma_0)^{-1} \Sigma_0\}^2].
\end{align*}
Note that the term $Q_2(S_1R^{-1},S_1'R^{-1})$ is decomposed into
\begin{align*}
&Q_2(S_1R^{-1},S_1'R^{-1})\\
&=\text{tr}[S_1R^{-1}] \text{tr}[S_1'R^{-1}] +2\text{tr}[S_1R^{-1} S_1'R^{-1}]
\end{align*}
and we try to obtain useful expression for each term. 
It is easy to confirm that
\begin{align*}
&\text{tr}[S_1R^{-1}] \text{tr}[S_1'R^{-1}]\\
&=\{
-2\sigma \text{tr}[V^{-1} \Sigma_0] 
+\sigma \text{tr}[(V+2\sigma \Sigma_0)^{-1} \Sigma_0]
\}\\
&~~~~~\times 
\{
-2\sigma \text{tr}[(V+2\sigma \Sigma_0)^{-1} \Sigma_0]
-2\sigma^2\{\text{tr}[(V+2\sigma \Sigma_0)^{-1} \Sigma_0]\}^2
-2\sigma^2 \text{tr}[\{(V+2\sigma \Sigma_0)^{-1} \Sigma_0\}^2]\\
&~~~~~~~~~~+2\sigma^3 \text{tr}[(V+2\sigma \Sigma_0)^{-1} \Sigma_0]\text{tr}[\{(V+2\sigma \Sigma_0)^{-1} \Sigma_0\}^2]
+4\sigma^3 \text{tr}[\{(V+2\sigma \Sigma_0)^{-1} \Sigma_0\}^3]
\}\\
&=4\sigma^2 \text{tr}[V^{-1} \Sigma_0] \text{tr}[(V+2\sigma \Sigma_0)^{-1} \Sigma_0]
-2\sigma^2 \{\text{tr}[(V+2\sigma \Sigma_0)^{-1} \Sigma_0]\}^2
+4\sigma^3\text{tr}[V^{-1} \Sigma_0] \{\text{tr}[(V+2\sigma \Sigma_0)^{-1} \Sigma_0]\}^2\\
&~~~~~
-2\sigma^3 \{\text{tr}[(V+2\sigma \Sigma_0)^{-1} \Sigma_0]\}^3
+4\sigma^3 \text{tr}[V^{-1} \Sigma_0] \text{tr}[\{(V+2\sigma \Sigma_0)^{-1} \Sigma_0\}^2]\\
&~~~~~
-2\sigma^3 \text{tr}[(V+2\sigma \Sigma_0)^{-1} \Sigma_0] \text{tr}[\{(V+2\sigma \Sigma_0)^{-1} \Sigma_0\}^2]\\
&~~~~~
-4\sigma^4 \text{tr}[V^{-1} \Sigma_0] \text{tr}[(V+2\sigma \Sigma_0)^{-1} \Sigma_0] \text{tr}[\{(V+2\sigma \Sigma_0)^{-1} \Sigma_0\}^2]\\
&~~~~~
+2\sigma^4\{\text{tr}[(V+2\sigma \Sigma_0)^{-1} \Sigma_0]\}^2 \text{tr}[\{(V+2\sigma \Sigma_0)^{-1} \Sigma_0\}^2]
-8\sigma^4 \text{tr}[V^{-1} \Sigma_0] \text{tr}[\{(V+2\sigma \Sigma_0)^{-1} \Sigma_0\}^3] \\
&~~~~~
+4\sigma^4 \text{tr}[(V+2\sigma \Sigma_0)^{-1} \Sigma_0] \text{tr}[\{(V+2\sigma \Sigma_0)^{-1} \Sigma_0\}^3].
\end{align*}
On the other hand, it follows that
\begin{align*}
&\text{tr}[S_1R^{-1} S_1'R^{-1}]\\
&=2\sigma^3 \text{tr}[T_3' \Sigma_0 V^{-1} \Sigma_0 (V+2\sigma \Sigma_0)^{-1} \Sigma_0]
-4\sigma^4 \text{tr}[T_3' \Sigma_0 (V+2\sigma \Sigma_0)^{-1} \Sigma_0 V^{-1} \Sigma_0 (V+2\sigma \Sigma_0)^{-1} \Sigma_0]\\
&~~~~~+\sigma^2\text{tr}[V^{-1} \Sigma_0 (V+2\sigma \Sigma_0)^{-1} \Sigma_0]
-2\sigma^3 \text{tr}[V^{-1} \Sigma_0 \{(V+2\sigma \Sigma_0)^{-1} \Sigma_0\}^2]\\
&~~~~~+\sigma^2 \text{tr}[\{(V+2\sigma \Sigma_0)^{-1} \Sigma_0\}^2] +\sigma^2 \text{tr}[T_3' \Sigma_0V^{-1} \Sigma_0]\\
&=2\sigma^2 \text{tr}[T_3'\Sigma_0 V^{-1} \Sigma_0] 
-\sigma^2 \text{tr}[T_3'\Sigma_0 (V+2\sigma \Sigma_0)^{-1} \Sigma_0]
-2\sigma^3 \text{tr}[T_3'\Sigma_0 (V+2\sigma \Sigma_0)^{-1} \Sigma_0 V^{-1} \Sigma_0]\\
&~~~~~
+2\sigma^3 \text{tr}[T_3'\Sigma_0 \{(V+2\sigma \Sigma_0)^{-1} \Sigma_0\}^2]
+\frac{1}{2}\sigma \text{tr}[V^{-1} \Sigma_0]
-\frac{1}{2} \sigma \text{tr}[(V+2\sigma \Sigma_0)^{-1} \Sigma_0]\\
&~~~~~
-\sigma^2 \text{tr}[V^{-1} \Sigma_0 (V+2\sigma \Sigma_0)^{-1} \Sigma_0]
+2\sigma^2 \text{tr}[\{(V+2\sigma \Sigma_0)^{-1} \Sigma_0\}^2]\\
&=2\sigma^2 \text{tr}[T_3'\Sigma_0 V^{-1} \Sigma_0] 
-\sigma^2 \text{tr}[T_3'\Sigma_0 (V+2\sigma \Sigma_0)^{-1} \Sigma_0]
-2\sigma^3 \text{tr}[T_3'\Sigma_0 (V+2\sigma \Sigma_0)^{-1} \Sigma_0 V^{-1} \Sigma_0]\\
&~~~~~
+2\sigma^3 \text{tr}[T_3'\Sigma_0 \{(V+2\sigma \Sigma_0)^{-1} \Sigma_0\}^2]
+2\sigma^2 \text{tr}[\{(V+2\sigma \Sigma_0)^{-1} \Sigma_0\}^2].
\end{align*}
Here we need to proceed calculations for the terms involving the matrix $T'_{3}=T_{3}^{(4)}$ separately.
We have from long but straightforward evaluations that 
\begin{align*}
&2\sigma^2 \text{tr}[T_3'\Sigma_0 V^{-1} \Sigma_0]\\
&=2\sigma^2 \text{tr}[\{
(V+2\sigma \Sigma_0)^{-1} \Sigma_0+\sigma \text{tr}[(V+2\sigma \Sigma_0)^{-1} \Sigma_0] (V+2\sigma \Sigma_0)^{-1} \Sigma_0+2\sigma\{(V+2\sigma \Sigma_0)^{-1} \Sigma_0\}^2\}\\
&~~~~~\times V^{-1} \Sigma_0]\\
&=2\sigma^2 \text{tr}[V^{-1} \Sigma_0 (V+2\sigma \Sigma_0)^{-1} \Sigma_0]
+2\sigma^3 \text{tr}[(V+2\sigma \Sigma_0)^{-1} \Sigma_0] \text{tr}[V^{-1} \Sigma_0(V+2\sigma \Sigma_0)^{-1} \Sigma_0]\\
&~~~~~
+4\sigma^3 \text{tr}[V^{-1} \Sigma_0 \{(V+2\sigma \Sigma_0)^{-1} \Sigma_0\}^2]\\
&=\sigma \text{tr}[V^{-1} \Sigma_0]
-\sigma \text{tr}[(V+2\sigma \Sigma_0)^{-1} \Sigma_0]
+\sigma^2 \text{tr}[V^{-1} \Sigma_0] \text{tr}[(V+2\sigma \Sigma_0)^{-1} \Sigma_0]
\\
&~~~~~
-\sigma^2 \{\text{tr}[(V+2\sigma \Sigma_0)^{-1} \Sigma_0]\}^2
+2\sigma^2 \text{tr}[V^{-1} \Sigma_0 (V+2\sigma \Sigma_0)^{-1} \Sigma_0]
-2\sigma^2 \text{tr}[\{(V+2\sigma \Sigma_0)^{-1} \Sigma_0\}^2]\\
&=2\sigma \text{tr}[V^{-1} \Sigma_0]
-2\sigma \text{tr}[(V+2\sigma \Sigma_0)^{-1} \Sigma_0]
+\sigma^2 \text{tr}[V^{-1} \Sigma_0] \text{tr}[(V+2\sigma \Sigma_0)^{-1} \Sigma_0]\\
&~~~~~-\sigma^2 \{\text{tr}[(V+2\sigma \Sigma_0)^{-1} \Sigma_0]\}^2
-2\sigma^2 \text{tr}[\{(V+2\sigma \Sigma_0)^{-1} \Sigma_0\}^2],\\
&-\sigma^2 \text{tr}[T_3' \Sigma_0 (V+2\sigma \Sigma_0)^{-1} \Sigma_0]\\
&=-\sigma^2 \text{tr}[\{
(V+2\sigma \Sigma_0)^{-1} \Sigma_0+\sigma \text{tr}[(V+2\sigma \Sigma_0)^{-1} \Sigma_0] (V+2\sigma \Sigma_0)^{-1} \Sigma_0+2\sigma\{(V+2\sigma \Sigma_0)^{-1} \Sigma_0\}^2\}\\
&~~~~~\times (V+2\sigma \Sigma_0)^{-1} \Sigma_0]\\
&=-\sigma^2 \text{tr}[\{(V+2\sigma \Sigma_0)^{-1} \Sigma_0\}^2]
-\sigma^3 \text{tr}[(V+2\sigma \Sigma_0)^{-1} \Sigma_0] \text{tr}[\{(V+2\sigma \Sigma_0)^{-1} \Sigma_0\}^2]\\
&~~~~~
-2\sigma^3 \text{tr}[\{(V+2\sigma \Sigma_0)^{-1} \Sigma_0\}^3],\\
&-2\sigma^3 \text{tr}[T_3' \Sigma_0 (V+2\sigma \Sigma_0)^{-1} \Sigma_0 V^{-1} \Sigma_0]\\
&=-2\sigma^3  \text{tr}[\{
(V+2\sigma \Sigma_0)^{-1} \Sigma_0+\sigma \text{tr}[(V+2\sigma \Sigma_0)^{-1} \Sigma_0] (V+2\sigma \Sigma_0)^{-1} \Sigma_0+2\sigma\{(V+2\sigma \Sigma_0)^{-1} \Sigma_0\}^2\}\\
&~~~~~\times (V+2\sigma \Sigma_0)^{-1} \Sigma_0 V^{-1} \Sigma_0]\\
&=-2\sigma^3 \text{tr}[V^{-1}\Sigma_0 \{(V+2\sigma \Sigma_0)^{-1} \Sigma_0\}^2]
-2\sigma^4 \text{tr}[(V+2\sigma \Sigma_0)^{-1} \Sigma_0] \text{tr}[V^{-1}\Sigma_0 \{(V+2\sigma \Sigma_0)^{-1} \Sigma_0\}^2]\\
&~~~~~
-4\sigma^4 \text{tr}[V^{-1}\Sigma_0 \{(V+2\sigma \Sigma_0)^{-1} \Sigma_0\}^3]\\
&=-\sigma^2 \text{tr}[V^{-1}\Sigma_0 (V+2\sigma \Sigma_0)^{-1} \Sigma_0] 
+\sigma^2 \text{tr}[\{(V+2\sigma \Sigma_0)^{-1} \Sigma_0\}^2]\\
&~~~~~
-\sigma^3 \text{tr}[(V+2\sigma \Sigma_0)^{-1} \Sigma_0] \text{tr}[V^{-1}\Sigma_0 (V+2\sigma \Sigma_0)^{-1} \Sigma_0] 
\\
&~~~~~
+\sigma^3 \text{tr}[(V+2\sigma \Sigma_0)^{-1} \Sigma_0]\text{tr}[\{(V+2\sigma \Sigma_0)^{-1} \Sigma_0\}^2]\\
&~~~~~
-2\sigma^3\text{tr}[V^{-1}\Sigma_0\{(V+2\sigma \Sigma_0)^{-1} \Sigma_0\}^2]
+2\sigma^3 \text{tr}[\{(V+2\sigma \Sigma_0)^{-1} \Sigma_0\}^3]\\
&=-\frac{1}{2}\sigma \text{tr}[V^{-1}\Sigma_0]
+\frac{1}{2}\sigma \text{tr}[(V+2\sigma \Sigma_0)^{-1} \Sigma_0]
+2\sigma^2 \text{tr}[\{(V+2\sigma \Sigma_0)^{-1} \Sigma_0\}^2]\\
&~~~~~
-\frac{1}{2}\sigma^2\text{tr}[V^{-1}\Sigma_0] \text{tr}[(V+2\sigma \Sigma_0)^{-1} \Sigma_0]
+\frac{1}{2}\sigma^2 \{\text{tr}[(V+2\sigma \Sigma_0)^{-1} \Sigma_0]\}^2\\
&~~~~~
+\sigma^3\text{tr}[(V+2\sigma \Sigma_0)^{-1} \Sigma_0] \text{tr}[\{(V+2\sigma \Sigma_0)^{-1} \Sigma_0\}^2]\\
&~~~~~
-\sigma^2 \text{tr}[V^{-1}\Sigma_0 (V+2\sigma \Sigma_0)^{-1} \Sigma_0] 
+2\sigma^3 \text{tr}[\{(V+2\sigma \Sigma_0)^{-1} \Sigma_0\}^3]\\
&=-\sigma \text{tr}[V^{-1}\Sigma_0] 
+\sigma \text{tr}[(V+2\sigma \Sigma_0)^{-1} \Sigma_0]
+2\sigma^2 \text{tr}[\{(V+2\sigma \Sigma_0)^{-1} \Sigma_0\}^2]\\
&~~~~~
-\frac{1}{2}\sigma^2 \text{tr}[V^{-1}\Sigma_0] \text{tr}[(V+2\sigma \Sigma_0)^{-1} \Sigma_0]
+\frac{1}{2}\sigma^2 \{\text{tr}[(V+2\sigma \Sigma_0)^{-1} \Sigma_0]\}^2\\
&~~~~~
+\sigma^3 \text{tr}[(V+2\sigma \Sigma_0)^{-1} \Sigma_0] \text{tr}[\{(V+2\sigma \Sigma_0)^{-1} \Sigma_0\}^2]
+2\sigma^3 \text{tr}[\{(V+2\sigma \Sigma_0)^{-1} \Sigma_0\}^3]
\end{align*}
and
\begin{align*}
&2\sigma^3 \text{tr}[T_3' \Sigma_0\{(V+2\sigma \Sigma_0)^{-1} \Sigma_0\}^2]\\
&=2\sigma^3 \text{tr}[\{
(V+2\sigma \Sigma_0)^{-1} \Sigma_0+\sigma \text{tr}[(V+2\sigma \Sigma_0)^{-1} \Sigma_0] (V+2\sigma \Sigma_0)^{-1} \Sigma_0+2\sigma\{(V+2\sigma \Sigma_0)^{-1} \Sigma_0\}^2\}\\
&~~~~~\times \{(V+2\sigma \Sigma_0)^{-1} \Sigma_0 \}^2]\\
&=2\sigma^3 \text{tr}[\{(V+2\sigma \Sigma_0)^{-1} \Sigma_0\}^3] 
+2\sigma^4 \text{tr}[(V+2\sigma \Sigma_0)^{-1} \Sigma_0] \text{tr}[\{(V+2\sigma \Sigma_0)^{-1} \Sigma_0\}^3] \\
&~~~~~+4\sigma^4 \text{tr}[\{(V+2\sigma \Sigma_0)^{-1} \Sigma_0\}^4],
\end{align*}
from which we obtain
\begin{align*}
&\text{tr}[S_1R^{-1}S_1'R^{-1}]\\
&=
\sigma \text{tr}[V^{-1}\Sigma_0]
-\sigma \text{tr}[(V+2\sigma \Sigma_0)^{-1}\Sigma_0] 
+\sigma^2 \text{tr}[\{(V+2\sigma \Sigma_0)^{-1} \Sigma_0\}^2]\\
&~~~~~
+\frac{1}{2}\sigma^2 \text{tr}[V^{-1}\Sigma_0] \text{tr}[(V+2\sigma \Sigma_0)^{-1} \Sigma_0]
-\frac{1}{2}\sigma^2 \{\text{tr}[(V+2\sigma \Sigma_0)^{-1} \Sigma_0]\}^2 \\
&~~~~~
+2\sigma^3 \text{tr}[\{(V+2\sigma \Sigma_0)^{-1} \Sigma_0\}^3]
+2\sigma^4 \text{tr}[(V+2\sigma \Sigma_0)^{-1} \Sigma_0] \text{tr}[\{(V+2\sigma \Sigma_0)^{-1} \Sigma_0\}^3]\\
&~~~~~+4\sigma^4 \text{tr}[\{(V+2\sigma \Sigma_0)^{-1} \Sigma_0\}^4]
\end{align*}
and therefore
\begin{align*}
&Q_2(S_1R^{-1},S_1'R^{-1})\\
&=\text{tr}[S_1R^{-1}] \text{tr}[S_1'R^{-1}] +2\text{tr}[S_1R^{-1} S_1'R^{-1}]\\
&=2\sigma \text{tr}[V^{-1}\Sigma_0]
-2\sigma \text{tr}[(V+2\sigma \Sigma_0)^{-1} \Sigma_0]
+5\sigma^2 \text{tr}[V^{-1}\Sigma_0] \text{tr}[(V+2\sigma \Sigma_0)^{-1} \Sigma_0]
\\
&~~~~~
-3\sigma^2 \{\text{tr}[(V+2\sigma \Sigma_0)^{-1} \Sigma_0]\}^2
+2\sigma^2 \text{tr}[\{(V+2\sigma \Sigma_0)^{-1} \Sigma_0\}^2]
+4\sigma^3 \text{tr}[\{(V+2\sigma \Sigma_0)^{-1} \Sigma_0\}^3]
\\
&~~~~~
+4\sigma^3 \text{tr}[V^{-1}\Sigma_0] \{\text{tr}[(V+2\sigma \Sigma_0)^{-1} \Sigma_0]\}^2
-2\sigma^3\{\text{tr}[(V+2\sigma \Sigma_0)^{-1} \Sigma_0]\}^3
\\
&~~~~~
+4\sigma^3 \text{tr}[V^{-1}\Sigma_0] \text{tr}[\{(V+2\sigma \Sigma_0)^{-1} \Sigma_0\}^2]
-2\sigma^3 \text{tr}[(V+2\sigma \Sigma_0)^{-1} \Sigma_0] \text{tr}[\{(V+2\sigma \Sigma_0)^{-1} \Sigma_0\}^2]\\
&~~~~~
-4\sigma^4 \text{tr}[V^{-1}\Sigma_0] \text{tr}[(V+2\sigma \Sigma_0)^{-1} \Sigma_0] \text{tr}[\{(V+2\sigma \Sigma_0)^{-1} \Sigma_0\}^2]\\
&~~~~~
+2\sigma^4 \{\text{tr}[(V+2\sigma \Sigma_0)^{-1} \Sigma_0]\}^2\text{tr}[\{(V+2\sigma \Sigma_0)^{-1} \Sigma_0\}^2] 
-8\sigma^4 \text{tr}[V^{-1}\Sigma_0] \text{tr}[\{(V+2\sigma \Sigma_0)^{-1} \Sigma_0\}^3]\\
&~~~~~
+8\sigma^4 \text{tr}[(V+2\sigma \Sigma_0)^{-1} \Sigma_0] \text{tr}[\{(V+2\sigma \Sigma_0)^{-1} \Sigma_0\}^3]
+8\sigma^4 \text{tr}[\{(V+2\sigma \Sigma_0)^{-1} \Sigma_0\}^4].
\end{align*}
Next we go to the terms related to $Q_{3}$.
We can easily check that
\begin{align*}
& Q_3(S_1R^{-1} ,S_2' R^{-1} ,S_2'R^{-1})\\
&= \text{tr}[S_1R^{-1}] \{\text{tr}[S_2'R^{-1}]\}^2
+2\text{tr}[S_1R^{-1}] \text{tr}[\{S_2'R^{-1}\}^2]
+4\text{tr}[S_2'R^{-1}] \text{tr}[S_1R^{-1} S_2'R^{-1}] \\
&~~~~~+8 \text{tr}[S_1R^{-1} \{S_2'R^{-1}\}^2]\\
&=2\text{tr}[S_1R^{-1}] \text{tr}[\{S_2'R^{-1}\}^2]
+8 \text{tr}[S_1R^{-1} \{S_2'R^{-1}\}^2]\\
&=\{
-2\sigma \text{tr}[V^{-1} \Sigma_0] +\sigma \text{tr}[(V+2\sigma \Sigma_0)^{-1}\Sigma_0]
\}
\{
\text{tr}[\{(V+2\sigma \Sigma_0)^{-1}\Sigma_0\}^2]-2\sigma \text{tr}[\{(V+2\sigma \Sigma_0)^{-1}\Sigma_0\}^3]
\}\\
&~~~~~-2\sigma \text{tr}[V^{-1}\Sigma_0\{(V+2\sigma \Sigma_0)^{-1}\Sigma_0\}^2]
+8\sigma^3 \text{tr}[V^{-1}\Sigma_0 \{(V+2\sigma \Sigma_0)^{-1}\Sigma_0\}^4]\\
&=-2\sigma \text{tr}[V^{-1} \Sigma_0] \text{tr}[\{(V+2\sigma \Sigma_0)^{-1} \Sigma_0\}^2]
+\sigma \text{tr}[(V+2\sigma \Sigma_0)^{-1} \Sigma_0] \text{tr}[\{(V+2\sigma \Sigma_0)^{-1} \Sigma_0\}^2]\\
&~~~~~
+4\sigma^2 \text{tr}[V^{-1} \Sigma_0] \text{tr}[\{(V+2\sigma \Sigma_0)^{-1} \Sigma_0\}^3]
-2\sigma^2 \text{tr}[(V+2\sigma \Sigma_0)^{-1} \Sigma_0] \text{tr}[\{(V+2\sigma \Sigma_0)^{-1} \Sigma_0\}^3]\\
&~~~~~
-2\sigma \text{tr}[V^{-1}\Sigma_0\{(V+2\sigma \Sigma_0)^{-1} \Sigma_0\}^2] 
+4\sigma^2 \text{tr}[V^{-1}\Sigma_0\{(V+2\sigma \Sigma_0)^{-1} \Sigma_0\}^3]\\
&~~~~~
-4\sigma^2 \text{tr}[\{(V+2\sigma \Sigma_0)^{-1} \Sigma_0\}^4]\\
&=-2\sigma \text{tr}[V^{-1} \Sigma_0] \text{tr}[\{(V+2\sigma \Sigma_0)^{-1} \Sigma_0\}^2]
+\sigma \text{tr}[(V+2\sigma \Sigma_0)^{-1} \Sigma_0] \text{tr}[\{(V+2\sigma \Sigma_0)^{-1} \Sigma_0\}^2]\\
&~~~~~
+4\sigma^2 \text{tr}[V^{-1} \Sigma_0] \text{tr}[\{(V+2\sigma \Sigma_0)^{-1} \Sigma_0\}^3]
-2\sigma^2 \text{tr}[(V+2\sigma \Sigma_0)^{-1} \Sigma_0] \text{tr}[\{(V+2\sigma \Sigma_0)^{-1} \Sigma_0\}^3]\\
&~~~~~
-2\sigma \text{tr}[\{(V+2\sigma \Sigma_0)^{-1} \Sigma_0\}^3]
-4\sigma^2 \text{tr}[\{(V+2\sigma \Sigma_0)^{-1} \Sigma_0\}^4].	
\end{align*}
We see that the factorization
\begin{align*}
&Q_3(S_1'R^{-1},S_2R^{-1},S_2R^{-1})\\
&= \text{tr}[S_1'R^{-1}] \{\text{tr}[S_2R^{-1}]\}^2
+2\text{tr}[S_1'R^{-1}] \text{tr}[\{S_2R^{-1}\}^2]
+4\text{tr}[S_2R^{-1}] \text{tr}[S_1'R^{-1} S_2R^{-1}]\\ 
&~~~~~+8 \text{tr}[S_1'R^{-1} \{S_2R^{-1}\}^2]\\
&=2\text{tr}[S_1'R^{-1}] \text{tr}[\{S_2R^{-1}\}^2]
+8 \text{tr}[S_1'R^{-1} \{S_2R^{-1}\}^2],
\end{align*}
where
\begin{align*}
&2\text{tr}[S_1'R^{-1}] \text{tr}[\{S_2R^{-1}\}^2]\\
&=\Big\{
-2\sigma \text{tr}[(V+2\sigma \Sigma_0)^{-1} \Sigma_0]
-2\sigma^2\{\text{tr}[(V+2\sigma \Sigma_0)^{-1} \Sigma_0]\}^2
-2\sigma^2 \text{tr}[\{(V+2\sigma \Sigma_0)^{-1} \Sigma_0\}^2]\\
&~~~~~+2\sigma^3 \text{tr}[(V+2\sigma \Sigma_0)^{-1} \Sigma_0]\text{tr}[\{(V+2\sigma \Sigma_0)^{-1} \Sigma_0\}^2]
+4\sigma^3 \text{tr}[\{(V+2\sigma \Sigma_0)^{-1} \Sigma_0\}^3]
\Big\}\\
&~~~~~\times \left\{
\frac{1}{2\sigma} \text{tr}[V^{-1} \Sigma_0]-\frac{1}{2\sigma} \text{tr}[(V+2\sigma \Sigma_0)^{-1} \Sigma_0]
\right\}\\
&=-\text{tr}[V^{-1}\Sigma_0] \text{tr}[(V+2\sigma \Sigma_0)^{-1}\Sigma_0] 
+\{\text{tr}[(V+2\sigma \Sigma_0)^{-1}\Sigma_0]\}^2
-\sigma \text{tr}[V^{-1}\Sigma_0] \{\text{tr}[(V+2\sigma \Sigma_0)^{-1}\Sigma_0]\}^2\\
&~~~~~
+\sigma \{\text{tr}[(V+2\sigma \Sigma_0)^{-1}\Sigma_0]\}^3
-\sigma \text{tr}[V^{-1}\Sigma_0] \text{tr}[\{(V+2\sigma \Sigma_0)^{-1}\Sigma_0\}^2] \\
&~~~~~
+\sigma \text{tr}[(V+2\sigma \Sigma_0)^{-1}\Sigma_0] \text{tr}[\{(V+2\sigma \Sigma_0)^{-1}\Sigma_0\}^2]\\
&~~~~~
+\sigma^2 \text{tr}[V^{-1}\Sigma_0] \text{tr}[(V+2\sigma \Sigma_0)^{-1}\Sigma_0] \text{tr}[\{(V+2\sigma \Sigma_0)^{-1}\Sigma_0\}^2]\\
&~~~~~
-\sigma^2 \{\text{tr}[(V+2\sigma \Sigma_0)^{-1}\Sigma_0]\}^2 \text{tr}[\{(V+2\sigma \Sigma_0)^{-1}\Sigma_0\}^2]
+2\sigma^2 \text{tr}[V^{-1}\Sigma_0] \text{tr}[\{(V+2\sigma \Sigma_0)^{-1}\Sigma_0\}^3]\\
&~~~~~
-2\sigma^2 \text{tr}[(V+2\sigma \Sigma_0)^{-1}\Sigma_0] \text{tr}[\{(V+2\sigma \Sigma_0)^{-1}\Sigma_0\}^3]
\end{align*}
and
\begin{align*}
&8\text{tr}[S_1'R^{-1} \{S_2R^{-1}\}^2]\\
&=-2\sigma \text{tr}[T_3' \Sigma_0 V^{-1}\Sigma_0 (V+2\sigma \Sigma_0)^{-1} \Sigma_0]\\
&~~~~~+4\sigma^2 \text{tr}[T_3' \Sigma_0 (V+2\sigma \Sigma_0)^{-1} \Sigma_0 V^{-1} \Sigma_0 (V+2\sigma \Sigma_0)^{-1} \Sigma_0]\\
&~~~~~-2\sigma \text{tr}[T_3' \Sigma_0 (V+2\sigma \Sigma_0)^{-1} \Sigma_0 V^{-1} \Sigma_0]\\
&=-\text{tr}[T_3' \Sigma_0 V^{-1}\Sigma_0] 
+\text{tr}[T_3' \Sigma_0 (V+2\sigma \Sigma_0)^{-1} \Sigma_0]
-2\sigma \text{tr}[T_3' \Sigma_0 \{(V+2\sigma \Sigma_0)^{-1} \Sigma_0\}^2]\\
&=-\frac{1}{\sigma}\text{tr}[V^{-1} \Sigma_0] +\frac{1}{\sigma}\text{tr}[(V+2\sigma \Sigma_0)^{-1} \Sigma_0] 
-\frac{1}{2}\text{tr}[V^{-1}\Sigma_0] \text{tr}[(V+2\sigma \Sigma_0)^{-1} \Sigma_0] \\
&~~~~~
+\frac{1}{2}\{\text{tr}[(V+2\sigma \Sigma_0)^{-1} \Sigma_0]\}^2
+2\text{tr}[\{(V+2\sigma \Sigma_0)^{-1} \Sigma_0\}^2]\\
&~~~~~
+\sigma \text{tr}[(V+2\sigma \Sigma_0)^{-1} \Sigma_0] \text{tr}[\{(V+2\sigma \Sigma_0)^{-1} \Sigma_0\}^2]\\
&~~~~~
-2\sigma^2 \text{tr}[(V+2\sigma \Sigma_0)^{-1} \Sigma_0] \text{tr}[\{(V+2\sigma \Sigma_0)^{-1} \Sigma_0\}^3]
-4\sigma^2 \text{tr}[\{(V+2\sigma \Sigma_0)^{-1} \Sigma_0\}^4],
\end{align*}
so we get
\begin{align*}
&Q_3(S_1'R^{-1},S_2R^{-1},S_2R^{-1})\\
&=-\frac{1}{\sigma} \text{tr}[V^{-1}\Sigma_0]
+\frac{1}{\sigma} \text{tr}[(V+2\sigma \Sigma_0)^{-1} \Sigma_0]
-\frac{3}{2}\text{tr}[V^{-1}\Sigma_0] \text{tr}[(V+2\sigma \Sigma_0)^{-1}\Sigma_0] 
\\
&~~~~~
+\frac{3}{2}\{\text{tr}[(V+2\sigma \Sigma_0)^{-1}\Sigma_0]\}^2
+2\text{tr}[\{(V+2\sigma \Sigma_0)^{-1} \Sigma_0\}^2]
-\sigma \text{tr}[V^{-1}\Sigma_0] \{\text{tr}[(V+2\sigma \Sigma_0)^{-1}\Sigma_0]\}^2\\
&~~~~~
+\sigma \{\text{tr}[(V+2\sigma \Sigma_0)^{-1}\Sigma_0]\}^3
-\sigma \text{tr}[V^{-1}\Sigma_0] \text{tr}[\{(V+2\sigma \Sigma_0)^{-1}\Sigma_0\}^2] \\
&~~~~~
+2\sigma \text{tr}[(V+2\sigma \Sigma_0)^{-1}\Sigma_0] \text{tr}[\{(V+2\sigma \Sigma_0)^{-1}\Sigma_0\}^2]\\
&~~~~~
+\sigma^2 \text{tr}[V^{-1}\Sigma_0] \text{tr}[(V+2\sigma \Sigma_0)^{-1}\Sigma_0] \text{tr}[\{(V+2\sigma \Sigma_0)^{-1}\Sigma_0\}^2]\\
&~~~~~
-\sigma^2 \{\text{tr}[(V+2\sigma \Sigma_0)^{-1}\Sigma_0]\}^2 \text{tr}[\{(V+2\sigma \Sigma_0)^{-1}\Sigma_0\}^2]
+2\sigma^2 \text{tr}[V^{-1}\Sigma_0] \text{tr}[\{(V+2\sigma \Sigma_0)^{-1}\Sigma_0\}^3]\\
&~~~~~
-4\sigma^2 \text{tr}[(V+2\sigma \Sigma_0)^{-1}\Sigma_0] \text{tr}[\{(V+2\sigma \Sigma_0)^{-1}\Sigma_0\}^3]
-4\sigma^2 \text{tr}[\{(V+2\sigma \Sigma_0)^{-1} \Sigma_0\}^4].
\end{align*}
Further it is straightforward that
\begin{align*}
&Q_3(S_1R^{-1},S_3'R^{-1},S_4'R^{-1})\\
&=\text{tr}[S_1R^{-1}] \text{tr}[S_3'R^{-1}] \text{tr}[S_4'R^{-1}] 
+2\text{tr}[S_1R^{-1}] \text{tr}[S_3'R^{-1} S_4'R^{-1}] 
+2\text{tr}[S_3'R^{-1}]  \text{tr}[S_1R^{-1} S_4'R^{-1}] \\
&~~~~~
+2\text{tr}[S_4'R^{-1}] \text{tr}[S_1R^{-1} S_3'R^{-1}]
+8\text{tr}[S_1R^{-1} S_3'R^{-1} S_4'R^{-1}]\\
&=\text{tr}[S_1R^{-1}] \text{tr}[S_3'R^{-1}] \text{tr}[S_4'R^{-1}] 
+2\text{tr}[S_3'R^{-1}]  \text{tr}[S_1R^{-1} S_4'R^{-1}]
+2\text{tr}[S_4'R^{-1}] \text{tr}[S_1R^{-1} S_3'R^{-1}]\\
&=\{
-2\sigma \text{tr}[V^{-1} \Sigma_0] 
+\sigma \text{tr}[(V+2\sigma \Sigma_0)^{-1} \Sigma_0]
\}\\
&~~~~~\times 
\{
\{\text{tr}[(V+2\sigma \Sigma_0)^{-1} \Sigma_0]\}^2
-2\sigma \text{tr}[(V+2\sigma \Sigma_0)^{-1} \Sigma_0] \text{tr}[\{(V+2\sigma \Sigma_0)^{-1} \Sigma_0\}^2]
\}\\
&~~~~~-2\sigma\{
\text{tr}[(V+2\sigma \Sigma_0)^{-1} \Sigma_0]-2\sigma \text{tr}[\{(V+2\sigma \Sigma_0)^{-1} \Sigma_0\}^2]
\}
\text{tr}[V^{-1} \Sigma_0 (V+2\sigma \Sigma_0)^{-1} \Sigma_0]\\
&~~~~~+2 \text{tr}[(V+2\sigma \Sigma_0)^{-1} \Sigma_0] \{
-2\sigma^2 \text{tr}[V^{-1} \Sigma_0 \{(V+2\sigma \Sigma_0)^{-1} \Sigma_0\}^2]\\
&~~~~~+4\sigma^3 \text{tr}[V^{-1} \Sigma_0 \{(V+2\sigma \Sigma_0)^{-1} \Sigma_0\}^3]
\}\\
&=-2\sigma \text{tr}[V^{-1} \Sigma_0] \{\text{tr}[(V+2\sigma \Sigma_0)^{-1} \Sigma_0]\}^2
+\sigma \{\text{tr}[(V+2\sigma \Sigma_0)^{-1} \Sigma_0]\}^3\\
&~~~~~
+4\sigma^2 \text{tr}[V^{-1} \Sigma_0] \text{tr}[(V+2\sigma \Sigma_0)^{-1} \Sigma_0] \text{tr}[\{(V+2\sigma \Sigma_0)^{-1} \Sigma_0\}^2]\\
&~~~~~
-2\sigma^2 \{\text{tr}[(V+2\sigma \Sigma_0)^{-1} \Sigma_0]\}^2 \text{tr}[\{(V+2\sigma \Sigma_0)^{-1} \Sigma_0\}^2]
-\text{tr}[V^{-1} \Sigma_0] \text{tr}[(V+2\sigma \Sigma_0)^{-1} \Sigma_0] \\
&~~~~~
+\{\text{tr}[(V+2\sigma \Sigma_0)^{-1} \Sigma_0]\}^2
+2\sigma\text{tr}[V^{-1} \Sigma_0] \text{tr}[\{(V+2\sigma \Sigma_0)^{-1} \Sigma_0\}^2]\\
&~~~~~
-2\sigma \text{tr}[(V+2\sigma \Sigma_0)^{-1} \Sigma_0] \text{tr}[\{(V+2\sigma \Sigma_0)^{-1} \Sigma_0\}^2]\\
&~~~~~-4 \sigma^2\text{tr}[(V+2\sigma \Sigma_0)^{-1} \Sigma_0]  \text{tr}[\{(V+2\sigma \Sigma_0)^{-1} \Sigma_0\}^3].
\end{align*}
The terms related to $Q_{4}$ involve more terms that should be arranged clearly.
We have 
\begin{align*}
&Q_4(S_2R^{-1},S_2R^{-1},S_2'R^{-1},S_2'R^{-1})\\
&=\{\text{tr}[S_2 R^{-1}]\}^2 \{\text{tr}[S'_2 R^{-1}]\}^2 
+16\text{tr}[S_2 R^{-1}] \text{tr}[S_2 R^{-1} \{S'_2 R^{-1}\}^2] 
+16 \text{tr}[S'_2 R^{-1}] \text{tr}[S'_2 R^{-1} \{S_2 R^{-1}\}^2]\\
&~~~~~+4\text{tr}[\{S_2 R^{-1} \}^2] \text{tr}[\{S'_2 R^{-1}\}^2] 
+8 \{\text{tr}[S_2 R^{-1} S'_2 R^{-1}] \}^2
+2\{\text{tr}[S_2 R^{-1}] \}^2 \text{tr}[\{S'_2 R^{-1}\}^2]\\
&~~~~~
+8 \text{tr}[S_2 R^{-1}] \text{tr}[S'_2 R^{-1}] \text{tr}[S_2 R^{-1} S'_2 R^{-1}]
+2\{\text{tr}[S'_2 R^{-1}] \}^2 \text{tr}[\{S_2 R^{-1}\}^2]\\
&~~~~~
+32 \text{tr}[\{S_2 R^{-1}\}^2 \{S'_2 R^{-1}\}^2]
+16 \text{tr}[\{S_2 R^{-1} S'_2 R^{-1} \}^2]\\
&=4\text{tr}[\{S_2 R^{-1} \}^2] \text{tr}[\{S'_2 R^{-1}\}^2] 
+8 \{\text{tr}[S_2 R^{-1} S'_2 R^{-1}] \}^2
+32 \text{tr}[\{S_2 R^{-1}\}^2 \{S'_2 R^{-1}\}^2]\\
&~~~~~
+16 \text{tr}[\{S_2 R^{-1} S'_2 R^{-1} \}^2],
\end{align*}
where
\begin{align*}
&4\text{tr}[\{S_2R^{-1}\}^2]\text{tr}[\{S_2'R^{-1}\}^2]\\
&=\{\text{tr}[\{(V+2\sigma \Sigma_0)^{-1} \Sigma_0\}^2]-2\sigma \text{tr}[\{(V+2\sigma \Sigma_0)^{-1} \Sigma_0\}^3]\}
\left\{
\frac{1}{2\sigma} \text{tr}[V^{-1} \Sigma_0]-\frac{1}{2\sigma} \text{tr}[(V+2\sigma \Sigma_0)^{-1} \Sigma_0]
\right\}\\
&=\frac{1}{2\sigma} \text{tr}[V^{-1} \Sigma_0] \text{tr}[\{(V+2\sigma \Sigma_0)^{-1} \Sigma_0\}^2]
-\frac{1}{2\sigma} \text{tr}[(V+2\sigma \Sigma_0)^{-1} \Sigma_0] \text{tr}[\{(V+2\sigma \Sigma_0)^{-1} \Sigma_0\}^2]\\
&~~~~~- \text{tr}[V^{-1} \Sigma_0] \text{tr}[\{(V+2\sigma \Sigma_0)^{-1} \Sigma_0\}^3] 
+ \text{tr}[(V+2\sigma \Sigma_0)^{-1} \Sigma_0] \text{tr}[\{(V+2\sigma \Sigma_0)^{-1} \Sigma_0\}^3],\\
&8\{\text{tr}[S_2R^{-1}S_2'R^{-1}]\}^2\\
&=2\left\{ \frac{1}{2}\text{tr}[V^{-1} \Sigma_0 (V+2\sigma \Sigma_0)^{-1} \Sigma_0]
-\sigma \text{tr}[V^{-1} \Sigma_0 \{(V+2\sigma \Sigma_0)^{-1} \Sigma_0\}^2]
+\frac{1}{2}\text{tr}[\{(V+2\sigma \Sigma_0)^{-1} \Sigma_0\}^2]\right\}^2\\
&=2\{ \text{tr}[ \{(V+2\sigma \Sigma_0)^{-1} \Sigma_0\}^2]
\}^2\\
&32 \text{tr}[\{S_2R^{-1}\}^2 \{S_2'R^{-1}\}^2]\\ 
&=2\text{tr}[V^{-1} \Sigma_0 \{(V+2\sigma \Sigma_0)^{-1} \Sigma_0\}^3]
-4\sigma \text{tr}[V^{-1}\Sigma_0\{(V+2\sigma \Sigma_0)^{-1} \Sigma_0\}^4]\\
&~~~~~
+2\text{tr}[V^{-1} \Sigma_0\{(V+2\sigma \Sigma_0)^{-1} \Sigma_0\}^3] 
-4\sigma \text{tr}[V^{-1} \Sigma_0\{(V+2\sigma \Sigma_0)^{-1} \Sigma_0\}^4]\\
&=4\text{tr}[V^{-1} \Sigma_0\{(V+2\sigma \Sigma_0)^{-1} \Sigma_0\}^3] 
-8\sigma \text{tr}[V^{-1} \Sigma_0\{(V+2\sigma \Sigma_0)^{-1} \Sigma_0\}^4]\\
&=4 \text{tr}[\{(V+2\sigma \Sigma_0)^{-1}\Sigma_0\}^4]
\end{align*}
and
\begin{align*}
&16\text{tr}[\{S_2R^{-1} S_2'R^{-1}\}^2]\\
&=\text{tr}[\{V^{-1} \Sigma_0\}^2 \{(V+2\sigma \Sigma_0)^{-1} \Sigma_0\}^2]
+4\sigma^2 \text{tr}[\{V^{-1} \Sigma_0\}^2 \{(V+2\sigma \Sigma_0)^{-1} \Sigma_0\}^4]\\
&~~~~~
-4\sigma \text{tr}[\{V^{-1} \Sigma_0(V+2\sigma \Sigma_0)^{-1} \Sigma_0\}^2 \{(V+2\sigma \Sigma_0)^{-1} \Sigma_0\}]
+\text{tr}[\{(V+2\sigma \Sigma_0)^{-1} \Sigma_0\}^4]\\
&=-\text{tr}[\{V^{-1} \Sigma_0\}^2\{(V+2\sigma \Sigma_0)^{-1} \Sigma_0\}^2] 
+2\sigma \text{tr}[\{V^{-1} \Sigma_0\}^2 \{(V+2\sigma \Sigma_0)^{-1} \Sigma_0\}^3]\\
&~~~~~
+2 \text{tr}[V^{-1} \Sigma_0 \{(V+2\sigma \Sigma_0)^{-1} \Sigma_0\}^3]
-2\sigma \text{tr}[V^{-1} \Sigma_0 \{(V+2\sigma \Sigma_0)^{-1} \Sigma_0\}^4]\\
&~~~~~
+\text{tr}[\{(V+2\sigma \Sigma_0)^{-1} \Sigma_0\}^4]\\
&= 
2\text{tr}[\{(V+2\sigma \Sigma_0)^{-1} \Sigma_0\}^4], 	
\end{align*}
hence we obtain
\begin{align*}
&Q_4(S_2R^{-1},S_2R^{-1},S_2'R^{-1},S_2'R^{-1})\\
&=\frac{1}{2\sigma} \text{tr}[V^{-1}\Sigma_0] \text{tr}[\{(V+2\sigma \Sigma_0)^{-1} \Sigma_0\}^2]
-\frac{1}{2\sigma} \text{tr}[(V+2\sigma \Sigma_0)^{-1} \Sigma_0] \text{tr}[\{(V+2\sigma \Sigma_0)^{-1} \Sigma_0\}^2]\\
&~~~~~
-\text{tr}[V^{-1}\Sigma_0] \text{tr}[\{(V+2\sigma \Sigma_0)^{-1} \Sigma_0\}^3]
+\text{tr}[(V+2\sigma \Sigma_0)^{-1} \Sigma_0]\text{tr}[\{(V+2\sigma \Sigma_0)^{-1} \Sigma_0\}^3]\\
&~~~~~
+2\{\text{tr}[\{(V+2\sigma \Sigma_0)^{-1} \Sigma_0\}^2]\}^2
+6\text{tr}[\{(V+2\sigma \Sigma_0)^{-1} \Sigma_0\}^4].
\end{align*}
On the other side, we confirm that
\begin{align*}
&Q_4(S_2R^{-1},S_2R^{-1},S_3'R^{-1},S_4'R^{-1})\\
&=
\{\text{tr}[S_2 R^{-1}] \}^2 \text{tr}[S'_3 R^{-1}] \text{tr}[S'_4 R^{-1}]
+16\text{tr}[S_2 R^{-1}] \text{tr}[S_2 R^{-1} S'_3 R^{-1} S'_4 R^{-1}] \\
&~~~~~
+8 \text{tr}[S'_3 R^{-1} ] \text{tr}[\{S_2 R^{-1}\}^2 S'_4 R^{-1}]
+8 \text{tr}[S'_4 R^{-1}]  \text{tr}[\{S_2 R^{-1}\}^2 S'_3 R^{-1}]\\
&~~~~~
+4\text{tr}[\{S_2 R^{-1}\}^2] \text{tr}[S'_3 R^{-1}S'_4 R^{-1}] 
+8\text{tr}[S_2 R^{-1} S'_3 R^{-1}] \text{tr}[S_2 R^{-1} S'_4 R^{-1}]\\
&~~~~~
+2 \{\text{tr}[S_2 R^{-1}]\}^2 \text{tr}[S'_3 R^{-1} S'_4 R^{-1}] 
+4 \text{tr}[S_2 R^{-1} ] \text{tr}[S'_3 R^{-1} ] \text{tr}[S_2 R^{-1} S'_4 R^{-1}] \\
&~~~~~+4 \text{tr}[S_2 R^{-1}] \text{tr}[S'_4 R^{-1}] 
\text{tr}[S_2 R^{-1} S'_3 R^{-1}]
+2\text{tr}[S'_3 R^{-1}] \text{tr}[S'_4 R^{-1}] \text{tr}[\{S_2 R^{-1}\}^2]\\
&~~~~~+32\text{tr}[\{S_2 R^{-1}\}^2 S'_3 R^{-1} S'_4 R^{-1}]
+16\text{tr}[S_2 R^{-1} S'_3 R^{-1} S_2 R^{-1} S'_4 R^{-1}]\\
&=
8 \text{tr}[S'_3 R^{-1} ] \text{tr}[\{S_2 R^{-1}\}^2 S'_4 R^{-1}]
+8 \text{tr}[S'_4 R^{-1}]  \text{tr}[\{S_2 R^{-1}\}^2 S'_3 R^{-1}]\\
&~~~~~
+2\text{tr}[S'_3 R^{-1}] \text{tr}[S'_4 R^{-1}] \text{tr}[\{S_2 R^{-1}\}^2]
+16\text{tr}[S_2 R^{-1} S'_3 R^{-1} S_2 R^{-1} S'_4 R^{-1}],
\end{align*}
where
\begin{align*}
&8\text{tr}[S_3'R^{-1}] \text{tr}[\{S_2R^{-1}\}^2 S_4'R^{-1}]\\
&=2\{\text{tr}[(V+2\sigma \Sigma_0)^{-1} \Sigma_0]-2\sigma \text{tr}[\{(V+2\sigma \Sigma_0)^{-1} \Sigma_0\}^2]\}
\text{tr}[V^{-1}\Sigma_0 \{(V+2\sigma \Sigma_0)^{-1} \Sigma_0\}^2]\\
&=\frac{1}{\sigma}\{\text{tr}[(V+2\sigma \Sigma_0)^{-1} \Sigma_0]-2\sigma \text{tr}[\{(V+2\sigma \Sigma_0)^{-1} \Sigma_0\}^2]\}\\
&~~~~~\times
\{
\text{tr}[V^{-1}\Sigma_0(V+2\sigma \Sigma_0)^{-1} \Sigma_0]
-\text{tr}[\{(V+2\sigma \Sigma_0)^{-1} \Sigma_0\}^2]
\}
\\
&=\frac{1}{\sigma}\{\text{tr}[(V+2\sigma \Sigma_0)^{-1} \Sigma_0]-2\sigma \text{tr}[\{(V+2\sigma \Sigma_0)^{-1} \Sigma_0\}^2]\}\\
&~~~~~\times 
\left\{
\frac{1}{2\sigma}\text{tr}[V^{-1}\Sigma_0]
-\frac{1}{2\sigma} \text{tr}[(V+2\sigma \Sigma_0)^{-1} \Sigma_0]
-\text{tr}[\{(V+2\sigma \Sigma_0)^{-1} \Sigma_0\}^2]
\right\}\\
&=\frac{1}{2\sigma^2} \text{tr}[V^{-1}\Sigma_0] \text{tr}[(V+2\sigma \Sigma_0)^{-1} \Sigma_0]
-\frac{1}{2\sigma^2} \{\text{tr}[(V+2\sigma \Sigma_0)^{-1} \Sigma_0]\}^2\\
&~~~~~
-\frac{1}{\sigma} \text{tr}[V^{-1}\Sigma_0] \text{tr}[\{(V+2\sigma \Sigma_0)^{-1} \Sigma_0\}^2]
+2\{\text{tr}[\{(V+2\sigma \Sigma_0)^{-1} \Sigma_0\}^2]\}^2,\\
&8\text{tr}[S_4'R^{-1}] \text{tr}[\{S_2R^{-1}\}^2S_3'R^{-1}]\\
&=2 \text{tr}[(V+2\sigma \Sigma_0)^{-1} \Sigma_0]  \{
\text{tr}[V^{-1}\Sigma_0 \{(V+2\sigma \Sigma_0)^{-1} \Sigma_0\}^2]-2\sigma \text{tr}[V^{-1}\Sigma_0\{(V+2\sigma \Sigma_0)^{-1} \Sigma_0\}^3]
\}\\
&=2\text{tr}[(V+2\sigma \Sigma_0)^{-1} \Sigma_0] \text{tr}[\{(V+2\sigma \Sigma_0)^{-1} \Sigma_0\}^3],\\
&2\text{tr}[S_3'R^{-1}]\text{tr}[S_4'R^{-1}] \text{tr}[\{S_2R^{-1}\}^2]\\
&=\{
\{\text{tr}[(V+2\sigma \Sigma_0)^{-1} \Sigma_0]\}^2
-2\sigma \text{tr}[(V+2\sigma \Sigma_0)^{-1} \Sigma_0] \text{tr}[\{(V+2\sigma \Sigma_0)^{-1} \Sigma_0\}^2]
\}\\
&~~~~~\times
\left\{
\frac{1}{2\sigma} \text{tr}[V^{-1} \Sigma_0]-\frac{1}{2\sigma} \text{tr}[(V+2\sigma \Sigma_0)^{-1} \Sigma_0]
\right\}\\
&=\frac{1}{2\sigma} \text{tr}[V^{-1}\Sigma_0] \{\text{tr}[(V+2\sigma \Sigma_0)^{-1} \Sigma_0]\}^2
-\frac{1}{2\sigma} \{\text{tr}[(V+2\sigma \Sigma_0)^{-1} \Sigma_0]\}^3\\
&~~~~~
-\text{tr}[V^{-1}\Sigma_0] \text{tr}[(V+2\sigma \Sigma_0)^{-1} \Sigma_0] \text{tr}[\{(V+2\sigma \Sigma_0)^{-1} \Sigma_0\}^2]\\
&~~~~~
+\{\text{tr}[(V+2\sigma \Sigma_0)^{-1} \Sigma_0]\}^2 \text{tr}[\{(V+2\sigma \Sigma_0)^{-1} \Sigma_0\}^2],
\end{align*}
and
\begin{align*}
&16\text{tr}[S_2R^{-1}S_3'R^{-1}S_2R^{-1}S_4'R^{-1}]\\
&=4\text{tr}[V^{-1}\Sigma_0\{(V+2\sigma \Sigma_0)^{-1} \Sigma_0\}^3] 
-8\sigma \text{tr}[V^{-1}\Sigma_0 \{(V+2\sigma \Sigma_0)^{-1} \Sigma_0\}^4]\\
&=4\text{tr}[\{(V+2\sigma \Sigma_0)^{-1} \Sigma_0\}^4],
\end{align*}
so we have
\begin{align*}
&Q_4(S_2R^{-1},S_2R^{-1},S_3'R^{-1},S_4'R^{-1})\\
&=\frac{1}{2\sigma^2} \text{tr}[V^{-1}\Sigma_0] \text{tr}[(V+2\sigma \Sigma_0)^{-1} \Sigma_0]
-\frac{1}{2\sigma^2} \{\text{tr}[(V+2\sigma \Sigma_0)^{-1} \Sigma_0]\}^2
-\frac{1}{\sigma} \text{tr}[V^{-1}\Sigma_0] \text{tr}[\{(V+2\sigma \Sigma_0)^{-1} \Sigma_0\}^2]\\
&~~~~~
+\frac{1}{2\sigma} \text{tr}[V^{-1}\Sigma_0] \{\text{tr}[(V+2\sigma \Sigma_0)^{-1} \Sigma_0]\}^2
-\frac{1}{2\sigma} \{\text{tr}[(V+2\sigma \Sigma_0)^{-1} \Sigma_0]\}^3
+2\{\text{tr}[\{(V+2\sigma \Sigma_0)^{-1} \Sigma_0\}^2]\}^2\\
&~~~~~
+2\text{tr}[(V+2\sigma \Sigma_0)^{-1} \Sigma_0] \text{tr}[\{(V+2\sigma \Sigma_0)^{-1} \Sigma_0\}^3]
-\text{tr}[V^{-1}\Sigma_0] \text{tr}[(V+2\sigma \Sigma_0)^{-1} \Sigma_0] \text{tr}[\{(V+2\sigma \Sigma_0)^{-1} \Sigma_0\}^2] \\
&~~~~~
+\{\text{tr}[(V+2\sigma \Sigma_0)^{-1} \Sigma_0]\}^2 \text{tr}[\{(V+2\sigma \Sigma_0)^{-1} \Sigma_0\}^2]
+4\text{tr}[\{(V+2\sigma \Sigma_0)^{-1} \Sigma_0\}^4].
\end{align*}
All expressions obtained above furnish to reach 
\begin{align*}
&\mathcal{F}_{24}=|V+2\sigma \Sigma_0|^{-1} \Big\{
1+8\sigma^2 \text{tr}[\{(V+2\sigma \Sigma_0)^{-1} \Sigma_0\}^2]
+12\sigma^4 \{\text{tr}[\{(V+2\sigma \Sigma_0)^{-1} \Sigma_0\}^2]\}^2\\
&~~~~~~~~~~~~~~~~~~~~~~~~~~~~~+24\sigma^4 \text{tr}[\{(V+2\sigma \Sigma_0)^{-1} \Sigma_0\}^4]
\Big\},
\end{align*}
which is nothing other than $\mathcal{F}_{14}$ in (\ref{equation_F14}), which completes the proof of (\ref{J14=J24}).

Other equalities in (\ref{J22=J12}), (\ref{J44=J14}) and (\ref{J44=J23}) can be confirmed in the same way.
These equalities finally gives (\ref{VZ_SIMPLE}). \qed

\def\thesection{}
\addcontentsline{toc}{section}{\refname}
\bibliographystyle{econ}

\end{document}